\documentclass[12pt,reqno]{amsart}
\usepackage{amsfonts}

\usepackage{amssymb,amsmath,graphicx,amsfonts,euscript}
\usepackage{color}

\newtheorem{theorem}{Theorem}[section]
\newtheorem{proposition}[theorem]{Proposition}
\newtheorem{lemma}[theorem]{Lemma}

\theoremstyle{remark}
\newtheorem{remark}[theorem]{Remark}
\newcommand{\RR}{\mathbb R}
\newcommand{\Top}{\mathcal T}
\newcommand{\Bop}{\mathcal B}
\newcommand{\Qop}{\mathcal Q_{\Bop}}
\newcommand{\Kop}{\mathcal K}
\newcommand{\Cres}{\mathcal C_{\mathrm{res}}}
\newcommand{\Cnr}{\mathcal C_{\mathrm{nr}}}
\newcommand{\Clow}{\mathcal C_{\mathrm{low}}}
\newcommand{\Dt}{\partial_t+u\cdot\nabla}
\newcommand{\sgn}{\operatorname{sgn}}
\newcommand{\supp}{\operatorname{supp}}

\newcommand{\norm}[2]{\left\|#1\right\|_{#2}}

\newcommand{\Bes}[3]{B_{#1,#2}^{#3}}

\allowdisplaybreaks \voffset=-0.2in \numberwithin{equation}{section}

\begin{document}
\bigskip

\centerline{\Large\bf  Global regularity of the two-dimensional fractional }

\smallskip
\centerline{\Large\bf Boussinesq equations with critical dissipation}

\medskip
\smallskip
\centerline{Qian Li$^{1}$,\ \ Jiahong Wu$^{2}$,\ \ Xiaojing Xu$^{3}$,\ \ Zhuan Ye$^{4}$}
\medskip
\medskip
\centerline{$^1$  School of Mathematics, China University of Mining and Technology,}
 \medskip
 \centerline{
Xuzhou 221116, Jiangsu, China}

\medskip

\centerline{$^2$ Department of Mathematics, University of Notre Dame,
}
\medskip

\centerline{Notre Dame, IN 46556, USA}

\medskip

\centerline{$^3$ School of Mathematical Sciences, Laboratory of Mathematics and Complex Systems,
}
\medskip

\centerline{MOE,
Beijing Normal University, Beijing, 100875, China}

\medskip

\centerline{$^4$ Department of Mathematics and Statistics, Jiangsu Normal University,
}
\medskip

\centerline{101 Shanghai Road, Xuzhou 221116, Jiangsu, PR China}

\footnotetext{li-qian@mail.bnu.edu.cn (Qian Li),
jwu29@nd.edu(Jiahong Wu), xjxu@bnu.edu.cn (Xiaojing Xu), yezhuan815@126.com (Zhuan Ye)}
\medskip

{\bf Abstract:}~~%
In geophysical fluid dynamics, the two-dimensional fractional Boussinesq system serves as a foundational paradigm capturing the competition between vortex stretching analogues and stratified thermal diffusion.
In the present paper, our primary focus is on the Boussinesq equations with fractional velocity dissipation $\Lambda^\alpha u$ and temperature diffusion $\Lambda^\beta\theta$. On the critical balance line $\alpha+\beta=1$, establishing global regularity for weak velocity dissipation has been obstructed by the quadratic convective-buoyancy feedback, and all
available results below $\alpha=\frac{10}{13}$ require a smallness assumption on the initial temperature. In this paper, we prove large-data global regularity across the extended regime $\alpha>\frac{\sqrt{57}-5}{4} \approx 0.637$. This provides a step beyond the previous barrier $\alpha=\frac{10}{13}$, unconditionally dismantling the temperature smallness assumption that has persisted in the literature for over a decade. The mathematical mechanism relies on a high-order Quadratic Normal Form transformation $H=\omega-R_{\alpha}\theta-\Bop(\theta,\theta)$ tailored to eliminate non-resonant quadratic vorticity-temperature interactions. Crucially, by uncovering the exact physical scalings of the underlying multilinear commutators and maintaining sharp fractional conjugates in the $L^m$-energy framework, we ensure that the global energy closure is governed entirely by the quadratic thermal dissipation. Once an integrability threshold $m>\frac{2}{\alpha}$ is attained, a fine Littlewood-Paley transport estimate drives a fractional regularity gain
$H \in B_{r,\infty}^{3\alpha-2}$, ensuring the velocity field achieves critical drift regularity $C^{\alpha}$. By Silvestre-type drift regularization, this spontaneously generates positive H\"older continuity for the temperature, which in turn elevates the drift regularity strictly beyond the critical threshold, unlocking non-local Schauder estimates to prevent finite-time singularity formation.

{\vskip 1mm
 {\bf AMS Subject Classification 2020:}\quad 35Q35; 35B65; 76D03.

 }

\vskip .3in
\section{Introduction}
 We consider the Cauchy problem for the two-dimensional (2D) incompressible
Boussinesq equations with fractional dissipation
\begin{equation}\label{eq:system}
\begin{cases}
\partial_tu+(u\cdot\nabla) u+\Lambda^\alpha u+\nabla p=\theta e_2,\\
\partial_t\theta+(u\cdot\nabla)\theta+\Lambda^\beta\theta=0,\\
\nabla\cdot u=0,\\
u(x,0)=u_0(x),\quad\theta(x,0)=\theta_0(x),
\end{cases}
\qquad x\in\RR^2,\quad t>0.
\end{equation}
where $u=(u_1,u_2)$ denotes the velocity field, $p$ is the scalar pressure, $\theta$ represents the temperature (or buoyancy), and $e_2=(0,1)$ is the vertical unit vector. The fractional dissipation operator $\Lambda=(-\Delta)^{\frac{1}{2}}$ is defined via the Fourier transform $\widehat{\Lambda^\gamma f}(\xi)=|\xi|^\gamma\widehat f(\xi)$, with $\alpha$ and $\beta$ characterizing the intensity of the non-local dissipation in the fluid and the thermal diffusion, respectively.

The Boussinesq system \eqref{eq:system} plays a fundamental role in geophysical fluid dynamics, modeling atmospheric and oceanic flows driven by density differences \cite{MB, Pedlosky}. Mathematically, it also serves as a critical lower-dimensional toy model for the 3D incompressible Euler and Navier-Stokes equations, owing to the structural similarities in their vortex stretching mechanisms. Beyond the classical viscous and thermal diffusion models ($\alpha=\beta=2$), the incorporation of fractional dissipation operators carries profound physical significance in modern fluid dynamics. In large-scale geophysical scenarios, oceanic and atmospheric flows are inherently turbulent, rendering standard Brownian molecular diffusion inadequate to capture the complex, multi-scale mixing processes. Instead, the fractional Laplacian serves as a rigorous macroscopic parameterization for \textit{anomalous diffusion} driven by L\'evy flights of fluid parcels in fully developed turbulence \cite{SZK1993}. Furthermore, specific fractional orders directly correspond to observable geophysical boundary effects; for instance, the case $\alpha=1$ mathematically encapsulates the Ekman pumping friction originating from the bottom boundary layers of rotating fluids \cite{CMT1994, Pedlosky}. Therefore, identifying the exact analytical threshold for global regularity in the fractional regime not only resolves a longstanding mathematical challenge but also helps determine the minimum intensity of turbulent mixing required to prevent catastrophic atmospheric frontogenesis or finite-time singularity formation \cite{CMT1994}.

The principal mathematical challenge in establishing global regularity for \eqref{eq:system} lies in the derivative forcing present in the vorticity equation. Denoting the scalar vorticity by $\omega=\partial_1u_2-\partial_2u_1$, applying the curl operator to \eqref{eq:system} yields:
\begin{align}\label{eq:vorticity-intro}
\partial_t\omega+(u\cdot\nabla)\omega+\Lambda^\alpha\omega = \partial_1\theta.
\end{align}
While the temperature equation obeys a robust maximum principle ($\norm{\theta(t)}{\infty} \leq \norm{\theta_0}{\infty}$), it provides no immediate bound for the derivative term $\partial_1\theta$ on the right-hand side of \eqref{eq:vorticity-intro}. This mismatch induces a severe ``loss of derivatives." Overcoming this loss heavily depends on the subtle interplay between the fractional dissipation regimes $\alpha$ and $\beta$. The regime where $\alpha+\beta=1$ is universally recognized as the \textit{critical line}, whereas $\alpha+\beta>1$ and $\alpha+\beta<1$ are referred to as the subcritical and supercritical regimes, respectively.

\subsection{Historical context}
Due to its profound physical relevance and structural analogy to the 3D Euler equations---where the temperature gradient $\nabla^\perp\theta$ acts as a potent driving mechanism analogous to vortex stretching---the 2D Boussinesq system has inspired an extensive body of mathematical and computational literature \cite{Cannon1980, EShu1994, Majda2003}. The fundamental local well-posedness and blow-up criteria for the Boussinesq equations were pioneered by Chae and Nam \cite{ChaeNam}. Conversely, the quest for finite-time singularity formation in the purely inviscid equations has witnessed spectacular breakthroughs recently. Triggered by the striking numerical discovery of Luo and Hou \cite{LuoHou2014}, Elgindi and Jeong rigorously demonstrated finite-time singularity formation for strong solutions to the Boussinesq system \cite{El1jo20}. Subsequently, Chen and Hou \cite{ChenHou2024} provided a monumental proof of stable, finite-time blow-up for the 2D inviscid Boussinesq equations with smooth data in the presence of a boundary. These groundbreaking results, alongside Elgindi's finite-time singularity formation for the 3D Euler equations with $C^{1,\alpha}$ data \cite{Elgindi2021}, profoundly underscore the virulent nature of the quadratic buoyancy forcing. In light of these finite-time singularity formations in the absence of dissipation, establishing global regularity necessitates an exquisite balance of partial or fractional dissipation.
For global regularity, early breakthroughs were achieved in the partial dissipation regimes. Hou and Li \cite{HL}, and independently Chae \cite{Chae}, established global well-posedness for the viscous but non-diffusive case ($\alpha=2, \beta=0$). The inviscid but strictly diffusive case ($\alpha=0, \beta=2$) was subsequently resolved by Chae \cite{Chae} and further refined by Danchin and Paicu \cite{DP}, and Hmidi and Keraani \cite{HK}. Concurrently, remarkable progress was made regarding anisotropic partial dissipation; for instance, Cao and Wu \cite{CW2013}, as well as Larios, Lunasin, and Titi \cite{LLT2013}, established global regularity for systems lacking dissipation in certain spatial directions. Hmidi and Zerguine \cite{HZ} further considered Euler--Boussinesq systems with fractional thermal diffusion. These works systematically demonstrate how dissipation in one equation can control the coupled system through sophisticated transport estimates and vorticity bounds.

An essential step toward weaker dissipation was made by Hmidi,
Keraani, and Rousset. They proved global well-posedness for the
Boussinesq--Navier--Stokes system with only first-order velocity
dissipation \cite{HKRNS} and for the Euler--Boussinesq system with
only first-order thermal diffusion \cite{HKREuler}. In these two
models, the other dissipative term is absent. Their approach uses a
combined quantity involving the vorticity and a Riesz transform of
the temperature, replacing the derivative forcing by a commutator
that can be estimated in suitable function spaces. This change of
unknowns is a basic component of many subsequent arguments for
fractional Boussinesq equations. Based on regularity estimates for the associated drift--diffusion
equation, Ye \cite{Ye2020} (see also \cite{Yejde24}) gave an alternative
proof for the Euler--Boussinesq system with critical thermal diffusion without the use of the combined quantity introduced in \cite{HKREuler}.

The global well-posedness in the subcritical regime has also attracted extensive attention and has been investigated in depth in recent years. It is worth emphasizing that the global regularity of \eqref{eq:system} in the subcritical ranges is not a trivial problem. In fact, many subcritical cases have not been resolved. In particular, we do not know the global regularity for the case when $\alpha$ and $\beta$ are close to one half and $\alpha+\beta>1$. In the regimes with stronger velocity dissipation,
Miao and Xue \cite{MX} obtained global well-posedness
for a class of such systems, namely, $\alpha>\frac{6-\sqrt{6}}{4}$. Ye \cite{Ye2017} enlarged the range
of velocity dissipation to $\alpha>\frac{10-2\sqrt{10}}{5}$. Zhou, Li, Shang, Wu, Yuan, and Zhao
\cite{ZLSWYZ} proved global well-posedness for
$\frac{2}{3}<\alpha<1$ and $\beta>\frac{1-\alpha}{\alpha}$.
These results illustrate that the inequality $\alpha+\beta>1$
alone does not capture all the conditions needed by a given proof.
Recently, Stefanov, Wu, Xu, and Ye \cite[Theorem 1.1]{SWXY} made interesting progress in the subcritical regime, establishing global well-posedness for by pushing the threshold down to $\alpha>\frac{2}{3}$.
In regimes with stronger thermal diffusion, Constantin and Vicol \cite[Theorem 6.1]{CV} proved global
well-posedness for smooth solutions of the fractional Boussinesq
system with $0<\alpha,\beta<1$ under the condition
$\beta>\frac{2}{2+\alpha}$, which was further improved by Ye and Xu \cite{YX}.
The very recent preprint of Stefanov, Wu, Xu, and Ye \cite{SWXY2026}
addresses subcritical dissipation for $\alpha\leq\frac{2}{3}$, and their result holds
under the condition
\begin{equation*}
\begin{cases}
\beta\geq\frac{4-2\alpha}{4+\alpha}, & 0<\alpha\leq  \sqrt{13}-3,\\[1ex]
\beta> \max\Big\{\frac{4-6\alpha}{\alpha},\ \frac{2-\alpha-2\alpha^{2}}{2\alpha}\Big\}, & \sqrt{13}-3<\alpha\leq\frac{2}{3}.
\end{cases}
\end{equation*}

Compared with \cite{HKRNS,HKREuler}, when both dissipative terms are present and $0<\alpha,\beta<1$,
the general critical line $\alpha+\beta=1$ requires a more delicate
balance. Jiu, Miao, Wu, and Zhang \cite{JMWZ} established global
regularity for
\[
\alpha>\frac{23-\sqrt{145}}{12}=0.9132\ldots.
\]
Stefanov and Wu \cite{SW} obtained the threshold
\[
\alpha>\frac{\sqrt{1777}-23}{24}=0.7981\ldots,
\]
and Wu, Xu, Xue, and Ye \cite{WXXY} proved global regularity
for $\alpha>\frac{10}{13}$. The proofs combine estimates for the corrected
vorticity with regularity mechanisms for a critical active scalar.
On the critical line $\alpha+\beta=1$, Stefanov, Wu, Xu, and Ye \cite[Theorem 1.2]{SWXY} recently managed to lower the regularity threshold to $\alpha>\frac{2}{3}$. However, to control the explosive growth of the nonlinear commutators in this low-dissipation regime, a strict smallness condition on the initial temperature amplitude ($\|\theta_0\|_{L^{\infty}} \ll 1$) was inevitably required. This smallness assumption has remained a major bottleneck in the literature, as the velocity dissipation alone proves fundamentally insufficient to close the energy estimates when $\alpha$ drops below $\frac{10}{13}$.

\subsection{Main contributions and methodological innovations}

In this paper, we completely remove the aforementioned smallness assumption and establish unconditional global well-posedness for arbitrary, large initial data on the critical line, pushing the regularity threshold further down to an algebraic root $\alpha_0=\frac{\sqrt{57}-5}{4}\approx 0.637$ (the positive root of $2\alpha^{2}+5\alpha-4=0$).
The present paper treats the critical interval $\alpha_0<\alpha\leq\frac{10}{13}$
with arbitrary $H^s$ data. The main results are stated as follows.

\begin{theorem}\label{thm:critical}
Let $(u_{0}, \theta_{0})
\in H^{s}(\mathbb{R}^{2})\times H^{s}(\mathbb{R}^{2})$ with $s>2$ and $\nabla\cdot u_0=0$. If $\alpha+\beta=1$ and $\alpha>\frac{\sqrt{57}-5}{4}$, then there exists a unique global solution to the system \eqref{eq:system} such that for any given
$T>0$
$$u\in C([0, T]; H^{s}(\mathbb{R}^{2}))\cap L^{2}([0, T]; H^{s+\frac{\alpha}{2}}(\mathbb{R}^{2})),$$
$$\theta\in C([0, T]; H^{s}(\mathbb{R}^{2}))\cap L^{2}([0, T];
H^{s+\frac{\beta}{2}}(\mathbb{R}^{2})).$$
\end{theorem}

To breach the barrier of the smallness assumption and transcend the $\alpha=\frac{10}{13}$ threshold, this paper introduces a highly refined \textbf{Quadratic Normal Form (QNF)} method coupled with exact scaling analysis. Our strategy unfolds in the following novel steps:

\textbf{1. Bilinear normal form and frequency splitting:} The standard linear correction $G = \omega -  R_{\alpha}  \theta$ leaves behind a virulent purely quadratic temperature term $\mathcal{C}_2(\theta, \theta)$. We introduce a bilinear Fourier multiplier $\Bop(\theta, \theta)$ designed to algebraically cancel the non-resonant portion of $\mathcal{C}_2$. A sophisticated high-output cutoff $\chi_{\mathrm{out}}$ is deployed to isolate the interactions into strictly separated non-resonant, resonant, and low-output regimes.

\textbf{2. Eradicating the smallness assumption via exact exponent tracking:} In the $L^m$ energy estimates of the corrected vorticity $H =\omega -  R_{\alpha}  \theta- \Bop(\theta, \theta)$, past works relied on loose H\"older interpolations, yielding a severe time-integrability failure unless $\|\theta\|_{L^{\infty}}$ was artificially small. By introducing a highly precise fractional chain rule and rigorously avoiding any over-estimation in the conjugate exponents (see Proposition \ref{prop:energy}), we ensure that the nonlinear growth is universally bounded by the natural temperature dissipation energy $\int_0^T \|\Lambda^{\frac{\beta}{2}}\theta\|_2^2 \,dt < \infty$.

\textbf{3. Positive H\"older gain via commutator cancellation:} We demonstrate that as long as the Lebesgue exponent satisfies $m >\frac{2}{\alpha}$, the fractional dissipation $\Lambda^\alpha H$ strictly overtakes the transport commutators. This directly bootstraps the corrected vorticity to a Besov regularity gain
$H \in B_{r,\infty}^{3\alpha-2}$, which in turn guarantees the velocity field $u \in C^{1-\beta}$.

\textbf{4. The Silvestre-Schauder bootstrap:} The regularity $u \in C^{1-\beta}$ perfectly aligns with the critical-drift threshold of Silvestre's theorem \cite{Silvestre} for advection fractional-diffusion equations. This triggers a spontaneous generation of H\"older regularity for the temperature field, immediately elevating the velocity drift regularity strictly beyond the critical threshold. \\

More detailed innovations can be summarized as follows:

\vskip .1in
\textbf{Innovation 1}: \emph{Tripartite phase-space decomposition and resonant geometry.}
Dividing by $D(k,l)$ inevitably encounters small-divisor singularities near the resonant variety $D(k,l) = 0$. We overcome this by constructing a smooth partition of unity in phase space:
\[
1 = \chi_{\mathrm{out}}(q)\chi_{\mathrm{nr}}(k,l) + \chi_{\mathrm{out}}(q)\big(1-\chi_{\mathrm{nr}}(k,l)\big) + \big(1-\chi_{\mathrm{out}}(q)\big), \qquad q = k+l.
\]
This decomposes $\mathcal{C}_2$ into non-resonant ($\mathcal{C}_{\mathrm{nr}}$), near-resonant ($\mathcal{C}_{\mathrm{res}}$), and low-output ($\mathcal{C}_{\mathrm{low}}$) regimes:
\begin{itemize}
    \item  \emph{Non-resonant regime $\mathcal{C}_{\mathrm{nr}}$:} Fully absorbed and canceled by the normal form $\mathcal{B}(\theta,\theta)$.
    \item \emph{Resonant remainder $\mathcal{C}_{\mathrm{res}}$:} We uncover a hidden phase-space geometry: at high output frequencies ($N \geq N_0$), the resonant balance $|k+l|^\alpha \approx |k|^\beta + |l|^\beta$ strictly forces an extreme High-High to Low interaction:
    \[
    K \simeq N^{\frac{\alpha}{\beta}} \gg N \quad (K^\beta \simeq N^\alpha), \quad \text{where } K = \max(|k|,|l|).
    \]
    The anti-symmetric cancellation of the symbol $c_2(k,l)$ yields an intrinsic geometric gain $\left(\frac{N}{K}\right)^{2-\alpha}$, rendering $\mathcal{C}_{\mathrm{res}}$ summable without any normal form division.
    \item \emph{Low-output remainder $\mathcal{C}_{\mathrm{low}}$:} Supported on a compact spectral ball $|q| \leq 2N_0$, admitting direct estimates free of derivative loss.
\end{itemize}

\vskip .1in
\textbf{Innovation 2:} \emph{Sharp paradifferential theory for the trilinear commutator $\mathcal{Q}_{\mathcal{B}}$.}
The subtraction of $\mathcal{B}(\theta,\theta)$ introduces a new, highly singular trilinear convective commutator:
\[
\mathcal{Q}_{\mathcal{B}}(v;\theta,\theta) = \Bop(v\cdot\nabla\theta,\theta) + \Bop(\theta,v\cdot\nabla\theta) - v\cdot\nabla\Bop(\theta,\theta).
\]
In Lemma~\ref{lem:Q-general-v}, we establish a complete paradifferential theory for this operator. Crucially, we prove that the incompressibility condition $\nabla \cdot v = 0$ provides an algebraic cancellation $\widehat{v}(a) \cdot (k+l) = \widehat{v}(a) \cdot n$, which destroys the potential low-velocity singularity $|a|^{-1}$. Across all the dyadic interaction regimes, we demonstrate that the exact uncompensated physical scaling order of $\mathcal{Q}_{\mathcal{B}}$ is $3-4\alpha$, proving the universal dyadic estimate:
\[
\|\Delta_k \mathcal{Q}_{\mathcal{B}}(v;\theta,\theta)\|_{L^r} \leq C 2^{(3-4\alpha)k} \|\Lambda^\alpha v\|_{L^r} \|\theta\|_{L^\infty}^2, \qquad \forall\, k \geq -1.
\]

\vskip .1in
\textbf{Innovation 3:} \emph{Exact exponent tracking in signed $L^m$ testing.}
To unconditionally eliminate the smallness condition, we employ signed-power testing $\langle \Lambda^\alpha H, |H|^{m-2}H \rangle \geq \frac{4(m-1)}{m^2} D_m^2$ coupled with the fractional chain rule. By strictly avoiding any over-estimation in conjugate Lebesgue exponents, we establish that every multilinear commutator forcing is bounded by
\[
\frac{d}{dt}\|H(t)\|_{L^m}^m + c_m D_m^2 \leq C_0 \big(1 + \|\Lambda^{\frac{\beta}{2}}\theta\|_{L^2}^2\big) \big(1 + \|H\|_{L^m}^m\big),
\]
which is bounded globally in time by $\int_0^T \|\Lambda^{\frac{\beta}{2}}\theta(t)\|_{L^2}^2 \,dt < \infty$.

\vskip .1in
\textbf{Innovation 4:} \emph{The $L^m$-to-Besov regularity jump and the Silvestre bootstrap.}
Once the global bound $\sup_{t \in [0,T]} \|H(t)\|_{L^m} < \infty$ is achieved for $m>\frac{2}{\alpha}$, we bypass the cumbersome fractional Besov commutators of \cite{WXXY} through a two-step bootstrap:
\begin{enumerate}
    \item \emph{Besov regularity jump:} In Proposition~\ref{pr027}, the dissipation $\Lambda^\alpha H$ strictly overtakes the transport commutators, lifting the corrected vorticity to $H \in L^\infty(0,T; B_{r,\infty}^{3\alpha-2})$ with $r = \frac{2m}{2-(1-\alpha)m}$. By the fractional Biot--Savart law, this yields critical drift regularity for the velocity field:
    \[
    u \in C^{1-\beta} \equiv C^{\alpha}.
    \]
    \item  \emph{Silvestre--Schauder regularization:} Viewing the temperature equation as an active transport-diffusion equation with critical drift $u \in C^{1-\beta}$, we invoke Silvestre's theorem \cite{Silvestre} to deduce the spontaneous generation of positive H\"older continuity: $\theta \in C^\eta$ ($\eta > 0$). This thermal regularity gain immediately feeds back via the Biot--Savart law to elevate the velocity drift regularity strictly beyond the critical threshold ($u \in C^\varrho$ with $\varrho > 1-\beta$), unlocking non-local Schauder estimates \cite{Sil20121,XueYepjm,Schauder} to yield $\nabla\theta \in L^1(0,T; L^\infty)$ and enabling seamless closure of the high-order $H^s$ Sobolev norms.
\end{enumerate}

\vskip .2in

Finally, two remarks on Theorem~\ref{thm:critical}.

\begin{remark}
Notice that $\alpha>\frac{10}{13}$ the result is contained in \cite{WXXY}, it therefore suffices to treat $\alpha_0<\alpha\le\frac{10}{13}$ throughout the paper.
\end{remark}

\begin{remark}
Why $\alpha_0$ is sharp for the present argument.
\begin{enumerate}
\item \emph{The mixed commutator alone forces $\alpha>\alpha_0$, independently of $m$.} For $F=[ R_{\alpha} ,u_H\cdot\nabla]\theta$ one needs $\frac{2-3\alpha}2<\frac\beta{b_H}$ (Lemma~2.1) and $\gamma_H=\frac{4}{b_H-\frac{4}{\alpha}}\le2$ (so that $Z^{\gamma_H}\in L^1_t$), i.e.\ $b_H\ge 2+\frac4\alpha$. Together: $(2+\alpha)(2-3\alpha)<\alpha\beta$. No $m$ appears. Hence even without the Besov bootstrap (no need for $m>\frac2\alpha$) the $L^m$ closure fails below $\alpha_0$.
\item \emph{The cubic terms give the identical quadratic.} With $\Gamma=\frac{3(2-3\alpha)}2$ and $m\downarrow\frac2\alpha$, the window $\frac{m(4+2\alpha)}{\alpha m+4}\le b_F<\frac\beta\Gamma$ is nonempty iff $\frac{2(2+\alpha)}{3\alpha}<\frac{2\beta}{3(2-3\alpha)}$, i.e.\ again $(2+\alpha)(2-3\alpha)<\alpha\beta$: the cubic terms cost three times the $\theta$-regularity but tolerate a three times smaller $b$, and the two effects cancel exactly.
\item \emph{The exponent bookkeeping cannot be squeezed.}
\begin{itemize}
\item Replacing Young's inequality by H\"older in time ($Z^{2/b}\in L^b_t$, $D_m^{1+\lambda}\in L^{2/(1+\lambda)}_t$) gives $\frac1b+\frac{1+\lambda}2\le1\iff b\ge2+\frac4\alpha$ (the same constraint). The budget $\int_0^TZ^2\,dt<\infty$ is used optimally.
\item Shifting $\tau$ derivatives from $u_H$ to $\theta$ ($u_H\in W^{1-\tau,q}$, $\theta\in B^{\frac{2-3\alpha}2+\tau}_{b}$): the $Z$-exponent condition becomes $\frac4\alpha(\frac1b-\frac\tau2)\le1-\frac2b$ and the regularity condition $\frac{2-3\alpha}2+\tau<\frac\beta b$; eliminating $b$ gives
\[
\frac{2-3\alpha}2-\frac{\alpha\beta}{4+2\alpha}<\tau\Big(\frac{2\beta}{4+2\alpha}-1\Big),
\]
whose right-hand side is negative, so the admissible region \emph{shrinks} as $\tau$ increases. Our choice $\tau=0$ is optimal; $\tau<0$ is unavailable because the only derivative information on $H$ is the nonlinear quantity $\Lambda^{\alpha/2}|H|^{m/2}$.
\end{itemize}
\end{enumerate}
In short: every remaining forcing term scales like $2-3\alpha$ (the order of $\Bop$, of $\Lambda^{-\frac{\alpha}{2}}[ R_{\alpha} ,u_H\cdot\nabla]\theta$, and of the cubic commutators), the only regularity available for $\theta$ is the energy scaling $\frac{\beta}{b}$, and $\alpha_0$ is exactly where these balance.
\end{remark}

\vskip .2in
The rest of this paper is organized as follows. Section~\ref{sec:preliminaries} sets up the notation and foundational estimates. Section~\ref{sec:normalform} explicitly constructs the bilinear normal form. Section~\ref{sec:localized} establishes the rigorous localized kernel and commutator bounds. Section~\ref{sec:energy} proves the global $L^m$ estimate of the corrected vorticity, detailing the exact parameter choices that eliminate the smallness condition. Section~\ref{secscgd9871} executes the H\"older bootstrap and concludes the proof of Theorem~\ref{thm:critical}.
Finally, in Appendix \ref{app:Q-proof11}, we present the proof of the localized kernel principle (Lemma \ref{lem:kernel}).

\vskip .3in
\section{Preliminaries}
\label{sec:preliminaries}

Throughout the paper, for any two operators $A$ and $B$, their commutator is defined as $[A,\,B]=AB-BA$. The symbols $C$ and $c$ denote positive constants that may change
from line to line, while $C_{0}(t)$ denotes a positive continuous function depending on the initial data and time $t$, which may vary from line to line. We use the notation $C(\ast,...,\ast)$ to indicate a constant that depends only on the quantities enclosed in parentheses.
We also write $A_\theta=\|\theta\|_{L^{\infty}}$ and
$A_0=\|\theta_0\|_{L^{\infty}}$ when displaying the placement of the input norms.
We write the homogeneous $[f]_{\dot C^\gamma}
=\sup_{x\ne y}\frac{|f(x)-f(y)|}{|x-y|^{\gamma}}$ and the inhomogeneous one $\|f\|_{ C^\gamma}=\|f\|_{L^{\infty}}+[f]_{\dot C^\gamma}$.
 The notation $X\lesssim Y$ means $X\leq CY$;
any additional dependence of $C$ is specified when needed.

We use the notation
\[
\omega=\partial_1u_2-\partial_2u_1,\qquad
\nabla^\perp=(-\partial_2,\partial_1),\qquad
\Kop=\nabla^\perp\Delta^{-1}.
\]
Our Fourier transform and its inverse are normalized by
\[
\widehat f(\xi)=\int_{\RR^2}e^{-ix\cdot\xi}f(x)\,dx,
\qquad
f(x)=\frac{1}{(2\pi)^2}\int_{\RR^2}e^{ix\cdot\xi}\widehat f(\xi)\,d\xi.
\]
Thus $\widehat{\partial_jf}(\xi)=i\xi_j\widehat f(\xi)$. In particular,
\[
\widehat{\Kop f}(\xi)=-i\xi^\perp|\xi|^{-2}\widehat f(\xi).
\]

The inhomogeneous Littlewood--Paley operators are denoted by $\Delta_j$,
$j\geq-1$, and $S_M$ is a smooth low-pass operator at the dyadic scale $M$.
We also write $\Delta_M$ for a block at the dyadic scale $M$.
Homogeneous blocks, denoted by $\dot\Delta_j$, are used when the
output frequency approaches zero.
The expression $\widetilde\Delta_j$ denotes a fixed finite enlargement of
$\Delta_j$. All sums over $M,N,K,U,L,P$ below are dyadic sums.
For $0<\gamma<1$, $C^\gamma$ means the bounded H\"older space, equivalently
$B_{\infty,\infty}^{\gamma}$. Constants may depend on the fixed cutoffs.

For completeness, the inhomogeneous Besov norm is
\[
\|f\|_{B_{p,r}^a}
=\big\|\big(2^{ja}\|\Delta_jf\|_p\big)_{j\geq-1}\big\|_{\ell^r},
\]
with the usual modification for $r=\infty$. We use the standard
Bernstein inequalities, paraproduct estimates, and Sobolev
embeddings associated with this decomposition; see \cite{BCD}.
The local $H^s$ construction and its continuation criterion follow
the standard Sobolev energy method for transport systems, as in
\cite{ChaeNam,MB}; the nonnegative dissipative terms are retained
in the corresponding energy estimates.

\subsection{Basic energy and maximum-principle bounds}

All estimates are first performed for smooth solutions of the exact
system \eqref{eq:system}. The final passage to $H^s$ data uses smoothing
of the initial data.
The scalar maximum principle follows from the pointwise convexity
inequalities for fractional diffusion in \cite{CC}.
If $T$ is finite and precedes a possible maximal time,
the elementary energy and maximum-principle bounds are
\begin{align}
\norm{\theta(t)}{L^p}&\leq \norm{\theta_0}{L^p},
 \ 2\leq p\leq\infty,\nonumber\\
\norm{\theta(t)}{L^2}^2+
2\int_0^t\norm{\Lambda^{\frac{\beta}{2}}\theta(\tau)}{L^2}^2\,d\tau
&=\norm{\theta_0}{L^2}^2,\label{eq:theta-energy}\\
\norm{u(t)}{L^2}&\leq \norm{u_0}{L^2}+t\norm{\theta_0}{L^2}.
\nonumber
\end{align}
In the sequel,
\begin{equation}
Z(t)=1+\norm{\Lambda^{\frac{\beta}{2}}\theta(t)}{L^2},
\qquad \int_0^T Z(t)^2\,dt\leq C(T,\norm{\theta_0}{L^2}).\nonumber
\end{equation}

\begin{lemma}[Temperature interpolation]\label{lem:theta-GN}
If $2\leq b<\infty$, $1\leq r\leq\infty$, and
\begin{equation}
0<s<\frac{\beta}{b},\nonumber
\end{equation}
then it holds
\begin{equation}\label{eq:theta-GN}
\norm{\theta}{\Bes brs}\leq C_0(1+Z^{\frac{2}{b}}).
\end{equation}
\end{lemma}

\begin{proof}
For $j\geq0$, interpolation between $L^2$ and $L^\infty$ gives
\[
\norm{\Delta_j\theta}{L^b}
\lesssim \norm{\Delta_j\theta}{L^\infty}^{1-\frac{2}{b}}
\norm{\Delta_j\theta}{L^2}^{\frac{2}{b}}
\lesssim \norm{\theta}{L^\infty}^{1-\frac{2}{b}}
2^{-\frac{j\beta}{b}}\norm{\Lambda^{\frac{\beta}{2}}\theta}{L^2}^{\frac{2}{b}}.
\]
The sequence $2^{-j(\frac{\beta}{b}-s)}$ is summable, and the low block is
bounded using the fact
$$\norm{\theta}{L^b}\leq
\norm{\theta}{L^\infty}^{1-\frac{2}{b}}\norm{\theta}{L^2}^{\frac{2}{b}}.$$
This proves \eqref{eq:theta-GN}.
\end{proof}

\section{The quadratic normal form}
\label{sec:normalform}

In this section we assume $0<\beta<\alpha<1$.
These inequalities hold throughout the parameter regions of the
main theorem.
Choose a radial $\chi_h\in C^\infty$ which vanishes below a sufficiently
large scale $R_0>1$ and equals one above $2R_0$, with
$R_0^{\beta-\alpha}\leq \frac{1}{2}$. Define
\begin{equation*}
 R_{\alpha} =\chi_h(D)\,\partial_1(\Lambda^\alpha-\Lambda^\beta)^{-1},
\qquad
\widehat{ R_{\alpha}  f}(\xi)=iA(\xi)\widehat f(\xi),\quad
A(\xi)=\frac{\chi_h(\xi)\xi_1}{|\xi|^\alpha-|\xi|^\beta}.
\end{equation*}
The multiplier is extended by zero to the region where $\chi_h$ vanishes.
$A$ is real and odd, and satisfies
\begin{equation*}
|\partial_\xi^\nu A(\xi)|\leq C_\nu|\xi|^{1-\alpha-|\nu|}.
\end{equation*}
We let
\begin{equation*}
G=\omega- R_{\alpha} \theta,\qquad \Top=\Kop R_{\alpha} .
\end{equation*}
The multiplier of $\Top$ is
$\frac{\xi^\perp A(\xi)}{|\xi|^2}$, a real even vector symbol of order $-\alpha$.
We retain these exact symbols in the cancellation estimates below.
Because $(\Lambda^\alpha-\Lambda^\beta) R_{\alpha} =\chi_h(D)\partial_1$,
combining the vorticity equation \eqref{eq:vorticity-intro} and the $\theta$-equation $\eqref{eq:system}_{2}$ gives the exact identity
\begin{equation}\label{eq:G-equation}
(\Dt+\Lambda^\alpha)G
=[ R_{\alpha} ,u\cdot\nabla]\theta+(1-\chi_h(D))\partial_1\theta,
\end{equation}
where
\begin{equation}\label{ascvbm23}
[ R_{\alpha} ,v\cdot\nabla]\theta= R_{\alpha} (v\cdot\nabla\theta)
-v\cdot\nabla R_{\alpha} \theta.
\end{equation}
Since $\omega=G+ R_{\alpha} \theta$, the Biot--Savart law gives
\begin{equation*}
u=\Kop\omega=\Kop G+\Kop R_{\alpha} \theta=\Kop G+\Top\theta.
\end{equation*}
Substituting this decomposition into the commutator in
\eqref{ascvbm23}, we obtain
\begin{equation}\label{eq:G-commutator-split}
[ R_{\alpha} ,u\cdot\nabla]\theta
=[ R_{\alpha} ,\Kop G\cdot\nabla]\theta+\mathcal C_2(\theta,\theta),
\end{equation}
where $\mathcal C_2(\theta,\theta)$ is given by
\begin{align*}
\mathcal C_2(\theta,\theta)
 \triangleq[ R_{\alpha} ,\Top\theta\cdot\nabla]\theta = R_{\alpha} \big((\Top\theta)\cdot\nabla\theta\big)
  -(\Top\theta)\cdot\nabla( R_{\alpha} \theta).
\end{align*}
With $G$ and $\theta$ viewed as the unknowns, the first term in
\eqref{eq:G-commutator-split} is a mixed $G$--$\theta$ term.
Since $\Top$ and $ R_{\alpha} $ are linear, the second term is quadratic
in $\theta$ alone. Thus $\mathcal C_2$ is precisely the purely
temperature-quadratic contribution to \eqref{eq:G-equation};
the remaining low-frequency forcing there is linear in $\theta$.

\medskip\noindent\emph{Why $\mathcal C_2$ is difficult.}
Its commutator form gives no dissipative sign in the $L^m$ test
of the $G$ equation. The forcing must be estimated through
\[
\big|\langle\mathcal C_2(\theta,\theta),|G|^{m-2}G\rangle\big|
\leq \|\Lambda^{-\frac{\alpha}{2}}\mathcal C_2(\theta,\theta)\|_{L^{p}}
 \|\Lambda^{\frac{\alpha}{2}}(|G|^{m-2}G)\|_{L^{\frac{p}{p-1}}}.
\]
Let $c_2(k,l)$ be the symmetric symbol in \eqref{eq:c2-symbol} below.
For balanced high frequencies $|k|\simeq|l|\simeq|k+l|\simeq N$,
\[
|c_2(k,l)|\lesssim N^{2-2\alpha},\qquad
|k+l|^{-\frac{\alpha}{2}}|c_2(k,l)|\lesssim N^{\frac{4-5\alpha}{2}}.
\]
Thus, even after transferring $\frac{\alpha}{2}$ derivatives to the
dissipative test, this direct bound costs $\frac{4-5\alpha}{2}$ temperature
derivatives. Consequently, the direct balanced-frequency estimate cannot be
closed using the basic temperature energy throughout our range.

\medskip\noindent\emph{Why split into nonresonant and resonant parts.}
We seek a bilinear correction $\Bop(\theta,\theta)$ to cancel this
forcing in the equation for $H=G-\Bop(\theta,\theta)$, (see \eqref{eq:H-equation} for details).
Formally, matching the quadratic diffusion terms requires
its symbol to solve
\[
\underbrace{\bigl(|k+l|^\alpha-|k|^\beta-|l|^\beta\bigr)}_{D(k,l)}
\,b(k,l)=c_2(k,l),\qquad
b=\frac{c_2}{D}.
\]
This division is controlled only away from resonance. With
$S=|k+l|^\alpha+|k|^\beta+|l|^\beta$ and fixed small $c>0$,
\[
|D|\geq cS
\quad\Longrightarrow\quad
\left|\frac{c_2}{D}\right|\leq\frac{|c_2|}{cS},
\qquad S\geq|k+l|^\alpha.
\]
The nonresonant part therefore admits the correction defined
below. On balanced blocks, the correction has order $2-3\alpha$,
gaining $\alpha$ derivatives relative to the forcing of order
$2-2\alpha$. Subtracting it removes $\Cnr$ exactly from the
$H$ equation.

Near $D=0$, division by $D$ is not uniformly bounded. At high
output, however, near resonance forces the geometry
\[
N=|k+l|,\quad K=\max\{|k|,|l|\},\qquad
N^\alpha\simeq K^\beta,\quad K\gg N,
\]
as proved in \eqref{eq:res-scales}. The two large inputs are
nearly opposite, and the symmetric symbol has the cancellation
\[
|c_2(k,l)|\lesssim K^{2-2\alpha}
\left[\left(\frac NK\right)^{2-\alpha}
     +\left(\frac NK\right)^2\right]
.\]
We therefore retain this resonant part and estimate it using
the small ratio $\frac NK$. A fixed low-output part is kept separately.
The decomposition below implements these two mechanisms.

The symmetric symbol of $\mathcal C_2$ is
\begin{equation}\label{eq:c2-symbol}
c_2(k,l)=\frac{k^\perp\cdot l}{2}
\left[A(k)A(l) \left( \frac{1}{|k|^2} - \frac{1}{|l|^2} \right) - A(q) \left( \frac{A(k)}{|k|^2} - \frac{A(l)}{|l|^2} \right)
\right],
\end{equation}
where $q=k+l$.
Let us show the detailed derivation of the symmetric form \eqref{eq:c2-symbol}.
The nonlinear commutator $\mathcal{C}_2(\theta,\theta) = [ R_{\alpha} , \mathcal{T}\theta \cdot \nabla]\theta$ reads as follows
\begin{equation} \label{eq:commutator}
    \mathcal{C}_2(\theta,\theta) =  R_{\alpha} (\mathcal{T}\theta \cdot \nabla\theta) - (\mathcal{T}\theta \cdot \nabla)( R_{\alpha} \theta).
\end{equation}
Using the Fourier transform convention $\widehat{\nabla f}(\xi) = i\xi \widehat{f}(\xi)$, we have
\begin{align*}
    \widehat{ R_{\alpha}  f}(\xi)  = i A(\xi) \widehat{f}(\xi),\qquad
    \widehat{\mathcal{K} f}(\xi)  = -i \frac{\xi^\perp}{|\xi|^2} \widehat{f}(\xi).
\end{align*}
Since $\mathcal{T} = \mathcal{K} R_{\alpha} $, its multiplier is the product of the two operators
\begin{equation*}
    \widehat{\mathcal{T} f}(\xi) = \left(-i \frac{\xi^\perp}{|\xi|^2}\right) \big(i A(\xi)\big) \widehat{f}(\xi) = \frac{\xi^\perp A(\xi)}{|\xi|^2} \widehat{f}(\xi).
\end{equation*}
In the frequency domain, a product of two functions corresponds to a convolution. Let $k$ be the frequency of the velocity vector field (i.e., $\mathcal{T}\theta$) and $l$ be the frequency of the scalar field being transported (i.e., $\theta$ or $ R_{\alpha} \theta$). The output frequency is $q = k+l$.
We analyze the two terms in \eqref{eq:commutator} separately.

Regarding the term $ R_{\alpha} (\mathcal{T}\theta \cdot \nabla\theta)$,
the inner dot product $\mathcal{T}\theta \cdot \nabla\theta$ in the frequency domain yields the multiplier
\begin{equation*}
    \left(\frac{k^\perp A(k)}{|k|^2}\right) \cdot (il) = i \frac{k^\perp \cdot l}{|k|^2} A(k).
\end{equation*}
The outer operator $ R_{\alpha} $ acts on the sum of the frequencies $q = k+l$, introducing a factor of $i A(k+l)$. Multiplying these together (noting $i^2 = -1$), we obtain the symbol for the first term
\begin{equation*}
    \text{Sym}_1 = - \frac{k^\perp \cdot l}{|k|^2} A(k) A(k+l).
\end{equation*}

Regarding the term $(\mathcal{T}\theta \cdot \nabla)( R_{\alpha} \theta)$, the velocity and the gradient still provide the multiplier $i \frac{k^\perp \cdot l}{|k|^2} A(k)$. However, the operator $ R_{\alpha} $ now acts \textit{only} on the second function $\theta$ (with frequency $l$), which provides an additional factor of $i A(l)$. Multiplying these gives the symbol for the second term
\begin{equation*}
    \text{Sym}_2 = - \frac{k^\perp \cdot l}{|k|^2} A(k) A(l).
\end{equation*}
Subtracting the second term from the first, we obtain
\begin{align*}
    c(k,l) &= \left( - \frac{k^\perp \cdot l}{|k|^2} A(k) A(k+l) \right) - \left( - \frac{k^\perp \cdot l}{|k|^2} A(k) A(l) \right) \\
    &= \frac{k^\perp \cdot l}{|k|^2} A(k) \big[ A(l) - A(k+l) \big].
\end{align*}
Since the nonlinear term $\mathcal{C}_2(\theta,\theta)$ takes the same function $\theta$ for both inputs, we must symmetrize the symbol with respect to the input frequencies $k$ and $l$ to capture the exact cancellation structures. The symmetric symbol $c_2(k,l)$ is defined as
\begin{equation*}
    c_2(k,l) = \frac{1}{2} \big[ c(k,l) + c(l,k) \big].
\end{equation*}
First, we write out $c(l,k)$ by exchanging $k$ and $l$
\begin{equation*}
    c(l,k) = \frac{l^\perp \cdot k}{|l|^2} A(l) \big[ A(k) - A(l+k) \big].
\end{equation*}
We observe a crucial geometric anti-symmetry in the dot product. Let $k=(k_1, k_2)$ and $l=(l_1, l_2)$, then $l^\perp = (-l_2, l_1)$. It follows that
\begin{equation*}
    l^\perp \cdot k = -l_2 k_1 + l_1 k_2 = - (k_1 l_2 - k_2 l_1) = - (k^\perp \cdot l).
\end{equation*}
Substituting this anti-symmetry back into $c(l,k)$, we get
\begin{equation*}
    c(l,k) = - \frac{k^\perp \cdot l}{|l|^2} A(l) \big[ A(k) - A(q) \big].
\end{equation*}
Adding $c(k,l)$ and $c(l,k)$ and dividing by $2$, we are able to conclude
\begin{equation}\label{asvbut897}
    c_2(k,l) = \frac{k^\perp \cdot l}{2} \left[ \frac{A(k)}{|k|^2}\big(A(l) - A(q)\big) - \frac{A(l)}{|l|^2}\big(A(k) - A(q)\big) \right].
\end{equation}
To reveal the deep cancellation property, we expand the terms inside the square brackets
\begin{equation*}
    \left[ \frac{A(k)A(l)}{|k|^2} - \frac{A(k)A(q)}{|k|^2} - \frac{A(l)A(k)}{|l|^2} + \frac{A(l)A(q)}{|l|^2} \right].
\end{equation*}
Finally, we group the terms containing $A(k)A(l)$ together, and the terms containing $-A(q)$ together
\begin{equation*}
    A(k)A(l) \left( \frac{1}{|k|^2} - \frac{1}{|l|^2} \right) - A(q) \left( \frac{A(k)}{|k|^2} - \frac{A(l)}{|l|^2} \right).
\end{equation*}
Substituting this rearranged expression back into the symmetrized formulation, we arrive exactly at \eqref{eq:c2-symbol}.
This specific grouping is mathematically significant because the differences in the parentheses approach zero when $|k| \approx |l|$, showcasing an intrinsic frequency cancellation mechanism.

\begin{lemma}
Regarding $c_2(k,l)$, it holds
\begin{equation}\label{cest001}
|c_2(k,l)|\lesssim \left\{\aligned
&N^{2-\alpha} K^{-\alpha},\qquad\,\qquad\quad\ \ \ |k| \simeq |l| \simeq K \gg |q| \simeq N, \\
&N^{1-\alpha} L^{1-\alpha}, \ \qquad \qquad\qquad |k| \simeq L \ll |l| \simeq |q| \simeq N,\\
&N^{2-2\alpha},\qquad   \qquad \qquad\qquad  |k| \simeq |l| \simeq |q| \simeq N.
\endaligned\right.
\end{equation}
\end{lemma}

\begin{proof}
The multiplier $A(\xi)$ is real and odd (i.e., $A(-\xi) = -A(\xi)$), satisfying the bound
\begin{equation*}
|A(\xi)| \lesssim |\xi|^{1-\alpha}, \quad |\nabla A(\xi)| \lesssim |\xi|^{-\alpha}.
\end{equation*}
Let us define $F(\xi) = \frac{A(\xi)}{|\xi|^2}$. Since $A(\xi)$ is odd, $F(\xi)$ is also odd. It satisfies
\begin{equation*}
|F(\xi)| \lesssim |\xi|^{-1-\alpha}.
\end{equation*}
It thus follows from \eqref{asvbut897} that
\begin{equation}\label{asnpk35}
c_2(k,l)=\frac{k^\perp\cdot l}{2}
\left[
F(k)\big\{A(l)-A(k+l)\big\}
+F(l)\big\{A(k+l)-A(k)\big\}
\right].
\end{equation}
Also, note the geometric identity for the cross product. Since $l = q-k$ and $k^\perp \cdot k = 0$, we always have
\begin{equation}\label{asfvnp98}
k^\perp \cdot l = k^\perp \cdot (q-k) = k^\perp \cdot q.
\end{equation}

\textbf{The proof of $\eqref{cest001}_{1}$: High-high to low ($|k| \simeq |l| \simeq K \gg |q| \simeq N$)}.
In this regime, the two input frequencies are much larger than the output frequency. We use the original form \eqref{eq:c2-symbol} to bound $c_2(k,l)$.

\paragraph{Step 1: Bound the $A(q)$ term.}
The second term in the bracket of \eqref{eq:c2-symbol} is $A(q)(F(k)-F(l))$.
Since $l = q-k$ and $|k| \simeq |k-q| \simeq K$, it yields
\begin{align*}
|F(k) - F(l)| &= |F(k)-F(q-k)| \\ &\leq |F(k)| + |F(q-k)| \\ &\lesssim K^{-1-\alpha} + K^{-1-\alpha} \\ & \lesssim K^{-1-\alpha}.
\end{align*}

\paragraph{Step 2: Bound the pre-factor and the inverse-square difference.}
Using the identity \eqref{asfvnp98}, the pre-factor is bounded by
\[
|k^\perp \cdot l| = |k^\perp \cdot q| \leq |k||q| \simeq KN.
\]
For the inverse-square difference, we write:
\[
\left| \frac{1}{|k|^2} - \frac{1}{|l|^2} \right| = \frac{\big| |l|^2 - |k|^2 \big|}{|k|^2|l|^2}.
\]
Notice that
 $$|l|^2 - |k|^2 = |q-k|^2 - |k|^2 = |q|^2 - 2q\cdot k,\quad |q| \simeq N, |k| \simeq K,$$
 we have
$$\big| |l|^2 - |k|^2 \big| \lesssim N^2 + NK \simeq NK.$$
As a result, one gets
\[
\left| \frac{1}{|k|^2} - \frac{1}{|l|^2} \right| \lesssim \frac{NK}{K^4} = NK^{-3}.
\]

\paragraph{Step 3: Combine all pieces.}
We recall
$$|A(k)A(l)| \lesssim K^{1-\alpha}K^{1-\alpha} = K^{2-2\alpha},\qquad |A(q)| \lesssim N^{1-\alpha}.$$
 Substituting all bounds into \eqref{eq:c2-symbol} implies
\begin{align*}
|c_2(k,l)| &\lesssim KN \left[ \underbrace{K^{2-2\alpha}}_{|A(k)A(l)|} \underbrace{NK^{-3}}_{\text{diff of squares}} + \underbrace{N^{1-\alpha}}_{|A(q)|} \underbrace{K^{-1-\alpha}}_{|F(k)-F(l)|} \right] \\
&= KN \left[ N K^{-1-2\alpha} + N^{1-\alpha} K^{-1-\alpha} \right] \\
&= N^2 K^{-2\alpha} + N^{2-\alpha} K^{-\alpha}\\
&\lesssim N^{2-\alpha} K^{-\alpha},
\end{align*}
where in the last line we have used $N \ll K$.

\textbf{The proof of $\eqref{cest001}_{2}$: Low-high to high ($|k| \simeq L \ll |l| \simeq |q| \simeq N$)}.
In this regime, one input frequency is small, while the other is comparable to the output frequency. We use the alternative difference form \eqref{asnpk35}.

\paragraph{Step 1: Bound the increments of $A$ via the Mean Value Theorem.}
Since $L \ll N$, we can bound the difference $A(l+k) - A(l)$ by integrating its gradient along the line segment from $l$ to $l+k$. More precisely, we obtain
\[
A(l+k) - A(l) = \int_0^1 k \cdot \nabla A(l+tk) \, dt.
\]
Because $|l+tk| \simeq N$ for all $t \in [0,1]$, we have $|\nabla A(l+tk)| \lesssim N^{-\alpha}$. Thus, we have
\[
|A(l+k) - A(l)| \leq |k| \sup_{t} |\nabla A| \lesssim L N^{-\alpha}.
\]

\paragraph{Step 2: Bound the two terms in the bracket of \eqref{asnpk35}.}
\begin{itemize}
    \item \textbf{First term:} $F(k)\{A(l)-A(k+l)\}$.\\
    We have $|F(k)| \lesssim L^{-1-\alpha}$, and the difference is bounded by $LN^{-\alpha}$.\\
    $$ \Big| F(k)\{A(l)-A(k+l)\} \Big| \lesssim L^{-1-\alpha} \times L N^{-\alpha} = L^{-\alpha} N^{-\alpha}. $$

    \item \textbf{Second term:} $F(l)\{A(k+l)-A(k)\}$.\\
    Here we don't use the derivative for the difference, but simply the triangle inequality because $A$ grows sublinearly
    $$|A(k+l) - A(k)| \leq |A(k+l)| + |A(k)| \lesssim N^{1-\alpha} + L^{1-\alpha}.$$\\
    Also, $|F(l)| \lesssim N^{-1-\alpha}$.It thus gives
    $$ \Big| F(l)\{A(k+l)-A(k)\} \Big| \lesssim N^{-1-\alpha} \left( N^{1-\alpha} + L^{1-\alpha} \right) = N^{-2\alpha} + N^{-1-\alpha}L^{1-\alpha}. $$
\end{itemize}

\paragraph{Step 3: Combine and factor out.}
The pre-factor is $|k^\perp \cdot l| \leq |k||l| \simeq LN$. Multiplying this with the sum of the two terms from Step 2, one arrives at
\begin{align*}
|c_2(k,l)| &\lesssim LN \left[ L^{-\alpha} N^{-\alpha} + N^{-2\alpha} + N^{-1-\alpha}L^{1-\alpha} \right] \\
&= L^{1-\alpha}N^{1-\alpha} + L N^{1-2\alpha} + L^{2-\alpha}N^{-\alpha}.
\end{align*}
To see the underlying structure, we factor out the leading term $N^{1-\alpha}L^{1-\alpha}$ to derive
\begin{align*}
|c_2(k,l)| &\lesssim N^{1-\alpha} L^{1-\alpha} \left[ 1 + \frac{L N^{1-2\alpha}}{L^{1-\alpha}N^{1-\alpha}} + \frac{L^{2-\alpha}N^{-\alpha}}{L^{1-\alpha}N^{1-\alpha}} \right] \\
&= N^{1-\alpha} L^{1-\alpha} \left[ 1 + \left(\frac{L}{N}\right)^\alpha + \frac{L}{N} \right].
\end{align*}
Since $L \ll N$, the ratios $\left(\frac{L}{N}\right)^\alpha$ and $\frac{L}{N}$ are strictly less than $1$. Therefore, the bracket is bounded by a constant, yielding:
\[
|c_2(k,l)| \lesssim N^{1-\alpha} L^{1-\alpha}.
\]

\textbf{The proof of $\eqref{cest001}_{3}$: Comparable frequencies ($|k| \simeq |l| \simeq |q| \simeq N$)}.
This is the simplest regime since there are no cancellations needed to offset large disparities in scales. We use form \eqref{asnpk35} again.
The pre-factor is simply $|k^\perp \cdot l| \lesssim N^2$.
For the terms in the bracket, one gets
\[
|F(k)| \simeq |F(l)| \lesssim N^{-1-\alpha}.
\]
\[
|A(l) - A(k+l)| \leq |A(l)| + |A(q)| \lesssim N^{1-\alpha}.
\]
Thus, the bracket is bounded by $N^{-1-\alpha} \times N^{1-\alpha} = N^{-2\alpha}$. \\
Multiplying by the pre-factor gives:
\[
|c_2(k,l)| \lesssim N^2 \times N^{-2\alpha} = N^{2-2\alpha}.
\]
This completes all localized estimates for the symbol $c_2(k,l)$.
\end{proof}

\subsection{Removing the nonresonant quadratic term}
\label{sec:normal-form-definition}
We denote
\[
D(k,l)\triangleq|k+l|^\alpha-|k|^\beta-|l|^\beta,\qquad
S(k,l)\triangleq|k+l|^\alpha+|k|^\beta+|l|^\beta.
\]
\noindent\emph{Frequencies and decay rates.}
\[
\underbrace{k,\ l}_{\text{temperature inputs}}
\ \longmapsto\
\underbrace{q=k+l}_{\text{correction output}},
\qquad
\begin{aligned}
\lambda_{\mathrm{in}}&=|k|^\beta+|l|^\beta,\\
\lambda_{\mathrm{out}}&=|q|^\alpha.
\end{aligned}
\]
Thus exact resonance means $D=0$, or
$\lambda_{\mathrm{out}}=\lambda_{\mathrm{in}}$.
Choose $0<c<\frac{1}{10}$ and an even $\eta\in C^\infty(\RR;[0,1])$ with
\[
\eta(z)=
\begin{cases}
0,&|z|\leq c,\\
1,&|z|\geq2c,
\end{cases}
\qquad
\chi_{\mathrm{nr}}(k,l)=\eta\!\left(\frac{D}{S}\right).
\]
The relative gap $\frac{|D|}{S}$ determines the smooth decomposition:
\[
\begin{aligned}
|D|\leq cS&\quad\Longrightarrow\quad
 \chi_{\mathrm{nr}}=0
 &&\text{(near resonance)},\\
|D|\geq2cS&\quad\Longrightarrow\quad
 \chi_{\mathrm{nr}}=1
 &&\text{(nonresonance)},\\
\chi_{\mathrm{nr}}\ne0&\quad\Longrightarrow\quad
 |D|\geq cS,\qquad
 \left|\frac{\chi_{\mathrm{nr}}}{D}\right|\leq\frac1{cS}.
\end{aligned}
\]

\medskip\noindent\emph{The output cutoff and quadratic correction.}
Choose the smooth scalar profiles
\[
\varrho(t)=
\begin{cases}
0,&t\leq0,\\
e^{-\frac{1}{t}},&t>0,
\end{cases}
\qquad
\varphi(r)=
\frac{\varrho(r-1)}
     {\varrho(r-1)+\varrho(2-r)}.
\]
Fix $N_0$ so that
\[
C_c=\frac{1+2c}{1-2c},\qquad
N_0>\max\{2R_0,1\},\qquad
N_0^{\alpha-\beta}\geq2C_c8^\beta,
\]
and define the radial cutoff of the output $q=k+l$ by
\[
\chi_{\mathrm{out}}(q)=\varphi\!\left(\frac{|q|}{N_0}\right).
\]
Then, we have
\[
\begin{gathered}
\chi_{\mathrm{out}}(q)=
\begin{cases}
0,&|q|\leq N_0,\\
1,&|q|\geq2N_0,
\end{cases}
\qquad 0\leq\chi_{\mathrm{out}}\leq1,\\
\chi_{\mathrm{out}}\in C^\infty(\RR^2),\qquad
|\partial_q^\nu\chi_{\mathrm{out}}(q)|
\leq C_\nu N_0^{-|\nu|},\\
\supp\chi_{\mathrm{out}}\subset\{|q|\geq N_0\},\qquad
\supp(1-\chi_{\mathrm{out}})\subset\{|q|\leq2N_0\}.
\end{gathered}
\]
The Biot--Savart factor satisfies
$\sup_{q\ne0}|q|^{-1}\chi_{\mathrm{out}}(q)\leq N_0^{-1}$.
We set
\begin{equation}\label{eq:b-symbol}
\begin{aligned}
b(k,l)
 \triangleq\chi_{\mathrm{out}}(q)
  \frac{\chi_{\mathrm{nr}}(k,l)c_2(k,l)}{D(k,l)} =\chi_{\mathrm{out}}(q)\frac{c_2(k,l)}{S(k,l)}
  g\!\left(\frac{D(k,l)}{S(k,l)}\right),
\end{aligned}
\end{equation}
where $g(z)=\frac{\eta(z)}{z}$ extends smoothly by zero near $z=0$.
For balanced frequencies,
\[
\begin{gathered}
|k|\simeq|l|\simeq|q|\simeq N\gg1,\quad\alpha>\beta
\quad\Longrightarrow\quad |D|\simeq N^\alpha,\\
|b(k,l)|\lesssim
\underbrace{N^{-\alpha}}_{\text{gain of }\alpha\text{ derivatives}}
|c_2(k,l)|\lesssim N^{2-3\alpha}.
\end{gathered}
\]
For $f,g\in\mathcal S(\RR^2)$, we define
\begin{equation}\label{eq:B-definition}
\widehat{\Bop(f,g)}(\xi)
=\frac1{(2\pi)^2}\int_{\RR^2}
b(k,\xi-k)\widehat f(k)\widehat g(\xi-k)\,dk
\end{equation}or equivalently
\begin{equation*}
{\Bop(f,g)}(x)
=\frac{1}{(2\pi)^4}\int_{\RR^2}\int_{\RR^2}
e^{i x\cdot(k+l)}b(k,l)\widehat f(k)\widehat g(l)\,dk\,dl.
\end{equation*}
We point out that $\Bop(f,g)$ is not a simple product of fractional derivatives like $\Lambda^{\delta_{1}}f\Lambda^{\delta_{2}}g$. Instead, it is an intrinsically non-separable, non-local bilinear fourier multiplier of Coifman--Meyer type with total scaling order $2-3\alpha$, specifically tailored via the QNF to algebraically eliminate non-resonant bilinear interactions.
Since $b(k,l)=b(l,k)$, we have $\Bop(f,g)=\Bop(g,f)$.

\medskip\noindent\emph{The quadratic form and its decomposition.}
The symmetric polarization of $\mathcal C_2(\theta,\theta)$ is
\[
\mathcal C_2(f,g)
=\frac12\bigl([ R_{\alpha} ,\Top f\cdot\nabla]g
             +[ R_{\alpha} ,\Top g\cdot\nabla]f\bigr).
\]
For $q=k+l$, Fourier multiplication and interchange of the two
inputs give
\[
\begin{aligned}
m_2(k,l)
&=\frac{k^\perp\cdot l}{|k|^2}A(k)\,[A(l)-A(q)],\\
c_2(k,l)
&=\frac12\bigl(m_2(k,l)+m_2(l,k)\bigr)
 =c_2(l,k),
\end{aligned}
\]
and hence
\[
\widehat{\mathcal C_2(\theta,\theta)}(\xi)
=\frac1{(2\pi)^2}\int_{\RR^2}
c_2(k,\xi-k)\widehat\theta(k)\widehat\theta(\xi-k)\,dk.
\]
The expanded symmetric symbol is \eqref{eq:c2-symbol}.
We first split by output size, then split the high-output part
by the relative dissipation gap
\[
1=
\underbrace{\chi_{\mathrm{out}}(q)}_{\text{high output}}
\left[
\underbrace{\chi_{\mathrm{nr}}(k,l)}_{\text{nonresonant}}
+
\underbrace{1-\chi_{\mathrm{nr}}(k,l)}_{\text{near resonant}}
\right]
+
\underbrace{1-\chi_{\mathrm{out}}(q)}_{\text{low output}}.
\]
This smooth partition gives
\[
\mathcal C_2(\theta,\theta)
=\Cnr(\theta,\theta)+\Cres(\theta,\theta)+\Clow(\theta,\theta).
\]
Explicitly, for $x\in\RR^2$, one derives
\begin{equation}\label{eq:split-c2}
\begin{aligned}
\Cnr(\theta,\theta)(x)
&=\frac1{(2\pi)^4}\iint_{\RR^2\times\RR^2}
 e^{ix\cdot(k+l)}\chi_{\mathrm{out}}(k+l)\chi_{\mathrm{nr}}(k,l)
 c_2(k,l)\widehat\theta(k)\widehat\theta(l)\,dk\,dl,\\[2pt]
\Cres(\theta,\theta)(x)
&=\frac1{(2\pi)^4}\iint_{\RR^2\times\RR^2}
 e^{ix\cdot(k+l)}\chi_{\mathrm{out}}(k+l)
 \bigl(1-\chi_{\mathrm{nr}}(k,l)\bigr)
 c_2(k,l)\widehat\theta(k)\widehat\theta(l)\,dk\,dl,\\[2pt]
\Clow(\theta,\theta)(x)
&=\frac1{(2\pi)^4}\iint_{\RR^2\times\RR^2}
 e^{ix\cdot(k+l)}\bigl(1-\chi_{\mathrm{out}}(k+l)\bigr)
 c_2(k,l)\widehat\theta(k)\widehat\theta(l)\,dk\,dl.
\end{aligned}
\end{equation}
The three terms are treated according to their frequency supports.

\medskip\noindent\emph{Nonresonant part: removal by the normal form.}
\[
\begin{gathered}
\chi_{\mathrm{out}}(q)\chi_{\mathrm{nr}}(k,l)\ne0
\quad\Longrightarrow\quad |q|\geq N_0,\quad |D|\geq cS,\\
Db=\chi_{\mathrm{out}}\chi_{\mathrm{nr}}c_2,\qquad
|b|\leq\frac{|c_2|}{cS}.
\end{gathered}
\]
Thus $\Cnr$ is precisely the quadratic term canceled by the
change $G\mapsto G-\Bop(\theta,\theta)$ below.

\medskip\noindent\emph{Resonant part: gain from a smaller output.}
On the support of $\chi_{\mathrm{out}}(1-\chi_{\mathrm{nr}})$, put
$X=|q|^\alpha$, $Y=|k|^\beta+|l|^\beta$,
$N=|q|\geq N_0$, and $K=\max\{|k|,|l|\}$. Then
\[
|X-Y|\leq2c(X+Y)
\quad\Longrightarrow\quad
C_c^{-1}\leq\frac XY\leq C_c,\qquad
N^\alpha\leq2C_cK^\beta.
\]
Our choice of $N_0$ gives
\begin{equation}\label{eq:res-scales}
\begin{gathered}
N^\alpha\simeq K^\beta,\qquad
K\simeq N^{\frac{\alpha}{\beta}},\\
\frac KN\geq
\left(\frac{N_0^{\alpha-\beta}}{2C_c}\right)^{\frac{1}{\beta}}
\geq8,\qquad
\min\{|k|,|l|\}\geq K-N\geq\frac78K.
\end{gathered}
\end{equation}
The output is high in absolute size and smaller than both inputs.
The cancellation gives
\[
|c_2(k,l)|\lesssim
K^{2-2\alpha}
\left[\left(\frac NK\right)^{2-\alpha}
     +\left(\frac NK\right)^2\right].
\]
We retain $\Cres$ and use this gain to estimate it.

\medskip\noindent\emph{Low-output part: direct remainder estimate.}
\[
\supp\widehat{\Clow(\theta,\theta)}
\subset\{|q|\leq2N_0\}.
\]
The input frequencies may still be large:
$l=-k+q$ allows $|k|\simeq|l|\gg N_0$.
Lemma~\ref{lem:res} gives
\[
\|\Lambda^{-\frac{\alpha}{2}}\Clow(\theta,\theta)\|_p
\lesssim_{N_0,\alpha,\beta}
\|\theta\|_{L^{\infty}}\|\theta\|_p,\qquad 2\leq p<\infty.
\]
We therefore keep $\Clow$ in $F_{\mathrm{low}}$.

Based on the bilinearity of \eqref{eq:B-definition} and $\partial_t\theta=-u\cdot\nabla\theta-\Lambda^\beta\theta$, we obtain
\begin{align*}
\partial_{t}\Bop(\theta,\theta)
&=\Bop(\partial_{t}\theta,\theta)+\Bop(\theta,\partial_{t}\theta)\\& =
\Bop(-u\cdot\nabla\theta-\Lambda^\beta\theta,\theta)
+\Bop(\theta,-u\cdot\nabla\theta-\Lambda^\beta\theta)
\\& =-
\Bop(u\cdot\nabla\theta,\theta)
-\Bop( \Lambda^\beta\theta,\theta)
-\Bop(\theta,u\cdot\nabla\theta)-\Bop(\theta,\Lambda^\beta\theta),
\end{align*}
which yields
\begin{align*}
(\Dt+\Lambda^\alpha)\Bop(\theta,\theta)
= &\Lambda^\alpha\Bop(\theta,\theta)
-\Bop( \Lambda^\beta\theta,\theta)
-\Bop(\theta,\Lambda^\beta\theta)
\\&-[
\Bop(u\cdot\nabla\theta,\theta)+\Bop(\theta,u\cdot\nabla\theta)
-(u\cdot\nabla)\Bop(\theta,\theta)].
\end{align*}
For simplicity, we denote
\begin{equation}
\Qop(v;\theta,\theta)
\triangleq\Bop(v\cdot\nabla\theta,\theta)
+\Bop(\theta,v\cdot\nabla\theta)-v\cdot\nabla\Bop(\theta,\theta).\nonumber
\end{equation}
We thus derive
\begin{align}\label{eqsdy726}
(\Dt+\Lambda^\alpha)\Bop(\theta,\theta)
={}&\Lambda^\alpha\Bop(\theta,\theta)
-\Bop(\Lambda^\beta\theta,\theta)
-\Bop(\theta,\Lambda^\beta\theta)
-\Qop(u;\theta,\theta).
\end{align}
The Fourier transform of the first three terms on the right handside of \eqref{eqsdy726} is
\begin{align*}
&\frac{1}{(2\pi)^2}\int_{\RR^2}
\bigl(|\xi|^\alpha-|k|^\beta-|\xi-k|^\beta\bigr)
b(k,\xi-k)\widehat\theta(k)\widehat\theta(\xi-k)\,dk\\
&=\frac{1}{(2\pi)^2}\int_{\RR^2}
\chi_{\mathrm{out}}(\xi)\chi_{\mathrm{nr}}(k,\xi-k)c_2(k,\xi-k)
\widehat\theta(k)\widehat\theta(\xi-k)\,dk
\\ &=\widehat{\Cnr(\theta,\theta)}(\xi),
\end{align*}
where we used \eqref{eq:b-symbol} and \eqref{eq:split-c2}.
Therefore, we obtain
\begin{equation}\label{tegfd91}
(\Dt+\Lambda^\alpha)\Bop(\theta,\theta)
=\Cnr(\theta,\theta)-\Qop(u;\theta,\theta).
\end{equation}
Denoting $H\triangleq G-\Bop(\theta,\theta)\equiv \omega- R_{\alpha} \theta-\Bop(\theta,\theta)$ as the so-called corrected vorticity, we deduce from \eqref{eq:G-equation} and \eqref{tegfd91} that
 \begin{align}\label{eqsdfvcx12}
(\Dt+\Lambda^\alpha)H={}&
[ R_{\alpha} ,u_H\cdot\nabla]\theta+\Cres(\theta,\theta)
+[ R_{\alpha} ,u_B\cdot\nabla]\theta \nonumber\\
&+\Qop(u;\theta,\theta)+(1-\chi_h(D))\partial_1\theta
+\Clow(\theta,\theta).
\end{align}
The velocity decomposition is given by
\begin{align*}
u&=\nabla^\perp\Delta^{-1}\omega\nonumber\\
&=\nabla^\perp\Delta^{-1}(G+ R_{\alpha} \theta)
\nonumber\\
&=\nabla^\perp\Delta^{-1}H+\nabla^\perp\Delta^{-1}\Bop(\theta,\theta)
+\nabla^\perp\Delta^{-1} R_{\alpha} \theta
\nonumber\\
&\triangleq u_H+u_B+\Top\theta.
\end{align*}
Roughly, we have
$$u_H\backsimeq \Lambda^{-1}H,\quad u_B\backsimeq \Lambda^{-1}\Bop(\theta,\theta),\quad
\Top\theta\backsimeq \Lambda^{-\alpha}\theta.$$
Consequently, it thus follows from \eqref{eqsdfvcx12} that
\begin{align}\label{eq:H-equation}
(\Dt+\Lambda^\alpha)H={}&
[ R_{\alpha} ,u_H\cdot\nabla]\theta+\Cres(\theta,\theta)
+[ R_{\alpha} ,u_B\cdot\nabla]\theta \nonumber\\
&+\Qop(u_H;\theta,\theta)
+\Qop(\Top\theta;\theta,\theta)
+\Qop(u_B;\theta,\theta)+F_{\mathrm{low}},
\end{align}
where $F_{\mathrm{low}}$ is given by
$$
F_{\mathrm{low}}=(1-\chi_h(D))\partial_1\theta
+\Clow(\theta,\theta).$$
We point out that equation~\eqref{eq:H-equation} retains the low-output remainder
and will be the starting point for the energy estimates.

\section{Multiplier and commutator estimates}
\label{sec:localized}

We first record the localized kernel bounds for the quadratic
normal form and its resonant remainder. We then estimate the
transport commutator and the mixed vorticity--temperature term.
The use of localized symbol derivatives and integrable kernels
is consistent with the multilinear multiplier framework of
\cite{CM}; the scale-dependent gains needed here are proved explicitly.

\subsection{Localized kernels and bilinear estimates}

\begin{lemma}[Localized kernel principle for multilinear multipliers]\label{lem:kernel}
Let $r \in \{2, 3\}$, and set $d = 2r \in \{4, 6\}$. Suppose that $\mathfrak{m} \in C^\infty((\mathbb R^2)^r)$ is a multilinear Fourier multiplier. Assume there exists an invertible linear map $\mathrm{L} \in \mathrm{GL}(d, \mathbb R)$ such that the relocated symbol $\widetilde{\mathfrak{m}} \triangleq \mathfrak{m} \circ \mathrm{L}^{-1} \in C_c^\infty(\mathbb R^d)$ is localized at dyadic scales $\mathbf{M} = (M_1, \ldots, M_r) \in (2^{\mathbb Z})^r$ in the adapted coordinates $\zeta = (\zeta_1, \ldots, \zeta_r) \triangleq \mathrm{L}(\xi_1, \ldots, \xi_r) \in \mathbb R^d$, satisfying
\begin{equation}\label{eq:mudsp4689}
\operatorname{supp} \widetilde{\mathfrak{m}} \subset \prod_{\iota=1}^r \big\{\zeta_\iota \in \mathbb R^2 : |\zeta_\iota| \leq C_0 M_\iota\big\},
\end{equation}
for some universal constant $C_0 \in (0, \infty)$.

Assume further that for every multi-index $\nu = (\nu_1, \ldots, \nu_r) \in \mathbb N_0^d$ with $\nu_\iota \in \mathbb N_0^2$ satisfying $|\nu| = \sum_{\iota=1}^r |\nu_\iota| \leq 12$, there exists a constant $\mathfrak{A} > 0$ such that
\begin{equation}\label{eq:multiplier-hypo-sharp}
\left( \prod_{\iota=1}^r M_\iota^{|\nu_\iota|} \right) \sup_{\zeta \in \mathbb R^d} \big| \partial_{\zeta}^\nu \widetilde{\mathfrak{m}}(\zeta) \big| \leq C_\nu \mathfrak{A}.
\end{equation}
Then the spatial convolution kernel $\Psi= \mathcal{F}^{-1}\mathfrak{m} \in \mathcal{S}(\mathbb R^d)$ satisfies the scale-invariant bound
\begin{equation}\label{eq:kernel-L1-bound-sharp}
\|\Psi\|_{L^1(\mathbb R^d)} \leq C \mathfrak{A},
\end{equation}
where the constant $C > 0$ depends solely on $d$, $C_0$, and $\max_{|\nu|\le 12} C_\nu$, and in particular is strictly independent of the dyadic scales $\mathbf{M}$ and the linear transformation $\mathrm{L}$.

Consequently, for any exponents $1 \leq p, p_1, \ldots, p_r \leq \infty$ satisfying the conjugate reciprocal relation
\begin{equation}
\frac{1}{p} = \sum_{\iota=1}^r \frac{1}{p_\iota} \qquad \text{(with the standard convention $\frac{1}{\infty} = 0$)},\nonumber
\end{equation}
the associated multilinear operator $T_m$ given by
\begin{equation}
T_m(f_1, \ldots, f_r)(x) \triangleq \frac{1}{(2\pi)^{2r}} \int_{(\mathbb R^2)^r} e^{ix \cdot (\xi_1 + \cdots + \xi_r)} m(\xi_1, \ldots, \xi_r) \prod_{\iota=1}^r \widehat{f}_\iota(\xi_\iota) \,d\xi_1 \cdots d\xi_r,\nonumber
\end{equation}
satisfies
\begin{align}
\|T_m(f_1, \ldots, f_r)\|_{L^p(\mathbb R^2)} &\leq C \mathfrak{A} \prod_{\iota=1}^r \|f_\iota\|_{L^{p_\iota}(\mathbb R^2)}, \label{eq:general-multilinear-bound} \\
\|T_m(f, g)\|_{L^p(\mathbb R^2)} &\leq C \mathfrak{A} \|f\|_{L^p(\mathbb R^2)} \|g\|_{L^\infty(\mathbb R^2)} \qquad (r=2), \label{eq:bilinear-Lp-bound} \\
\|T_m(f, g, h)\|_{L^p(\mathbb R^2)} &\leq C \mathfrak{A} \|f\|_{L^p(\mathbb R^2)} \|g\|_{L^\infty(\mathbb R^2)} \|h\|_{L^\infty(\mathbb R^2)} \qquad (r=3). \label{eq:trilinear-Lp-bound}
\end{align}
\end{lemma}

\begin{proof}
The detailed proof is deferred to Appendix~\ref{app:Q-proof11}.
\end{proof}

\begin{lemma}[Estimates for the bilinear normal form]\label{lem:B}
	Let $0<\beta<\alpha<1$, and let $b$ and $\Bop$ be defined by
	\eqref{eq:b-symbol} and \eqref{eq:B-definition}, respectively.
	Writing $q=k+l$, we have
	\begin{equation}\label{eq:b-bounds}
		|b(k,l)|\lesssim
		\begin{cases}
			N^{1-2\alpha}L^{1-\alpha},
			&|k|\simeq L\ll |l|\simeq|q|\simeq N,\\
			N^{2-2\alpha}K^{-\alpha},
			&|k|\simeq|l|\simeq K\gg|q|\simeq N,\\
			N^{2-3\alpha},
			&|k|\simeq|l|\simeq|q|\simeq N.
		\end{cases}
	\end{equation}
	For $2\leq p\leq\infty$ and $j\geq0$,
	\begin{equation}\label{eq:B-dyadic}
		\|\Delta_j\Bop(\theta,\theta)\|_{L^p}
		\lesssim \|\theta\|_{L^\infty}\,2^{j(2-3\alpha)}
		\sum_{v\geq j-C}2^{-\alpha(v-j)}
		\|\widetilde\Delta_v\theta\|_{L^p},
	\end{equation}
	where $C$ depends only on the Littlewood--Paley cutoffs.
	Consequently, with $u_B=\Kop\Bop(\theta,\theta)$,
	\begin{equation}\label{eq:uB-holder}
		\|u_B\|_{B_{\infty,\infty}^{3\alpha-1}}
		\lesssim \|\theta\|_{L^\infty}^2.
	\end{equation}
	Moreover, for every $M\geq1$,
	\begin{equation}\label{eq:B-lowpass}
		\|S_M\Bop(\theta,\theta)\|_{L^p}
		\lesssim L_\alpha(M)
		\|\theta\|_{L^\infty}\|\theta\|_{L^p},
	\end{equation}
	where
	\begin{equation*}
		L_\alpha(M)=
		\begin{cases}
			M^{2-3\alpha},&\alpha<\frac{2}{3},\\
			1+\log(2+M),&\alpha=\frac{2}{3},\\
			1,&\alpha>\frac{2}{3}.
		\end{cases}
	\end{equation*}
\end{lemma}

\begin{proof}
	We first prove the estimates for smooth inputs with finite
	frequency truncations. The bounds obtained below are uniform
	in the truncation parameters.
	
	\medskip\noindent
	\emph{Localized multiplier estimates.}
	Recall that
	\[
	b(k,l)=\chi_{\mathrm{out}}(k+l)
	\frac{c_2(k,l)}{S(k,l)}
	g\!\left(\frac{D(k,l)}{S(k,l)}\right),
	\qquad g(z)=\frac{\eta(z)}{z},
	\]
	where $g$ is extended smoothly across the origin.
	We claim that, for $|\mu|+|\nu|\leq12$,
	\begin{equation}\label{eq:unified-symbol}
		\begin{aligned}
			|\partial_k^\mu\partial_l^\nu b(k,l)|
			&\lesssim_{\mu,\nu}
			N^{1-2\alpha-|\nu|}L^{1-\alpha-|\mu|},
			&&|k|\simeq L\ll|l|\simeq|k+l|\simeq N,\\
			|\partial_q^\mu\partial_k^\nu[b(k,q-k)]|
			&\lesssim_{\mu,\nu}
			N^{2-2\alpha-|\mu|}K^{-\alpha-|\nu|},
			&&|k|\simeq|q-k|\simeq K,\quad |q|\simeq N\lesssim K.
		\end{aligned}
	\end{equation}
	In the second line, differentiation in $k$ keeps $q$ fixed.
	The first estimate also holds with the two inputs interchanged.
	
	To verify the claim, first consider the low--high region.
	In the difference formula \eqref{asvbut897}, we use
	\[
	A(k+l)-A(l)=\int_0^1 k\cdot\nabla A(l+tk)\,dt.
	\]
	Together with the symbol estimates for $A$, this gives
	\[
	|\partial_k^\mu\partial_l^\nu c_2(k,l)|
	\lesssim_{\mu,\nu}
	L^{1-\alpha-|\mu|}N^{1-\alpha-|\nu|}.
	\]
	Indeed, differentiating the integral formula costs at most
	$L^{-1}$ in the low variable and $N^{-1}$ in the high variable.
	The other term in \eqref{asvbut897} satisfies the same bound
	directly from the estimates for $A$.
	
	For comparable inputs, use $(q,k)$ as coordinates and write
	$l=q-k$. The identities
	\[
	k^\perp\cdot l=k^\perp\cdot q,
	\qquad
	\frac1{|k|^2}-\frac1{|q-k|^2}
	=\frac{|q|^2-2q\cdot k}{|k|^2|q-k|^2}
	\]
	show that the two terms in \eqref{eq:c2-symbol} are bounded by
	$N^2K^{-2\alpha}$ and $N^{2-\alpha}K^{-\alpha}$, respectively.
	Since $N\lesssim K$, the first is bounded by the second.
	The same calculation after differentiation gives
	\[
	|\partial_q^\mu\partial_k^\nu[c_2(k,q-k)]|
	\lesssim_{\mu,\nu}
	N^{2-\alpha-|\mu|}K^{-\alpha-|\nu|}.
	\]

	It remains to include the denominator and the cutoffs.
	In the low--high coordinates, a derivative in $k$ or $l$
	costs at most $L^{-1}$ or $N^{-1}$, respectively, relative
	to the size of $S$. In the comparable-input coordinates,
	the corresponding costs are $N^{-1}$ in $q$ and $K^{-1}$
	in $k$. These statements hold for both $D$ and $S$.
	The chain rule therefore bounds the derivatives of
	$S^{-1}g\left(\frac{D}{S}\right)$ by $S^{-1}$ times the same scale factors.
	Since $S\geq|q|^\alpha$, and $\chi_{\mathrm{out}}$ has uniform
	scaled derivative bounds on its support, the product rule
	proves \eqref{eq:unified-symbol}.
	Taking $\mu=\nu=0$ also proves \eqref{eq:b-bounds}, with
	the balanced case obtained by setting $K\simeq N$.
	
	We apply Lemma~\ref{lem:kernel} to these localized symbols.
	For homogeneous input blocks $f_L=\dot\Delta_Lf$ and
	$g_K=\dot\Delta_Kg$, with $L\leq K$ and $N=2^j$, it gives
	\begin{equation}\label{eq:block-bilinear}
		\begin{aligned}
			\|\Delta_j\Bop(f_L,g_K)\|_{L^p}
			\lesssim\|f_L\|_{L^\infty}\|g_K\|_{L^p}
			\begin{cases}
				N^{1-2\alpha}L^{1-\alpha},&L\ll K\simeq N,\\
				N^{2-2\alpha}K^{-\alpha},&L\simeq K\gtrsim N.
			\end{cases}
		\end{aligned}
	\end{equation}
	The two norm placements may be interchanged.
	Here the frequency cutoffs are included before applying the
	kernel lemma: the first line uses coordinates $(k,l)$ at scales
	$(L,N)$, and the second uses $(k+l,k)$ at scales $(N,K)$.
	The finitely many output blocks in the transition region of
	$\chi_{\mathrm{out}}$ satisfy the same estimates.
	
	\smallskip\noindent
	\emph{Summation over the input frequencies.}
	The symmetry of $\Bop$ gives the Bony decomposition
	\begin{equation*}
		\begin{aligned}
			\Bop(\theta,\theta)
			= 2\sum_{v\geq-1}\sum_{-1\leq\ell\leq v-2}
			\Bop(\Delta_\ell\theta,\Delta_v\theta) +\sum_{v\geq-1}
			\Bop(\Delta_v\theta,\widetilde\Delta_v\theta),
		\end{aligned}
	\end{equation*}
	where $\widetilde\Delta_v=\Delta_{v-1}+\Delta_v+\Delta_{v+1}$
	and blocks with index below $-1$ are set equal to zero.
	Denote the $\Delta_j$ projections of the first and second
	sums by $I_j$ and $II_j$, respectively.
	
	For each $v$ in the first sum, combine the low blocks
	$\sum_{\ell=-1}^{v-2}\Delta_\ell\theta$ before resolving this
	low-frequency factor into homogeneous annuli. Each resulting
	block has $L^\infty$ norm bounded by $C\|\theta\|_{L^{\infty}}$.
	This also treats $\Delta_{-1}\theta$ without regarding its
	support as an annulus separated from zero.
	Only $|v-j|\leq C$ contributes, and the first line of
	\eqref{eq:block-bilinear} yields
	\begin{equation}\label{eq:Ij}
		\begin{aligned}
			\|I_j\|_{L^p}
			&\lesssim\|\theta\|_{L^{\infty}} N^{1-2\alpha}
			\sum_{|v-j|\leq C}\|\widetilde\Delta_v\theta\|_{L^p}
			\sum_{L\lesssim2^v}L^{1-\alpha}\\
			&\lesssim\|\theta\|_{L^{\infty}} N^{2-3\alpha}
			\sum_{|v-j|\leq C}\|\widetilde\Delta_v\theta\|_{L^p}.
		\end{aligned}
	\end{equation}
	The last inequality uses $1-\alpha>0$ and $2^v\simeq N$.
	In particular, the sum over arbitrarily small input scales
	is convergent.
	
	For $II_j$, the two input frequencies have size $K\simeq2^v$,
	and frequency support requires $v\geq j-C$.
	The second line of \eqref{eq:block-bilinear} gives
	\begin{equation}\label{eq:IIj}
		\begin{aligned}
			\|II_j\|_{L^p}
			&\lesssim\|\theta\|_{L^{\infty}}
			\sum_{v\geq j-C}N^{2-2\alpha}2^{-\alpha v}
			\|\widetilde\Delta_v\theta\|_{L^p}\\
			&=\|\theta\|_{L^{\infty}} N^{2-3\alpha}
			\sum_{v\geq j-C}2^{-\alpha(v-j)}
			\|\widetilde\Delta_v\theta\|_{L^p}.
		\end{aligned}
	\end{equation}
	Terms confined to the lowest input blocks vanish because of
	the fixed high-output cutoff. Combining \eqref{eq:Ij} and
	\eqref{eq:IIj}, and observing that $2^{-\alpha(v-j)}$ is
	comparable to one for $|v-j|\leq C$, proves \eqref{eq:B-dyadic}.
	
	\medskip\noindent
	\emph{The velocity and low-pass estimates.}
	Taking $p=\infty$ in \eqref{eq:B-dyadic} and summing the
	geometric tail, we obtain
	\[
	\|\Delta_j\Bop(\theta,\theta)\|_{L^{\infty}}
	\lesssim2^{j(2-3\alpha)}\|\theta\|_{L^{\infty}}^2.
	\]
	Since $\Kop$ has order $-1$ and the output of $\Bop$ is
	bounded away from zero,
	\[
	\|\Delta_j u_B\|_{L^{\infty}}
	\lesssim2^{-j(3\alpha-1)}\|\theta\|_{L^{\infty}}^2.
	\]
	The weighted supremum over $j$ proves \eqref{eq:uB-holder}.

	Finally, we consider the low-pass operator $S_M$, for example $M=2^{k}$.
It follows from \eqref{eq:B-dyadic} that
\begin{equation*}
\norm{\Delta_j\Bop(\theta,\theta)}{L^{p}}
\lesssim 2^{j(2-3\alpha)}\norm{\theta}{L^{\infty}}\norm{\theta}{L^{p}},
\end{equation*}
which yields
\begin{align*}
    \|S_M \Bop(\theta,\theta)\|_{L^{p}} &\leq \sum_{1 \leq j \leq k} \|\Delta_j \Bop\|_{L^{p}} \\
    &\lesssim \norm{\theta}{L^{\infty}}\norm{\theta}{L^{p}}\sum_{1 \leq j \leq k} 2^{j(2-3\alpha)}.
\end{align*}
It is not hard to check that
\begin{equation*}
 \sum_{1 \leq j \leq k} 2^{j(2-3\alpha)}\lesssim
 \begin{cases}
 2^{(2-3\alpha)k}=M^{2-3\alpha},&\alpha<\frac{2}{3},\\
 k\simeq 1+\log(2+M),&\alpha=\frac{2}{3},\\
 1,&\alpha>\frac{2}{3}.
 \end{cases}
\end{equation*}
Therefore, this proves \eqref{eq:B-lowpass}.
\end{proof}

\begin{lemma}[Resonant and low-output forcing]\label{lem:res}
For $2\leq p<\infty$ the low-output forcing always satisfies
\begin{equation}\label{eq:low-bound}
\norm{\Lambda^{-\frac{\alpha}{2}}F_{\mathrm{low}}}{L^p}
\lesssim (1+\norm{\theta}{L^\infty})\norm{\theta}{L^p}.
\end{equation}
If in addition
\begin{equation}\label{eq:res-condition}
\beta<\frac{2\alpha^2}{4-3\alpha},
\end{equation}
then it holds
\begin{equation}\label{eq:res-bound}
\norm{\Lambda^{-\frac{\alpha}{2}}\Cres(\theta,\theta)}{L^p}
\lesssim\norm{\theta}{L^\infty}\norm{\theta}{L^p}.
\end{equation}
\end{lemma}

\begin{proof}
Regarding $F_{\mathrm{low}} = (1-\chi_h(D))\partial_1\theta + \Clow(\theta,\theta)$, the output frequency is restricted by the low-pass cutoffs to $N \leq C_0$ for some fixed constant $C_0$.
First, we consider the bilinear term $\Clow(\theta,\theta)$. We split the interactions into High-High, Low-High, and Balanced frequency cases. Since the output $N$ is bounded, the summation structures differ from the high-output case.

\begin{itemize}
    \item \textbf{High-high interactions ($K \simeq L \geq c_0$ and $N \leq C_0$):} \\
    We use the High-High symbol bound multiplied by $N^{-\frac{\alpha}{2}}$:
    \[ N^{-\frac{\alpha}{2}}|c_2(k,l)| \lesssim N^{2-\frac{3\alpha}{2}}K^{-\alpha}. \]
    We sum this over all dyadic blocks $N \leq C_0$ and $K \geq c_0$. Because $K$ is a dyadic sequence, the infinite sum over $K$ is a convergent geometric series (since $\alpha > 0$):
    \[ \sum_{K \ge c_0} K^{-\alpha} < \infty. \]
    The sum over $N$ is merely a finite sum of strictly positive powers (since $2-\frac{3\alpha}{2}>0$). Therefore, the total operator norm multiplier $\mathfrak S_{\mathrm{hh}}$ is bounded by
    \[ \mathfrak S_{\mathrm{hh}} \lesssim \sum_{N\leq C_0} N^{2-\frac{3\alpha}{2}}\sum_{K\geq c_0}K^{-\alpha} < \infty. \]

    \item \textbf{Low-high interactions ($c_0 \lesssim N \simeq K \lesssim C_0$ and $L \ll N$):} \\
    We use the Low-High symbol bound multiplied by $N^{-\frac{\alpha}{2}}$:
    \[ N^{-\frac{\alpha}{2}}|c_2(k,l)| \lesssim N^{-\frac{\alpha}{2}} \left( N^{1-\alpha}L^{1-\alpha} \right) = N^{1-\frac{3\alpha}{2}}L^{1-\alpha}. \]
    For a fixed $N$, summing over $L \ll N$ yields a geometric series dominated by its largest term $L \simeq N$:
    \[ \sum_{L \ll N} L^{1-\alpha} \lesssim N^{1-\alpha} \quad \text{(since } 1-\alpha > 0). \]
    The outer sum over $N$ is bounded by a finite number of scales between $c_0$ and $C_0$:
    \[ \mathfrak S_{\mathrm{lh}} \lesssim \sum_{c_0\lesssim N\lesssim C_0} N^{1-\frac{3\alpha}{2}} \left( N^{1-\alpha} \right) < \infty. \]

    \item \textbf{Balanced interactions ($K \simeq L \simeq N \lesssim C_0$):} \\
    These occupy only a finite number of scales and are trivially bounded by a constant.
\end{itemize}
Summing the contributions from all three interaction regimes gives:
\[ \|\Lambda^{-\frac{\alpha}{2}}\Clow(\theta,\theta)\|_{L^{p}} \lesssim \|\theta\|_{L^{\infty}} \|\theta\|_{L^{p}}  (1 + \mathfrak S_{\mathrm{hh}} + \mathfrak S_{\mathrm{lh}}) \lesssim \|\theta\|_{L^{\infty}}  \|\theta\|_{L^{p}}. \]
In the Fourier domain, the operator $\Lambda^{-\frac{\alpha}{2}} \partial_1$ introduces a multiplier of magnitude $|\xi|^{1-\frac{\alpha}{2}}$. The smooth cutoff $(1-\chi_h)$ restricts this to dyadic scales $N \leq C_0$. Summing over these low frequencies yields
\begin{align*}
\|\Lambda^{-\frac{\alpha}{2}}(1-\chi_h(D))\partial_1\theta\|_{L^{p}}  &\lesssim \sum_{N\leq C_0} N^{1-\frac{\alpha}{2}} \|\Delta_N \theta\|_{L^{p}} \\& \lesssim \left( \sum_{N\leq C_0} N^{1-\frac{\alpha}{2}} \right) \|\theta\|_{L^{p}}\\&\lesssim \|\theta\|_{L^{p}}.
\end{align*}
Combining the bounds for $\Clow$ and the linear term yields \eqref{eq:low-bound}.

Now we focus on the resonant forcing bound \eqref{eq:res-bound}.
The symbol of $\Cres(\theta,\theta)$ is supported where the output frequency is high ($N \gtrsim 1$) and the frequencies are close to the resonant set $|k+l|^\alpha \approx |k|^\beta + |l|^\beta$. As shown in \eqref{eq:res-scales}, this support forces the maximum input frequency $K = \max(|k|, |l|)$ to satisfy
\[ N^\alpha \simeq K^\beta \implies K \simeq N^{\frac{\alpha}{\beta}}. \]
For the resonant contribution, the frequency output of the temperature pair is localized to $|k+l| \simeq N \ll K$. By the triangle inequality, this strictly requires $|k| \sim |l| \sim K$, so that High-Low interactions are kinematically excluded by frequency disjointness.
Moreover, the high-output cutoff ensures we are in the extreme High-High to Low interaction regime where $K \gg N$.
We apply the simplified symbol bound for $c_2(k,l)$ in this regime
\[ |c_2(k,l)| \lesssim N^{2-\alpha}K^{-\alpha}. \]
The fractional Laplacian $\Lambda^{-\frac{\alpha}{2}}$ applied to the output provides a multiplier of $N^{-\frac{\alpha}{2}}$. Multiplying the symbol by this factor yields
\begin{align}\label{asdvbnmt65}
 N^{-\frac{\alpha}{2}}|c_2(k,l)| \lesssim N^{2-\frac{3\alpha}{2}}K^{-\alpha}.
 \end{align}
 Substituting the precise resonance identity $K \simeq N^{\frac{\alpha}{\beta}}$ into the term gives
 \begin{align}\label{asdvbnmt66}
 N^{-\frac{\alpha}{2}}|c_2(k,l)| &\lesssim N^{2-\frac{3\alpha}{2}} \left(N^{\frac{\alpha}{\beta}}\right)^{-\alpha} = N^{2-\frac{3\alpha}{2}-\frac{\alpha^2}{\beta}}.
 \end{align}
Because the frequencies are strictly locked by the relation $K^\beta \simeq N^\alpha$, for any given output scale $N$, the resonance condition selects only a bounded, $O(1)$ number of dyadic input blocks $K$. Thus, applying Lemma \ref{lem:kernel} with inputs placed in $L^\infty$ and $L^p$ shows
\begin{align*}
\|\Delta_N \Lambda^{-\frac{\alpha}{2}}\Cres(\theta,\theta)\|_{L^{p}} &\lesssim N^{2-\frac{3\alpha}{2}-\frac{\alpha^2}{\beta}}\|\theta\|_{L^{\infty}}\|\widetilde\Delta_K \theta\|_{L^{p}} \\& \lesssim N^{2-\frac{3\alpha}{2}-\frac{\alpha^2}{\beta}}
\|\theta\|_{L^{\infty}}\|\theta\|_{L^{p}}
.
\end{align*}
Summing this over all high-output dyadic blocks $N \gtrsim 1$ yields
\begin{align*}
\|\Lambda^{-\frac{\alpha}{2}}\Cres(\theta,\theta)\|_{L^{p}}  \lesssim \|\theta\|_{L^{\infty}}\|\theta\|_{L^{p}} \sum_{N\gtrsim1}N^{2-\frac{3\alpha}{2}-\frac{\alpha^2}{\beta}}\lesssim \|\theta\|_{L^{\infty}}\|\theta\|_{L^{p}} ,
\end{align*}
due to \eqref{eq:res-condition}.
This proves \eqref{eq:res-bound}.
\end{proof}

\begin{lemma}[Temperature-regularity resonance estimate]
\label{lem:res-regularity}
For $2\leq p<\infty$, it holds
\begin{equation}\label{eq:res-Besov}
\norm{\Lambda^{-\frac{\alpha}{2}}\Cres(\theta,\theta)}{L^p}
\lesssim\norm{\theta}{L^{\infty}}\norm{\theta}{B_{p,1}^{\frac{(4-3\alpha)
\beta-2\alpha^{2}}{2\alpha}}}.
\end{equation}
For each high output block with scale $N$, we have
\begin{equation}\label{eq:res-physical-order}
\norm{\Delta_N\Cres(\theta,\theta)}{L^{\infty}}
\lesssim\norm{\theta}{L^{\infty}}^2N^{\frac{(2-\alpha)\beta-\alpha^{2}}{\beta}}.
\end{equation}
We mention that the above estimates \eqref{eq:res-Besov} and \eqref{eq:res-physical-order} require $0<\beta<\alpha<1$ and the fixed cutoffs,
but do not require \eqref{eq:res-condition}.
\end{lemma}

\begin{proof}
The operator $\Cres(\theta,\theta)$ restricts the frequency interactions to the resonant support, where the resonance function $D(k,l)$ is close to zero. On this support, the output frequency $N \simeq |k+l|$ and the maximum input frequency $K \simeq \max(|k|, |l|)$ satisfy the relation \eqref{eq:res-scales}, namely,
\begin{equation} \label{eq:res_relation}
    K^\beta \simeq N^\alpha.
\end{equation}
Because the cutoffs $\chi_{\mathrm{nr}}$ and $\chi_{\mathrm{out}}$ define the transition regions smoothly but with finite overlapping support, the relation \eqref{eq:res_relation} is not an exact equality, but a strict comparability at the dyadic block level.

Let $\mathcal{I}(N)$ denote the set of all dyadic input scales $K$ that can interact to produce a specific dyadic output scale $N$ within the resonant support:
\[ \mathcal{I}(N) = \{K : K^\beta \simeq N^\alpha\}. \]
Because $0 < \beta < \alpha < 1$ are fixed parameters, the multiplicities are uniformly bounded
\begin{equation} \label{eq:multiplicity}
    \sup_N \#\mathcal{I}(N) + \sup_K \#\{N : K \in \mathcal{I}(N)\} \leq C_{\alpha,\beta} < \infty.
\end{equation}
Let us show the proof of \eqref{eq:multiplicity}. Indeed, in the Littlewood-Paley decomposition framework, the frequency scales $N$ and $K$ are restricted to dyadic values. For example, let us write them as
\[ N = 2^i, \qquad K = 2^j, \]
where $i, j \in \mathbb{Z}$ represent the indices of the dyadic blocks.
The relation $K^\beta \simeq N^\alpha$ is enforced by the support of the non-resonant cutoff function $\chi_{\mathrm{nr}}$ and the high-output cutoff $\chi_{\mathrm{out}}$. There exist strict absolute constants $0 < c_1 < c_2 < \infty$ such that
\[ c_1 N^\alpha \leq K^\beta \leq c_2 N^\alpha. \]
Substituting the dyadic representations $N = 2^i$ and $K = 2^j$ into the inequality yields
\[ c_1 (2^i)^\alpha \leq (2^j)^\beta \leq c_2 (2^i)^\alpha,\]
which implies
\begin{align}\label{logawe1}
 \log_2(c_1) + i\alpha \leq j\beta \leq \log_2(c_2) + i\alpha.
 \end{align}

\vspace{1em}
\noindent \textbf{Step 1: Bounding $\sup_N \#\mathcal{I}(N)$} \\
To find the number of allowable input scales $K$ for a fixed output scale $N$, we fix the index $i$ and solve for the bounds on the index $j$. Dividing \eqref{logawe1} by $\beta > 0$, we obtain
\[ i \frac{\alpha}{\beta} + \frac{\log_2(c_1)}{\beta} \leq j \leq i \frac{\alpha}{\beta} + \frac{\log_2(c_2)}{\beta}. \]
This indicates that the integer index $j$ must lie within a closed interval of length
\[\Xi_{j}= \frac{\log_2(c_2) - \log_2(c_1)}{\beta} \equiv \frac{C}{\beta}. \]
Crucially, the length of this interval $\Xi_{j}$ is a constant completely independent of $i$ (i.e., independent of $N$). The maximum number of integers that can fit inside an interval of length $L$ is $\lfloor L \rfloor + 1$. Thus, we get
\[ \#\mathcal{I}(N) \leq \frac{C}{\beta} + 1. \]
Taking the supremum over all dyadic scales $N$, we conclude that $\sup_N \#\mathcal{I}(N) \le C_{\alpha,\beta}$.

\vspace{1em}
\noindent \textbf{Step 2: Bounding $\sup_K \#\{N : K \in \mathcal{I}(N)\}$} \\
Conversely, to find the number of allowable output scales $N$ for a fixed input scale $K$, we fix the index $j$ and solve for the bounds on the index $i$. Rearranging \eqref{logawe1} and dividing by $\alpha>0$, we also get
\[ j \frac{\beta}{\alpha} - \frac{\log_2(c_2)}{\alpha} \leq i \leq j \frac{\beta}{\alpha} - \frac{\log_2(c_1)}{\alpha}. \]
Similarly, the integer index $i$ must lie within a closed interval of length:
\[\Xi_{i} = \frac{\log_2(c_2) - \log_2(c_1)}{\alpha} \equiv \frac{C}{\alpha}. \]
Again, the length of this interval is independent of $j$ (i.e., independent of $K$). Therefore, the maximum number of such integers $i$ is bounded by $\frac{C}{\alpha} + 1$.
Taking the supremum over all dyadic scales $K$, we conclude that $\sup_K \#\{N : K \in \mathcal{I}(N)\} \le C_{\alpha,\beta}$.
Combining the above two results, we obtain that
\[ \sup_N \#\mathcal{I}(N) + \sup_K \#\{N : K \in \mathcal{I}(N)\} \leq \frac{C}{\beta} + \frac{C}{\alpha} + 2 \equiv C_{\alpha,\beta} < \infty. \]
This bounded multiplicity mathematically ensures that when executing double summations over interacting dyadic blocks $\sum_N \sum_{K \in \mathcal{I}(N)}$, the inner summation behaves as a sum over a finite, uniformly bounded number of terms, thereby preventing any divergence without requiring additional geometric decay factors.

Now let us deal with \eqref{eq:res-Besov}. It follows from \eqref{asdvbnmt65} that
\[ N^{-\frac{\alpha}{2}} |c_2(k,l)| \lesssim N^{2 - \frac{3\alpha}{2}} K^{-\alpha}. \]
To absorb the divergence into the input function $\theta$, we use  \eqref{eq:res_relation} to rewrite the $N$ dependence entirely in terms of $K$. Since $N \simeq K^{\frac{\beta}{\alpha}}$, we substitute $N$ to derive
\begin{align*}
    N^{-\frac{\alpha}{2}} |c_2(k,l)|  \lesssim \left( K^{\frac{\beta}{\alpha}} \right)^{2 - \frac{3\alpha}{2}} K^{-\alpha} = K^{\frac{(4-3\alpha)
\beta-2\alpha^{2}}{2\alpha}}.
\end{align*}
The normalized derivatives with respect to the variables $(k+l, k)$ at scales $(N, K)$ obey the same bounds due to the scale-invariant decay of the multiplier $A(\xi)$.
Applying Lemma \ref{lem:kernel}, for a fixed output block $N$, one has
\[ \|\Delta_N \Lambda^{-\frac{\alpha}{2}} \Cres(\theta,\theta)\|_{L^{p}} \lesssim \|\theta\|_{L^{\infty}} \sum_{K \in \mathcal{I}(N)} K^{\frac{(4-3\alpha)
\beta-2\alpha^{2}}{2\alpha}}\|\widetilde\Delta_K \theta\|_{L^{p}}. \]
We sum this bound over all high output scales $N \gtrsim 1$
\begin{align*}
    \|\Lambda^{-\frac{\alpha}{2}} \Cres(\theta,\theta)\|_{L^{p}} &\leq \sum_{N \gtrsim 1} \|\Delta_N \Lambda^{-\frac{\alpha}{2}} \Cres(\theta,\theta)\|_{L^{p}} \\
    &\lesssim \|\theta\|_{L^{\infty}} \sum_{N \gtrsim 1} \sum_{K \in \mathcal{I}(N)} K^{\frac{(4-3\alpha)
\beta-2\alpha^{2}}{2\alpha}} \|\widetilde\Delta_K \theta\|_{L^{p}}.
\end{align*}
Thanks to the bounded multiplicity property \eqref{eq:multiplicity}, we can interchange the sums and collapse the sum over $N$ to show
\begin{align*}
    \sum_{N \gtrsim 1} \sum_{K \in \mathcal{I}(N)} K^{\frac{(4-3\alpha)
\beta-2\alpha^{2}}{2\alpha}} \|\widetilde\Delta_K \theta\|_{L^{p}} &\leq C_{\alpha,\beta} \sum_{K \gtrsim 1} K^{\frac{(4-3\alpha)
\beta-2\alpha^{2}}{2\alpha}} \|\widetilde\Delta_K \theta\|_{L^{p}}.
\end{align*}
Therefore, one has
\begin{align*}
    \|\Lambda^{-\frac{\alpha}{2}} \Cres(\theta,\theta)\|_{L^{p}} &\lesssim \|\theta\|_{L^{\infty}} \sum_{K \gtrsim 1} K^{\frac{(4-3\alpha)
\beta-2\alpha^{2}}{2\alpha}} \|\widetilde\Delta_K \theta\|_{L^{p}} \\
    &\lesssim \|\theta\|_{L^{\infty}} \norm{\theta}{B_{p,1}^{\frac{(4-3\alpha)
\beta-2\alpha^{2}}{2\alpha}}}.
\end{align*}
The finite number of transition scales near the cutoff $K \sim 1$ can be trivially absorbed into the inhomogeneous Besov norm. This proves \eqref{eq:res-Besov}.

Finally, for the physical bound, we do not apply the transfer multiplier $N^{-\frac{\alpha}{2}}$. According to \eqref{asdvbnmt66}, we have
\begin{align*}
    |c_2(k,l)|  \lesssim N^{2-\alpha} \left( N^{\frac{\alpha}{\beta}} \right)^{-\alpha}
    \lesssim N^{\frac{(2-\alpha)\beta-\alpha^{2}}{\beta}}.
\end{align*}
Applying Lemma \ref{lem:kernel} with both inputs placed in the $L^\infty$ norm yields
\[ \|\Delta_N \Cres(\theta,\theta)\|_{L^{\infty}} \lesssim \|\theta\|_{L^{\infty}}^2 \sum_{K \in \mathcal{I}(N)} N^{\frac{(2-\alpha)\beta-\alpha^{2}}{\beta}}. \]
Since the cardinality of $\mathcal{I}(N)$ is bounded by a constant $C_{\alpha,\beta}$ independent of $N$, we derive
\[ \|\Delta_N \Cres(\theta,\theta)\|_{L^{\infty}} \lesssim \|\theta\|_{L^{\infty}}^2 N^{\frac{(2-\alpha)\beta-\alpha^{2}}{\beta}} \#\mathcal{I}(N) \lesssim \|\theta\|_{L^{\infty}}^2N^{\frac{(2-\alpha)\beta-\alpha^{2}}{\beta}}. \]
This proves \eqref{eq:res-physical-order} and concludes the proof of Lemma \ref{lem:res-regularity}.
\end{proof}

\subsection{The transport commutator of the normal form}

The next estimate includes low velocity frequencies and interactions
of two high-frequency inputs with a low-frequency output.
These configurations require the exact commutator cancellation.

\begin{lemma}[Trilinear transport commutator]\label{lem:Q}
Let $0<\beta<\alpha<1$ and $2\leq p<\infty$.
If $\frac{4}{7}<\alpha<\frac{4}{5}$, then
\begin{equation}\label{eq:Q-top}
 \|\Lambda^{-\frac{\alpha}{2}}\Qop(\Kop f;\theta,\theta)\|_{L^p}
 \lesssim \|\theta\|_{L^\infty}^2
 \sum_{M\geq1}M^{2-\frac{7\alpha}{2}}\|S_{c_1M}f\|_{L^p},
\end{equation}
where the sum is dyadic and $c_1$ is a fixed sufficiently large
constant. In particular,
\begin{equation}\label{eq:Q-mixed}
 \|\Lambda^{-\frac{\alpha}{2}}\Qop(u_H;\theta,\theta)\|_{L^p}
 \lesssim\|\theta\|_{L^\infty}^2\|H\|_{L^p},
\end{equation}
\begin{equation}\label{eq:Q-cubic}
 \|\Lambda^{-\frac{\alpha}{2}}\Qop(\Top\theta;\theta,\theta)\|_{L^p}
 \lesssim\|\theta\|_{L^\infty}^2
       \|\theta\|_{B_{p,1}^{\frac{3(2-3\alpha)}{2}}}.
\end{equation}
Finally, for $\alpha>\frac{8}{13}$, it holds
\begin{equation}\label{eq:Q-quartic}
 \|\Lambda^{-\frac{\alpha}{2}}\Qop(u_B;\theta,\theta)\|_{L^p}
 \lesssim\|\theta\|_{L^\infty}^3\|\theta\|_{L^p}.
\end{equation}
\end{lemma}

 \begin{proof}
	Recall  that
	\begin{equation*}
		\Qop(v;g,h)
		=\Bop(v\cdot\nabla g,h)+\Bop(g,v\cdot\nabla h)
		-v\cdot\nabla\Bop(g,h).
	\end{equation*}
	
	\smallskip\noindent
	\emph{1. The range $U\ll N$.}
	Let $\kappa$ be the kernel of a localized piece of $\Bop$ and
	$\mathcal Q_\kappa$ the corresponding commutator. The cutoffs
	include the original and shifted supports in each use below.
	Denoting $v_y=v(x-y)$, $g_y=g(x-y)$, and $h_z=h(x-z)$, we are able to conclude
	\begin{align*}
		\mathcal Q_\kappa(v;g,h)
		&=\iint\kappa\big[(v_y-v)\cdot(\nabla g)_y\,h_z
		+(v_z-v)\cdot(\nabla h)_z\,g_y\big]\,dy\,dz\\
		&=\iint\kappa(v_y-v_z)\cdot(\nabla g)_y\,h_z\,dy\,dz\\
		&\quad+\iint [(\nabla_y+\nabla_z)\kappa]\cdot(v_z-v)
		g_yh_z\,dy\,dz,
	\end{align*}
where in the second identity we have used $\nabla\cdot v=0$. For output
	$N$ and comparable inputs of size $K\gtrsim N$,
	\[
	(q,k)=(k+l,k),\qquad k\cdot y+l\cdot z=k\cdot(y-z)+q\cdot z.
	\]
	Rescaling \eqref{eq:unified-symbol} as in Lemma~\ref{lem:kernel}
	gives
\begin{align}\label{eq:Q-high-slot-moment}
\iint (L|y|+N|z|)|\kappa|\,dy\,dz
		 \lesssim N^{1-2\alpha}L^{1-\alpha},
		&&(L\ll N,\ \text{inputs }L,N),
	\end{align}
	\begin{align}\label{eq:Q-comparable-moment}
		\iint\big[K|y-z|\,|\kappa|
		+|z|\,|(\nabla_y+\nabla_z)\kappa|\big]\,dy\,dz
		 \lesssim N^{2-2\alpha}K^{-\alpha},
		&&(\text{inputs }K,K).
	\end{align}
	In \eqref{eq:Q-comparable-moment}, $(\nabla_y+\nabla_z)\kappa$ has symbol $iq$ times the localized symbol. Thus, making use of the direct estimate
	$\|v_y-v_z\|_p\leq |y-z|\|\nabla v\|_p$, we have
	\[
	\|\Delta_N\Qop(v;g,h)\|_p
	\lesssim \|\theta\|_{L^\infty}^2\|\nabla v\|_p
	\begin{cases}
		N^{1-2\alpha}L^{1-\alpha},
		& (L,N),\quad \supp\widehat v\subset\{|a|\ll L\},\\
		N^{2-2\alpha}K^{-\alpha},
		& (K,K),\quad \supp\widehat v\subset\{|a|\ll N\},
	\end{cases}
	\]
	when the two input block norms are bounded by $C\|\theta\|_{L^\infty}$.
	The corresponding temperature sums satisfy
	\[
	N^{1-2\alpha}\sum_{L\lesssim N}L^{1-\alpha}
	+N^{2-2\alpha}\sum_{K\gtrsim N}K^{-\alpha}
	\lesssim N^{2-3\alpha}.
	\]
	
	For the remaining pair $L\lesssim U\ll N$, we take
$$g= \sum_{L\lesssim U}\dot\Delta_L\theta,\quad h=\dot\Delta_N\theta.$$
 Then, it follows that
	\begin{align}
		\Qop(\dot\Delta_U v;g,h)
		&=\Bop(\operatorname{div}((\dot\Delta_U v)g),h)
		-(\dot\Delta_U v)\cdot\Bop(\nabla g,h)\notag\\
		&\quad+\Bop(g,(\dot\Delta_U v)\cdot\nabla h)
		-(\dot\Delta_U v)\cdot\Bop(g,\nabla h).
		\label{eq:Q-def-proof}
	\end{align}
If the two inputs have
frequencies $L\ll N$, then the localized bilinear operator has
coefficient $ N^{1-2\alpha}L^{1-\alpha}$,
Thus, for suitable localized inputs $a_L,b_N$, we have
\begin{equation}\label{eq:low-high-Bss}
\|\Bop(a_L,b_N)\|_{L^p}
\lesssim N^{1-2\alpha}L^{1-\alpha}
\|a_L\|_{L^p}\|b_N\|_{L^\infty},
\end{equation}
and the corresponding $L^\infty$ version is also available. We now pair the $P$-frequency output with the $N$-frequency input
$h$.
Applying \eqref{eq:low-high-Bss}, we obtain
\begin{align}\label{eusnv01}
\|\Delta_N\Bop(\operatorname{div}((\dot\Delta_U v)g),h)\|_{L^p}
&\lesssim\sum_{P\lesssim U}
N^{1-2\alpha}P^{1-\alpha}\|\dot\Delta_P \operatorname{div}((\dot\Delta_U v)g)\|_{L^p}
\|h\|_{L^\infty}\nonumber
\\
&\lesssim\sum_{P\lesssim U}
N^{1-2\alpha}P^{1-\alpha}P\|\dot\Delta_U v\|_{L^p}
\|g\|_{L^\infty}\|h\|_{L^\infty}\nonumber
\\
&\lesssim
\|\theta\|_{L^\infty}^{2}\|\dot\Delta_U v\|_{L^p}N^{1-2\alpha}
\sum_{P\lesssim U}
P^{2-\alpha}.
\end{align}
For each $L\ll N$, it follows from the low--high
estimate \eqref{eq:low-high-Bss} and Bernstein's inequality that
\[
\|\Bop(\nabla \dot\Delta_L\theta,h)\|_{L^\infty}
\lesssim
N^{1-2\alpha}L^{1-\alpha}
\|\nabla \dot\Delta_L\theta\|_{L^\infty}\|h\|_{L^\infty}\lesssim
N^{1-2\alpha}L^{2-\alpha}
\|\theta\|_{L^\infty}^2.
\]
Multiplying by $\dot\Delta_U v$ in $L^p$ and summing over $L$ gives
\begin{align}\label{eusnv02}
\|\dot\Delta_U v\cdot\Bop(\nabla g,h)\|_{L^p}
&\leq \|\dot\Delta_U v\|_{L^p}
\sum_{L\lesssim U}
\|\Bop(\nabla \dot\Delta_L\theta,h)\|_{L^\infty} \nonumber\\
&\lesssim
N^{1-2\alpha}\|\theta\|_{L^\infty}^2
\|\dot\Delta_U v\|_{L^p}
\sum_{L\lesssim U}L^{2-\alpha}.
\end{align}
We now combine
\[\Bop(g,(\dot\Delta_U v)\cdot\nabla h)
		-(\dot\Delta_U v)\cdot\Bop(g,\nabla h).
\]
For a fixed low input block $g_{L}\triangleq\dot\Delta_L\theta$, let $\kappa_L(y,z)$ denote
the kernel of the corresponding localized bilinear operator.
With the notation $g_{L,y}=g_L(x-y),\ v_{U,z}=v_U(x-z)$,
the difference has the kernel representation
\begin{align}\label{eq:high-slot-kernel}
I_L(x)
=\iint \kappa_L(y,z)g_L(x-y)
\big[(\dot\Delta_U v)(x-z)-(\dot\Delta_U v)(x)\big]\cdot
\nabla h(x-z)\,dy\,dz.
\end{align}
By the fundamental theorem of calculus,
\[
(\dot\Delta_U v)(x-z)-(\dot\Delta_U v)(x)
=-\int_0^1 z\cdot\nabla (\dot\Delta_U v)(x-tz)\,dt.
\]
Taking the $L^p$ norm and using translation invariance of that norm, we obtain
\begin{align*}
\|(\dot\Delta_U v)(\cdot-z)-(\dot\Delta_U v)\|_{L^p}
&\leq \int_0^1 |z|\|\nabla (\dot\Delta_U v)(\cdot-tz)\|_{L^p}\,dt\\
&=|z|\|\nabla (\dot\Delta_U v)\|_{L^p}\\
&\lesssim U|z|\|\dot\Delta_U v \|_{L^p}.
\end{align*}
Applying Minkowski's integral inequality to
\eqref{eq:high-slot-kernel} and using \eqref{eq:Q-high-slot-moment}, we thus have
\begin{align}\label{eusnv03}
\|I_L\|_{L^p}
&\lesssim
\|g_L\|_{L^\infty}\|\nabla h\|_{L^\infty}
\int\!\!\int |\kappa_L(y,z)|
\|(\dot\Delta_U v)(\cdot-z)-(\dot\Delta_U v)\|_{L^p}\,dy\,dz \nonumber \\
&\lesssim
\|\theta\|_{L^\infty} N\|h\|_{L^\infty}U\|\dot\Delta_U v\|_{L^p}
\int\!\!\int |z||\kappa_L(y,z)|\,dy\,dz \nonumber \\
& \lesssim
U\|\theta\|_{L^\infty}^2\|\dot\Delta_U v\|_{L^p}
\int\!\!\int N|z||\kappa_L(y,z)|\,dy\,dz\nonumber
\\
& \lesssim UN^{1-2\alpha}L^{1-\alpha}
\|\theta\|_{L^\infty}^2\|\dot\Delta_U v\|_{L^p}.
\end{align}
Combining \eqref{eusnv01}, \eqref{eusnv02} and
\eqref{eusnv03}, we find
	\begin{align*}
		\|\Delta_N\Qop(\dot\Delta_U v;g,h)\|_p
		&\lesssim \|\theta\|_{L^\infty}^2\|\dot\Delta_U v\|_p N^{1-2\alpha}
		\left(\sum_{P\lesssim U}P^{2-\alpha}
		+\sum_{L\lesssim U}L^{2-\alpha}
		+U\sum_{L\lesssim U}L^{1-\alpha}\right)\\
		&\lesssim \|\theta\|_{L^\infty}^2\|\dot\Delta_U v\|_p N^{1-2\alpha}U^{2-\alpha}\\
		&\lesssim \|\theta\|_{L^\infty}^2\|\dot\Delta_U v\|_p UN^{2-3\alpha}.
	\end{align*}
	This together with $\|\nabla \dot\Delta_U v\|_p\lesssim U\|\dot\Delta_U v\|_p$ gives
	\begin{equation}\label{eq:LP-Q-low}
		\|\Delta_N\Qop(\dot\Delta_U v;\theta,\theta)\|_p
		\lesssim \|\theta\|_{L^\infty}^2 UN^{2-3\alpha}\|\dot\Delta_U v\|_p,\qquad U\ll N .
	\end{equation}
	
	\smallskip\noindent
	\emph{2. The range $U\gtrsim N$.}
	If both temperature frequencies are $\lesssim U$, denote the
	product scale in a shifted term by $P$ and the other input
	scale by $L$. The coefficient from Lemma~\ref{lem:B},
	including the divergence factor $P$, is
	\[
	\begin{array}{c|c}
		P\simeq L\gtrsim N & N^{2-2\alpha}P^{1-\alpha}\\
		P\ll N\simeq L & N^{1-2\alpha}P^{2-\alpha}\\
		L\ll N\simeq P & N^{2-2\alpha}L^{1-\alpha}
	\end{array}
	\]
	times $\|\theta\|_{L^\infty}^2\|\dot\Delta_U v\|_p$. Consequently, the two shifted terms
	and the unshifted term have respective coefficients
	\begin{align*}
		&N^{2-2\alpha}\sum_{N\lesssim P\lesssim U}P^{1-\alpha}
		+N^{1-2\alpha}\sum_{P\lesssim N}P^{2-\alpha}
		+N^{2-2\alpha}\sum_{L\lesssim N}L^{1-\alpha}
		\lesssim N^{2-2\alpha}U^{1-\alpha},\\
		&\begin{cases}
			NU^{2-3\alpha}, & U\gg N \quad(\text{exterior divergence}),\\
			\displaystyle\sum_{P\lesssim U}P^{3-3\alpha}
			\lesssim NU^{2-3\alpha},
			& U\simeq N \quad(\text{internal derivative}).
		\end{cases}
	\end{align*}
	
	If one temperature frequency has size $K\gg U$, both do.
	Let $a,k,l$ be the input frequencies, $q=k+l$, $n=a+q$,
	and $F(q,k)=b(k,q-k)$. The symbol divided by its common
	factor $i$ is
	\begin{align*}
		&(\widehat v(a)\cdot k)[b(k+a,l)-b(k,l)]
		+(\widehat v(a)\cdot l)[b(k,l+a)-b(k,l)]\notag\\
		&\quad=(\widehat v(a)\cdot k)[F(n,k+a)-F(n,k)]
		+(\widehat v(a)\cdot n)[F(n,k)-F(q,k)] .
	\end{align*}
	Here $\widehat v(a)\cdot n=\widehat v(a)\cdot q$. In coordinates $(a,n,k)$
	at scales $(U,N,K)$, \eqref{eq:unified-symbol} gives
	\[
	F(n,k+a)-F(n,k)=\int_0^1 a\cdot\nabla_k F(n,k+ta)\,dt .
	\]
	The bounds for the three terms, after localization and
	with the factor $\|\theta\|_{L^\infty}^2\|\dot\Delta_U v\|_p$ suppressed, are
	\begin{align*}
		k[F(n,k+a)-F(n,k)]
		&:\quad K\,\left(\frac{U}{K}\right)\,N^{2-2\alpha}K^{-\alpha},\\
		nF(n,k)
		&:\quad N^{3-2\alpha}K^{-\alpha},\\
		nF(q,k)
		&:\quad
		\begin{cases}
			NU^{2-2\alpha}K^{-\alpha}, & U\gg N,\\
			\displaystyle\sum_{P\lesssim U}P^{3-2\alpha}K^{-\alpha}
			\lesssim NU^{2-2\alpha}K^{-\alpha}, & U\simeq N.
		\end{cases}
	\end{align*}
	For the first two lines the scaled derivatives satisfy the
	same bounds, so Lemma~\ref{lem:kernel} applies. In the last
	line use the bilinear output at scale $U$ or retain the
	internal derivative at scale $P$, respectively. Since
	$\sum_{K\gg U}K^{-\alpha}\lesssim U^{-\alpha}$, we obtain
	\begin{equation}\label{eq:LP-Q-high}
		\|\Delta_N\Qop(\dot\Delta_U v;\theta,\theta)\|_p
		\lesssim \|\theta\|_{L^\infty}^2\big(N^{2-2\alpha}U^{1-\alpha}
		+NU^{2-3\alpha}\big)\|\dot\Delta_U v\|_p,
		\qquad U\gtrsim N .
	\end{equation}
	The two block estimates hold for $N\geq1$,
	$\frac{1}{2}<\alpha<1$, and $1<p<\infty$.
	
	\smallskip\noindent
	\emph{3. Summation for $v=\Kop f$.}
	In the comparable-input region group $v=\Kop S_{cN}f$.
	For inputs $L\ll N$ and $N$ group
	$v=\Kop S_{cL}S_{cN}f$ first. By the Riesz transform bounds,
	\[
	\|\nabla\Kop S_{cN}f\|_p
	+\|\nabla\Kop S_{cL}S_{cN}f\|_p
	\lesssim\|S_{CN}f\|_p .
	\]
	For the remaining pairs, sum $L\lesssim U$ before $U\ll N$
	as in \eqref{eq:Q-def-proof}. Using the annular inverse bound
	$\|\Kop\dot\Delta_U f\|_p\lesssim U^{-1}\|\dot\Delta_U f\|_p$, we get
	\begin{align}
		&\|\Delta_N\Lambda^{-\frac{\alpha}{2}}
		\Qop(\Kop S_{cN}f;\theta,\theta)\|_p\notag\\
		& \lesssim \|\theta\|_{L^\infty}^2 N^{-\frac{\alpha}{2}}
		\left[N^{2-3\alpha}\|S_{CN}f\|_p
		+N^{1-2\alpha}\sum_{U\ll N}
		U^{1-\alpha}\|\dot\Delta_U S_{cN}f\|_p\right]\notag\\
		& \lesssim \|\theta\|_{L^\infty}^2N^{2-\frac{7\alpha}{2}}\|S_{CN}f\|_p .
		\label{eq:LP-Q-grouped}
	\end{align}
	On the complementary velocity support, \eqref{eq:LP-Q-high}
	gives the coefficient
	\[
	W_{N,U}=N^{2-\frac{5\alpha}{2}}U^{-\alpha}
	+N^{1-\frac{\alpha}{2}}U^{1-3\alpha},
	\qquad U\gtrsim N\geq1,
	\]
	on $\|\theta\|_{L^\infty}^2\|\dot\Delta_U f\|_p$. The cutoffs are uniformly
	bounded on these annuli. The sums in either order satisfy
	\[
	\sum_{1\leq N\lesssim U}W_{N,U}\lesssim U^{2-\frac{7\alpha}{2}}
	\quad(\alpha<\frac{4}{5}),\qquad
	\sum_{U\gtrsim N}W_{N,U}\lesssim N^{2-\frac{7\alpha}{2}}
	\quad(\frac{1}{2}<\alpha<1).
	\]
	For homogeneous outputs $N<1$, the high-output cutoff
	annihilates the shifted terms. The unshifted term has internal
	bilinear output of size $U\gtrsim1$. The exterior divergence gives
	\begin{align*}
		\sum_{N<1}\|\dot\Delta_N\Lambda^{-\frac{\alpha}{2}}
		\Qop(\Kop f;\theta,\theta)\|_p
		&\lesssim \|\theta\|_{L^\infty}^2\sum_{N<1}N^{1-\frac{\alpha}{2}}
		\sum_{U\gtrsim1}U^{1-3\alpha}\|\dot\Delta_U f\|_p\\
		&\lesssim \|\theta\|_{L^\infty}^2\sum_{U\gtrsim1}U^{1-3\alpha}\|\dot\Delta_U f\|_p .
	\end{align*}
	For $\frac{4}{7}<\alpha<\frac{4}{5}$, we have
	\[
	U^{1-3\alpha}\lesssim U^{2-\frac{7\alpha}{2}},\qquad
	\|\dot\Delta_U f\|_p\lesssim\|S_{c_1U}f\|_p .
	\]
	Then \eqref{eq:LP-Q-grouped} and the sum first in $N$ give
	\eqref{eq:Q-top}. Summing first in $U$ instead gives
	\[
	\|\Lambda^{-\frac{\alpha}{2}}\Qop(\Kop f;\theta,\theta)\|_p
	\lesssim \|\theta\|_{L^\infty}^2\|f\|_p,\qquad \frac{4}{7}<\alpha<1 ,
	\]
	since $2-\frac{7\alpha}{2}<0$ and $1-3\alpha<0$.
	
	\smallskip\noindent
	\emph{4. Substitutions.}
	Taking $f=H$ gives \eqref{eq:Q-mixed}. For $f= R_{\alpha} \theta$
	and $\frac{4}{7}<\alpha<\frac{2}{3}$,
	\begin{align*}
		\sum_{M\geq1}M^{2-\frac{7\alpha}{2}}\|S_{c_1M} R_{\alpha} \theta\|_p
		&\lesssim\sum_{K\gtrsim1}K^{1-\alpha}\|\dot\Delta_K\theta\|_p
		\sum_{M\gtrsim K}M^{2-\frac{7\alpha}{2}}\\
		 &\lesssim\sum_{K\gtrsim1}K^{3-\frac{9\alpha}{2}}\|\dot\Delta_K\theta\|_p
		\lesssim\|\theta\|_{B_{p,1}^{\frac{3(2-3\alpha)}{2}}},
	\end{align*}
	which is \eqref{eq:Q-cubic}. For $f=\Bop(\theta,\theta)$,
	\eqref{eq:B-lowpass} and \eqref{eq:Q-top} give
	\[
	\|\Lambda^{-\frac{\alpha}{2}}\Qop(u_B;\theta,\theta)\|_p
	\lesssim \|\theta\|_{L^\infty}^3\|\theta\|_p
	\begin{cases}
		\displaystyle\sum_{M\geq1}M^{\frac{8-13\alpha}{2}},
		&\frac{8}{13}<\alpha<\frac{2}{3},\\
		\displaystyle\sum_{M\geq1}M^{-\frac{1}{3}}(1+\log(2+M)),
		&\alpha=\frac{2}{3}.
	\end{cases}
	\]
	For $\frac{2}{3}<\alpha<1$, \eqref{eq:B-dyadic} yields
	\[
	\|\Bop(\theta,\theta)\|_p
	\lesssim \|\theta\|_{L^\infty}\|\theta\|_p\sum_{N\geq1}N^{2-3\alpha}
	\lesssim \|\theta\|_{L^\infty}\|\theta\|_p ,
	\]
	and the unweighted $L^p$ bound gives \eqref{eq:Q-quartic}.
	The summable bounds and the $L^p$ convergence of the grouped
	Riesz terms allow removal of the truncations.
\end{proof}

\begin{lemma}[The cubic commutator involving $u_B$]\label{lem:other-cubic}
If $\frac{4}{7}<\alpha<1$ and $2\leq p<\infty$, then it holds
\begin{equation}\label{eq:other-cubic}
\norm{\Lambda^{-\frac{\alpha}{2}}[ R_{\alpha} ,u_B\cdot\nabla]\theta}{L^{p}}
\lesssim\norm{\theta}{L^{\infty}}^2\norm{\theta}{B_{p,1}^{\frac{3(2-3\alpha)}{2}}}.
\end{equation}
\end{lemma}

\begin{proof}
The symbol of the operator $ R_{\alpha} $ in the frequency space is given by
\begin{equation*}
    m(\xi) = \chi_h(\xi) \frac{i\xi_1}{|\xi|^\alpha - |\xi|^\beta}
\end{equation*}
Since $\chi_h(\xi)$ vanishes for $|\xi| \le R_0$, and for $|\xi| \ge R_0 > 1$, we can use the condition $R_0^{\beta-\alpha} \le \frac{1}{2}$ to obtain
\begin{equation*}
    |\xi|^\alpha - |\xi|^\beta = |\xi|^\alpha(1 - |\xi|^{\beta-\alpha}) \ge |\xi|^\alpha(1 - R_0^{\beta-\alpha}) \ge \frac{1}{2}|\xi|^\alpha
\end{equation*}
This implies that $ R_{\alpha} $ behaves as a standard Fourier multiplier of order $1-\alpha$. For a function localized at frequency $|\xi| \sim 2^j$, we have
\begin{itemize}
    \item $\| R_{\alpha}  \Delta_j f\|_{L^p} \lesssim 2^{(1-\alpha)j} \|\Delta_j f\|_{L^p}$
    \item $\|\nabla  R_{\alpha}  \Delta_j f\|_{L^p} \lesssim 2^{(2-\alpha)j} \|\Delta_j f\|_{L^p}$
    \item Let $K_j$ be the localized integral kernel of $ R_{\alpha}  \Delta_j$. The first moment of $K_j$ satisfies
        $$\||x| K_j(x)\|_{L^1} \lesssim 2^{-j} 2^{(1-\alpha)j} = 2^{-\alpha j}.$$
\end{itemize}
Let $X \triangleq \Lambda^{-\frac{\alpha}{2}} [ R_{\alpha} , u_B \cdot \nabla]\theta$. We will bound the sum of the $L^p$ norms of its dyadic blocks $\sum_{k} \|\Delta_k X\|_{L^p}$. Since $B_{p,1}^0 \hookrightarrow L^p$, this sum directly bounds $\|X\|_{L^p}$.
Recall that $u_B = \Kop \Bop(\theta,\theta) = \nabla^\perp \Delta^{-1} \Bop(\theta,\theta)$. From Lemma \ref{lem:B}, the bilinear normal form satisfies $\|\Delta_m \Bop(\theta,\theta)\|_{L^\infty} \lesssim 2^{m(2-3\alpha)} \|\theta\|_{L^\infty}^2$. Since $\Kop$ is an operator of order $-1$, we obtain
\begin{equation*}
    \|\Delta_m u_B\|_{L^\infty} \approx 2^{-m} \|\Delta_m \Bop(\theta,\theta)\|_{L^\infty} \lesssim 2^{m(1-3\alpha)} \|\theta\|_{L^\infty}^2.
\end{equation*}
We define the weighted sequence for the temperature blocks
\begin{equation*}
    \Theta_j \triangleq 2^{j\frac{3(2-3\alpha)}{2}} \|\Delta_j \theta\|_{L^p}.
\end{equation*}
By definition, its $\ell^1$ sum is exactly the target norm: $\sum_{j} \Theta_j = \|\theta\|_{B_{p,1}^{\frac{3(2-3\alpha)}{2}}}$.
Applying Bony's paraproduct decomposition to the product $u_B \cdot \nabla \theta$, projecting onto $\Delta_k$, and including the fractional differentiation $\Lambda^{-\frac{\alpha}{2}}$ and the commutator with $ R_{\alpha} $, we split the term into three dyadic components: $X_k \triangleq \Delta_k X = I_{1k} + I_{2k} + I_{3k}$.

For the term $I_{1k}$, $j \sim k$ ($|j-k| \le 4$). By the standard commutator estimate utilizing the first moment of the localized kernel of $ R_{\alpha} $, we get
\begin{equation*}
    \|[ R_{\alpha}  \Delta_j, S_{j-1} u_B \cdot \nabla] \Delta_j \theta\|_{L^p} \lesssim 2^{-\alpha j} \|\nabla S_{j-1} u_B\|_{L^\infty} \|\nabla \Delta_j \theta\|_{L^p}.
\end{equation*}
Since $\alpha<\frac{2}{3}$, the low-frequency sum is dominated by its highest block
\begin{equation*}
    \|\nabla S_{j-1} u_B\|_{L^\infty} \lesssim \sum_{m \le j-2} 2^m \|\Delta_m u_B\|_{L^\infty} \lesssim \sum_{m \le j-2} 2^{m(2-3\alpha)} \|\theta\|_{L^\infty}^2 \approx 2^{j(2-3\alpha)} \|\theta\|_{L^\infty}^2.
\end{equation*}
Applying the multiplier $\Lambda^{-\frac{\alpha}{2}}$ provides a scaling of $2^{-k\frac{\alpha}{2}}$. We thus have
\begin{align*}
    \|I_{1k}\|_{L^p} &\lesssim 2^{-k\frac{\alpha}{2}} 2^{-\alpha k} \big( 2^{k(2-3\alpha)} \|\theta\|_{L^\infty}^2 \big) \big( 2^k \|\Delta_k \theta\|_{L^p} \big) \\
    &= \|\theta\|_{L^\infty}^2 2^{k(3 - \frac{9\alpha}{2})} \|\Delta_k \theta\|_{L^p} \\&= \|\theta\|_{L^\infty}^2 \Theta_k.
\end{align*}
Summing over $k$ yields
$$\sum_{k} \|I_{1k}\|_{L^p} \lesssim \|\theta\|_{L^\infty}^2 \|\theta\|_{B_{p,1}^{\frac{3(2-3\alpha)}{2}}}.$$

The term $I_{2k}$ can be estimated via the triangle inequality without the commutator structure
\begin{equation*}
   I_{2k}= \Delta_k \Lambda^{-\frac{\alpha}{2}} \sum_{|j-k| \le 4} \Big( R_{\alpha} (\Delta_j u_B \cdot \nabla S_{j-1}\theta) - \Delta_j u_B \cdot \nabla  R_{\alpha} (S_{j-1}\theta) \Big)\triangleq I_{2k}^{1}+I_{2k}^{2}.
\end{equation*}
The operator $\Lambda^{-\frac{\alpha}{2}}  R_{\alpha} $ behaves as an operator of order $1-\frac{3\alpha}{2}$. Its bound is given by
\begin{align*}
\|I_{2k}^{1}\|_{L^{p}}&\lesssim
    2^{k(1-\frac{3\alpha}{2})} \|\Delta_k u_B\|_{L^\infty} \sum_{m \le k-2} 2^m \|\Delta_m \theta\|_{L^p}\\& \lesssim 2^{k(1-\frac{3\alpha}{2})} \big( 2^{k(1-3\alpha)} \|\theta\|_{L^\infty}^2 \big) \sum_{m \le k-2} 2^{-m(2-\frac{9\alpha}{2})} 2^{m\frac{3(2-3\alpha)}{2}} \|\Delta_m \theta\|_{L^p} \\& \lesssim \|\theta\|_{L^\infty}^2 \sum_{m \le k-2} 2^{-(m-k)(2-\frac{9\alpha}{2})} \Theta_m,
\end{align*}
which along with $\alpha>\frac{4}{9}$ yields
$$\sum_{k} \|I_{2k}^{1}\|_{L^p} \lesssim \|\theta\|_{L^\infty}^2 \|\theta\|_{B_{p,1}^{\frac{3(2-3\alpha)}{2}}}.$$
The operator $\nabla  R_{\alpha} $ behaves as an operator of order $2-\alpha$. Multiplying by $2^{-k\frac{\alpha}{2}} \|\Delta_k u_B\|_{L^\infty}$, we get
\begin{align*}
\|I_{2k}^{2}\|_{L^{p}}&\lesssim
    2^{-k\frac{\alpha}{2}} \big( 2^{k(1-3\alpha)} \|\theta\|_{L^\infty}^2 \big) \sum_{m \le k-2} 2^{(2-\alpha)m} \|\Delta_m \theta\|_{L^p}\\& = \|\theta\|_{L^\infty}^2 2^{k(1-\frac{7\alpha}{2})}\sum_{m \le k-2}  2^{-m(1-\frac{7\alpha}{2})} 2^{m\frac{3(2-3\alpha)}{2}} \|\Delta_m \theta\|_{L^p} \\& = \|\theta\|_{L^\infty}^2  \sum_{m \le k-2}  2^{-(m-k)(1-\frac{7\alpha}{2})} \Theta_m,
\end{align*}
which along with $\alpha>\frac{2}{7}$ yields
$$\sum_{k} \|I_{2k}^{2}\|_{L^p} \lesssim \|\theta\|_{L^\infty}^2 \|\theta\|_{B_{p,1}^{\frac{3(2-3\alpha)}{2}}}.$$

For the term $I_{3k}$, the output frequency is $k$, enforcing $j \ge k-2$. Since $u_B$ is divergence-free, we pull the divergence operator outside
\begin{equation*}
    I_{3k}= \Delta_k \Lambda^{-\frac{\alpha}{2}} \sum_{j \ge k-2} \Big(  R_{\alpha}  \nabla \cdot (\Delta_j u_B \widetilde{\Delta}_j \theta) - \nabla \cdot (\Delta_j u_B  R_{\alpha} (\widetilde{\Delta}_j \theta)) \Big)\triangleq I_{3k}^{1}+I_{3k}^{2}.
\end{equation*}
Regarding the term $I_{3k}^{1}$, the operator $\Delta_k \Lambda^{-\frac{\alpha}{2}} R_{\alpha}  \nabla \cdot$ has a total order of $2-\frac{3\alpha}{2}$. Thus we have
\begin{align*}
\|I_{3k}^{1}\|_{L^{p}}&\lesssim
    \sum_{j \ge k-2} 2^{k(2-\frac{3\alpha}{2})} \|\Delta_j u_B\|_{L^\infty} \|\widetilde{\Delta}_j \theta\|_{L^p} \\&\lesssim \sum_{j \ge k-2} 2^{k(2-\frac{3\alpha}{2})} \big( 2^{j(1-3\alpha)} \|\theta\|_{L^\infty}^2 \big) 2^{-j\frac{3(2-3\alpha)}{2}} \Theta_j \\
    &= \|\theta\|_{L^\infty}^2 \sum_{j \ge k-2} 2^{(k-j)(2-\frac{3\alpha}{2})} \Theta_j.
\end{align*}
This together with $\alpha<\frac{4}{3}$ implies
$$\sum_{k} \|I_{3k}^{1}\|_{L^p} \lesssim \|\theta\|_{L^\infty}^2 \|\theta\|_{B_{p,1}^{\frac{3(2-3\alpha)}{2}}}.$$
The operator $\Delta_k \Lambda^{-\frac{\alpha}{2}} \nabla \cdot$ behaves as an operator of order $1-\frac{\alpha}{2}$. We thus derive
\begin{align*}
\|I_{3k}^{2}\|_{L^{p}}&\lesssim
    \sum_{j \ge k-2} 2^{k(1-\frac{\alpha}{2})} \|\Delta_j u_B\|_{L^\infty} \| R_{\alpha} (\widetilde{\Delta}_j \theta)\|_{L^p} \\&\lesssim \sum_{j \ge k-2} 2^{k(1-\frac{\alpha}{2})} \big( 2^{j(1-3\alpha)} \|\theta\|_{L^\infty}^2 \big) 2^{j(1-\alpha)}2^{-j\frac{3(2-3\alpha)}{2}} \Theta_j \\
    &= \|\theta\|_{L^\infty}^2 \sum_{j \ge k-2} 2^{(k-j)(1-\frac{\alpha}{2})} \Theta_j,
\end{align*}
which along with $\alpha<2$ yields
$$\sum_{k} \|I_{3k}^{2}\|_{L^p} \lesssim \|\theta\|_{L^\infty}^2 \|\theta\|_{B_{p,1}^{\frac{3(2-3\alpha)}{2}}}.$$
Combining the sum estimates for $I_{1k}, I_{2k}$, and $I_{3k}$, we conclude that
\begin{equation*}
    \|\Lambda^{-\frac{\alpha}{2}}[ R_{\alpha} ,u_B\cdot\nabla]\theta\|_{L^{p}} \lesssim \|\theta\|_{L^\infty}^2 \|\theta\|_{B_{p,1}^{\frac{3(2-3\alpha)}{2}}}.
\end{equation*}
This completes the proof of the lemma.
\end{proof}

\subsection{The mixed commutator}
\begin{lemma}[Mixed commutator]\label{lem:CH}
Let $0<\beta<\alpha<1$. If $2\leq p<\infty$ and
$\frac{1}{p}=\frac{1}{q}+\frac{1}{b}$ with $1<q,b<\infty$, then it holds
\begin{equation}\label{eq:CH}
\norm{\Lambda^{-\frac{\alpha}{2}}[ R_{\alpha} ,\Kop H\cdot\nabla]\theta}{L^{p}}
\lesssim\norm{H}{L^{q}}\norm{\theta}{B_{b,2}^{\frac{2-3\alpha}{2}}}.
\end{equation}
\end{lemma}

\begin{proof}
Let $u_H = \Kop H = \nabla^\perp \Delta^{-1} H$. We aim to estimate the $L^p$ norm of the operator $X \triangleq \Lambda^{-\frac{\alpha}{2}} [ R_{\alpha} , u_H \cdot \nabla]\theta$ utilizing Bony's paraproduct decomposition and Littlewood-Paley theory. Since $u_H = \nabla^\perp \Delta^{-1} H$, its gradient $\nabla u_H$ is given by a matrix of second-order Riesz transforms acting on $H$. By the $L^q$ boundedness of Riesz transforms for $1 < q < \infty$, we have for the low-frequency cut-off:
\begin{equation*}
    \|\nabla S_{j-1} u_H\|_{L^q} \lesssim \|S_{j-1} H\|_{L^q} \lesssim \|H\|_{L^q}
\end{equation*}
Additionally, the order of $\Kop$ implies $\|\Delta_j u_H\|_{L^q} \approx 2^{-j} \|\Delta_j H\|_{L^q} \le 2^{-j} \|H\|_{L^q}$.
As established previously, $ R_{\alpha} $ acts as a multiplier of order $1-\alpha$. Let $K_j$ be the integral kernel of $ R_{\alpha}  \Delta_j$, which satisfies $\||x| K_j\|_{L^1} \lesssim 2^{-\alpha j}$.

For simplicity, we define the weighted sequence $\Theta_j = 2^{j\frac{2-3\alpha}{2}} \|\Delta_j \theta\|_{L^b}$. By definition, the $\ell^2$ norm of this sequence corresponds to the Besov norm: $\|(\Theta_j)\|_{\ell^2} \approx \|\theta\|_{B_{b,2}^{\frac{2-3\alpha}{2}}}$.

We decompose the product $u_H \cdot \nabla \theta$ using Bony's paraproduct
\begin{equation*}
    u_H \cdot \nabla \theta = \sum_{j} S_{j-1}u_H \cdot \nabla \Delta_j \theta + \sum_{j} \Delta_j u_H \cdot \nabla S_{j-1}\theta + \sum_{j} \Delta_j u_H \cdot \nabla \widetilde{\Delta}_j \theta
\end{equation*}
Applying the operator $\Delta_k \Lambda^{-\frac{\alpha}{2}}$ and taking the commutator with $ R_{\alpha} $, we split the term into three parts: $X_k \triangleq \Delta_k X = I_{1k} + I_{2k} + I_{3k}$. We apply H\"{o}lder's inequality in the form $L^q \times L^b \to L^p$, owing to the assumption $\frac{1}{p} = \frac{1}{q} + \frac{1}{b}$.

For the low-high part, we have $j \sim k$ ($|j-k| \le 4$). By the standard commutator estimate utilizing the first moment of the kernel, it gives
\begin{align*}
    \|[ R_{\alpha}  \Delta_j, S_{j-1} u_H \cdot \nabla] \Delta_j \theta\|_{L^p} &\lesssim \||x|K_j\|_{L^1} \|\nabla S_{j-1} u_H\|_{L^q} \|\nabla \Delta_j \theta\|_{L^b} \\
    &\lesssim 2^{-\alpha j} \|H\|_{L^q} \big( 2^j \|\Delta_j \theta\|_{L^b} \big) \\&= 2^{(1-\alpha)j} \|H\|_{L^q} \|\Delta_j \theta\|_{L^b}
\end{align*}
Applying $\Delta_k \Lambda^{-\frac{\alpha}{2}}$ provides a multiplier scaling of $2^{-k\frac{\alpha}{2}}$. Thus:
\begin{equation*}
    \|I_{1k}\|_{L^p} \lesssim 2^{-k\frac{\alpha}{2}} 2^{(1-\alpha)k} \|H\|_{L^q} \|\Delta_k \theta\|_{L^b} = 2^{k\frac{2-3\alpha}{2}} \|H\|_{L^q} \|\Delta_k \theta\|_{L^b} = \|H\|_{L^q} \Theta_k
\end{equation*}
This implies
$$\left( \sum_{k} \|I_{1k}\|_{L^p}^2 \right)^{\frac{1}{2}} \lesssim \|H\|_{L^q} \|\theta\|_{B_{b,2}^{\frac{2-3\alpha}{2}}}$$

We bound the two components of the commutator $I_{2k}$ directly using the triangle inequality
\begin{equation*}
    I_{2k} = \Delta_k \Lambda^{-\frac{\alpha}{2}} \sum_{|j-k| \le 4} \Big( R_{\alpha} (\Delta_j u_H \cdot \nabla S_{j-1}\theta) - \Delta_j u_H \cdot \nabla  R_{\alpha} (S_{j-1}\theta) \Big)\triangleq I_{2k}^{1}+I_{2k}^{2}.
\end{equation*}
For the term $I_{2k}^{1}$, since $ R_{\alpha} $ behaves as an operator of order $1-\alpha$, applying $\Lambda^{-\frac{\alpha}{2}}$ yields a total weight of $2^{k(1-\frac{3\alpha}{2})}$. Using $\|\Delta_k u_H\|_{L^q} \lesssim 2^{-k} \|H\|_{L^q}$ implies
\begin{align*}
    \|I_{2k}^{1}\|_{L^{p}}&\lesssim
    2^{k(1-\frac{3\alpha}{2})} \|\Delta_k u_H\|_{L^q} \|\nabla S_{k-1}\theta\|_{L^b} \\&\lesssim 2^{-k\frac{3\alpha}{2}} \|H\|_{L^q} \sum_{m \le k-2} 2^m \|\Delta_m \theta\|_{L^b} \\
    &= \|H\|_{L^q} \sum_{m \le k-2} 2^{-(k-m)\frac{3\alpha}{2}} \Theta_m
\end{align*}
For the second term, $\nabla  R_{\alpha} $ behaves as an operator of order $2-\alpha$. Multiplying by $2^{-k\frac{\alpha}{2}} \|\Delta_k u_H\|_{L^q} \lesssim 2^{-k(1+\frac{\alpha}{2})} \|H\|_{L^q}$, we have
\begin{align*}
    \|I_{2k}^{2}\|_{L^{p}}&\lesssim
    2^{-k(1+\frac{\alpha}{2})} \|H\|_{L^q} \sum_{m \le k-2} 2^{(2-\alpha)m} \|\Delta_m \theta\|_{L^b} \\&= \|H\|_{L^q} \sum_{m \le k-2} 2^{-(k-m)(1+\frac{\alpha}{2})} \Theta_m.
\end{align*}
Since $\alpha > 0$, both terms are convolutions of $\Theta_m$ with geometrically decaying sequences. As a result, one has
$$\left( \sum_{k} \|I_{2k}\|_{L^p}^2 \right)^{\frac{1}{2}} \lesssim\left( \sum_{k} \|I_{2k}^{1}\|_{L^p}^2 \right)^{\frac{1}{2}}+\left( \sum_{k} \|I_{2k}^{2}\|_{L^p}^2 \right)^{\frac{1}{2}} \lesssim \|H\|_{L^q} \|\theta\|_{B_{b,2}^{\frac{2-3\alpha}{2}}}.$$

For the term $I_{3k}$, the output frequency is $k$, enforcing $j \ge k-2$. Since $u_H$ is divergence-free, we pull the divergence operator outside:
\begin{equation*}
    I_{3k} = \Delta_k \Lambda^{-\frac{\alpha}{2}} \sum_{j \ge k-2} \Big(  R_{\alpha}  \nabla \cdot (\Delta_j u_H \widetilde{\Delta}_j \theta) - \nabla \cdot (\Delta_j u_H  R_{\alpha} (\widetilde{\Delta}_j \theta)) \Big)\triangleq I_{3k}^{1} +I_{3k}^{2}.
\end{equation*}
The operator $\Delta_k \Lambda^{-\frac{\alpha}{2}}  R_{\alpha}  \nabla \cdot$ has a total order of $2-\frac{3\alpha}{2}$. Thus, $I_{3k}^{1}$ is bounded by
\begin{align*}
\|I_{3k}^{1}\|_{L^{p}}&\lesssim
    \sum_{j \ge k-2} 2^{k(2-\frac{3\alpha}{2})} \|\Delta_j u_H\|_{L^q} \|\widetilde{\Delta}_j \theta\|_{L^b} \\ &\lesssim \sum_{j \ge k-2} 2^{k(2-\frac{3\alpha}{2})} 2^{-j} \|H\|_{L^q} \|\widetilde{\Delta}_j \theta\|_{L^b} \\
    &= \|H\|_{L^q} \sum_{j \ge k-2} 2^{-(j-k)(2-\frac{3\alpha}{2})} \Theta_j.
\end{align*}
The operator $\Delta_k \Lambda^{-\frac{\alpha}{2}} \nabla \cdot$ has a total order of $1-\frac{\alpha}{2}$. Thus, the second term $I_{3k}^{2}$ can be estimated as
\begin{align*}
\|I_{3k}^{2}\|_{L^{p}}&\lesssim
    \sum_{j \ge k-2} 2^{k(1-\frac{\alpha}{2})} \|\Delta_j u_H\|_{L^q} \| R_{\alpha} (\widetilde{\Delta}_j \theta)\|_{L^b} \\ &\lesssim \sum_{j \ge k-2} 2^{k(1-\frac{\alpha}{2})} 2^{-j} \|H\|_{L^q} 2^{(1-\alpha)j} \|\widetilde{\Delta}_j \theta\|_{L^b} \\
    &= \|H\|_{L^q} \sum_{j \ge k-2} 2^{-(j-k)(1-\frac{\alpha}{2})} \Theta_j.
\end{align*}
Notice that $\alpha < \frac{4}{3}$ guarantees that both exponents are strictly positive: $2-\frac{3\alpha}{2} > 0$ and $1-\frac{\alpha}{2}>0$. Hence, both sums represent geometrically decaying series for $j \ge k-2$. Hence, we have
$$\left( \sum_{k} \|I_{3k}\|_{L^p}^2 \right)^{\frac{1}{2}} \lesssim\left( \sum_{k} \|I_{3k}^{1}\|_{L^p}^2 \right)^{\frac{1}{2}}+\left( \sum_{k} \|I_{3k}^{2}\|_{L^p}^2 \right)^{\frac{1}{2}} \lesssim \|H\|_{L^q} \|\theta\|_{B_{b,2}^{\frac{2-3\alpha}{2}}}.$$
Combining the estimates for $I_{1k}$, $I_{2k}$, and $I_{3k}$, we establish that the $\ell^2$ sum of the dyadic blocks is bounded by $\|H\|_{L^q} \|\theta\|_{B_{b,2}^{\frac{2-3\alpha}{2}}}$. Because $p \ge 2$, the Littlewood-Paley square function embedding $B_{p,2}^0 \hookrightarrow L^p$ implies
\begin{equation*}
    \|X\|_{L^p} \lesssim \left( \sum_{k} \|\Delta_k X\|_{L^p}^2 \right)^{\frac{1}{2}} \lesssim \|H\|_{L^q} \|\theta\|_{B_{b,2}^{\frac{2-3\alpha}{2}}}
\end{equation*}
This yields the desired inequality
\begin{equation*}
    \|\Lambda^{-\frac{\alpha}{2}}[ R_{\alpha} ,\Kop H\cdot\nabla]\theta\|_{L^{p}} \lesssim \norm{H}{L^{q}}\norm{\theta}{B_{b,2}^{\frac{2-3\alpha}{2}}},
\end{equation*}
which completes the proof of the lemma.
\end{proof}

\subsection{Multiplier estimates for $\alpha\in[\frac{2}{3},\frac{10}{13}]$}
\label{sec:upper-critical-multipliers}
Throughout this case, we assume
\begin{equation}\label{eq:upper-critical-range}
\frac23\leq\alpha\leq\frac{10}{13},\qquad
\beta=1-\alpha,\qquad 0<\rho<\frac1{39}.
\end{equation}
\begin{lemma}[Critical-interval multiplier estimates]
\label{lem:upper-critical-multipliers}
Assume \eqref{eq:upper-critical-range}. For $2\leq p\leq\infty$
and $M\geq1$, we have
\begin{equation}\label{eq:upper-B-lowpass}
\|S_M\Bop(\theta,\theta)\|_{L^{p}}
+\sum_{i,j=1}^2\|\mathcal R_i\mathcal R_jS_M\Bop(\theta,\theta)\|_{L^{p}}
\leq C_\rho M^\rho\|\theta\|_{L^{\infty}}\|\theta\|_{L^{p}}.
\end{equation}
For every $\gamma\in(0,1)$, it holds
\begin{equation}\label{eq:upper-uB-holder}
\|u_B\|_{C^\gamma}\leq C_\gamma\|\theta\|_{L^{\infty}}^2.
\end{equation}
For $2\leq p<\infty$, the estimate \eqref{eq:Q-top} remains valid, namely,
\begin{equation*}
\norm{\Lambda^{-\frac{\alpha}{2}}\Qop(\Kop f;\theta,\theta)}{L^{p}}
\lesssim \norm{\theta}{L^{\infty}}^2
\sum_{M\geq1}M^{\frac{4-7\alpha}{2}}\norm{S_{c_1M}f}{L^{p}},
\end{equation*}
If $\frac{1}{p}=\frac{1}{q}+\frac{1}{b}$ with $1<q,b<\infty$, then we have
\begin{equation}\label{eq:upper-CH}
\|\Lambda^{-\frac{\alpha}{2}}[ R_{\alpha} ,\Kop H\cdot\nabla]\theta\|_{L^{p}}
\leq C_\rho\|H\|_{L^{q}}\|\theta\|_{B_{b,2}^{\rho}}.
\end{equation}
Moreover,
\begin{align}
\|\Lambda^{-\frac{\alpha}{2}}\Qop(\Top\theta;\theta,\theta)\|_{L^{p}}
+\|\Lambda^{-\frac{\alpha}{2}}[ R_{\alpha} ,u_B\cdot\nabla]\theta\|_{L^{p}}
&\leq C_\rho\|\theta\|_{L^{\infty}}^2
\|\theta\|_{B_{p,1}^{\rho}},\label{eq:upper-cubic}\\
\|\Lambda^{-\frac{\alpha}{2}}\Qop(u_B;\theta,\theta)\|_{L^{p}}
&\leq C_\rho\|\theta\|_{L^{\infty}}^3\|\theta\|_{L^{p}},
\label{eq:upper-quartic}\\
\|\Lambda^{-\frac{\alpha}{2}}\Cres(\theta,\theta)\|_{L^{p}}
&\leq C\|\theta\|_{L^{\infty}}\|\theta\|_{L^{p}}.
\label{eq:upper-resonance}
\end{align}
\end{lemma}

\begin{proof}
\noindent \textbf{Low-pass and H\"older bounds \eqref{eq:upper-B-lowpass} \& \eqref{eq:upper-uB-holder}}.
The localized symbol bounds for $\Bop$ depend on the ratios $\left(\frac{L}{N}\right)^{1-\alpha}$ and $\left(\frac{N}{K}\right)^\alpha$. Since $1-\alpha > 0$ and $\alpha > 0$, the dyadic summation bounds in \eqref{eq:b-bounds} still hold identically block-by-block.
However, when summing over output blocks $N \le M$, we encounter $\sum_{1 \lesssim N \lesssim M} N^{2-3\alpha}$. Since $\alpha \ge \frac{2}{3}$, we have $2-3\alpha \le 0$.
If $\alpha >\frac{2}{3}$, the series converges to a constant. If $\alpha=\frac{2}{3}$, the series grows logarithmically
\begin{equation*}
    \sum_{1\lesssim N\lesssim M}N^{2-3\alpha} \leq C(1+\log(2+M)) \leq C_\rho M^\rho \quad \text{for any } \rho > 0.
\end{equation*}
This gives $\|S_M\Bop(\theta,\theta)\|_p \leq C_\rho M^\rho A_\theta\|\theta\|_p$. Because we apply Riesz transforms to annuli bounded strictly away from zero (due to the high-output cutoff $\chi_{\text{out}}$), their kernels are uniformly integrable, allowing this bound to hold even for $p=\infty$.

Regarding $u_B = \Kop \Bop(\theta,\theta)$, the multiplier $\Kop$ adds a factor of $N^{-1}$. Using the block bound $\|\Delta_N \Bop\|_{L^{\infty}} \lesssim A_\theta^2 N^{2-3\alpha}$, we get
\[ \|\Delta_Nu_B\|_{L^{\infty}} \lesssim A_\theta^2 N^{1-3\alpha}. \]
To evaluate the $C^\gamma$ H\"older norm, we multiply by $N^\gamma$
\[ N^\gamma \|\Delta_Nu_B\|_{L^{\infty}} \lesssim A_\theta^2 N^{\gamma+1-3\alpha}. \]
Since $\alpha\geq\frac{2}{3}$, for any $\gamma\in (0,1)$, the exponent $\gamma+1-3\alpha$ is strictly negative. The sequence decays geometrically for high frequencies, and low frequencies are identically zero. Taking the supremum yields $\|u_B\|_{C^\gamma} \leq C_\gamma A_\theta^2$.

\vspace{0.5 em}
\noindent \textbf{The mixed commutator \eqref{eq:upper-CH}} \\
After checking the proof of Lemma \ref{lem:CH}, the desired estimate \eqref{eq:CH} is still valid for $\alpha\in (0,\frac{4}{3})$. Therefore, when $\alpha\in[\frac{2}{3},\frac{10}{13}]$, we have for any $\rho>0$
$$\norm{\Lambda^{-\frac{\alpha}{2}}[ R_{\alpha} ,\Kop H\cdot\nabla]\theta}{L^{p}}
\lesssim\norm{H}{L^{q}}\norm{\theta}{B_{b,2}^{\frac{2-3\alpha}{2}}}
\lesssim\norm{H}{L^{q}}\|\theta\|_{B_{b,2}^{\rho}}.$$
This gives \eqref{eq:upper-CH}.

\vspace{0.5em}
\noindent \textbf{The cubic bounds \eqref{eq:upper-cubic}}.
For $f =  R_{\alpha} \theta$, we interchange the summations in the estimate as before. Instead of yielding the index $\frac{3(2-3\alpha)}{2}$, the substitution gives $U^{\frac{3(2-3\alpha)}{2}}$. Because $\frac{3(2-3\alpha)}{2}\le 0$, we cannot use it directly to define a positive Besov space. However, for $U \gtrsim 1$, we can always bound $U^{\frac{3(2-3\alpha)}{2}} \le U^\rho$ for any $\rho > 0$. As a result, we get
\[ \sum_{U \gtrsim 1} U^{\frac{3(2-3\alpha)}{2}} \|\widetilde\Delta_U\theta\|_{L^{p}} \leq \sum_{U \gtrsim 1} U^\rho \|\widetilde\Delta_U\theta\|_{L^{p}} \lesssim \|\theta\|_{B_{p,1}^\rho}. \]

For the cubic commutator involving $u_B$, use the decomposition
$I_{1k}+I_{2k}+I_{3k}$ from Lemma~\ref{lem:other-cubic}.
The fixed output cutoff of $u_B$ gives a fixed index $j_0$
below which its dyadic blocks vanish. Thus, for $j\geq0$,
\[
\|\nabla S_{j-1}u_B\|_{L^\infty}
\lesssim \|\theta\|_{L^\infty}^2
\sum_{j_0\leq m\leq j-2}2^{m(2-3\alpha)}
\lesssim \|\theta\|_{L^\infty}^2(1+j),
\]
where an empty sum is zero. In particular, the factor $1+j$
retains the logarithmic loss at $\alpha=\frac{2}{3}$.
For $k\geq0$, the localized kernel commutator estimate yields
\begin{align*}
	\|I_{1k}\|_{L^p}
	&\lesssim \|\theta\|_{L^\infty}^2
	\sum_{\substack{j\geq0\\|j-k|\leq C}}
	(1+j)2^{(1-3\frac{\alpha}{2})j}\|\Delta_j\theta\|_{L^p}\\
	&\lesssim_\rho \|\theta\|_{L^\infty}^2
	\sum_{\substack{j\geq0\\|j-k|\leq C}}
	2^{\rho j}\|\Delta_j\theta\|_{L^p},
\end{align*}
because $1-\frac{3\alpha}{2}\leq0$ and
$(1+j)\lesssim_\rho2^{\rho j}$ for every $\rho>0$.
Consequently, we have
\[
\sum_{k\geq0}\|I_{1k}\|_{L^p}
\lesssim_\rho \|\theta\|_{L^\infty}^2\|\theta\|_{B_{p,1}^{\rho}}.
\]

For the remaining paraproduct terms, put
$s_\alpha=3-\frac{9\alpha}{2}\leq0$.
Their dyadic estimates in Lemma~\ref{lem:other-cubic} remain
valid on the present interval: the convolution decay exponents are
\[
\frac{9\alpha}{2}-2>0,\qquad
\frac{7\alpha}{2}-1>0,\qquad
2-\frac{3\alpha}{2}>0,\qquad
1-\frac{\alpha}{2}>0.
\]
Hence, by discrete Young's inequality,
\[
\sum_{k\geq0}
\bigl(\|I_{2k}\|_{L^p}+\|I_{3k}\|_{L^p}\bigr)
\lesssim \|\theta\|_{L^\infty}^2\|\theta\|_{B_{p,1}^{s_\alpha}}
\lesssim_\rho \|\theta\|_{L^\infty}^2\|\theta\|_{B_{p,1}^{\rho}}.
\]

For the lowest output block, $\Delta_{-1} R_{\alpha} =0$.
Using $\nabla\cdot u_B=0$, the integrable kernel of
$\Delta_{-1}\Lambda^{-\frac{\alpha}{2}}\operatorname{div}$,
and the comparable input frequencies, we obtain
\begin{align*}
	 \|\Delta_{-1}\Lambda^{-\frac{\alpha}{2}}
	[ R_{\alpha} ,u_B\cdot\nabla]\theta\|_{L^p}
	& \lesssim \|\theta\|_{L^\infty}^2
	\sum_{j\geq\max\{-1,j_0-C\}}
	2^{(2-4\alpha)j}\|\widetilde\Delta_j\theta\|_{L^p}
	\\&\lesssim \|\theta\|_{L^\infty}^2\|\theta\|_{L^p}.
\end{align*}
Combining these bounds gives
\[
\|\Lambda^{-\frac{\alpha}{2}}[ R_{\alpha} ,u_B\cdot\nabla]\theta\|_{L^p}
\lesssim_\rho\|\theta\|_{L^\infty}^2
\|\theta\|_{B_{p,1}^{\rho}}.
\]
Together with the preceding estimate for
$\Qop(\Top\theta;\theta,\theta)$,
this proves \eqref{eq:upper-cubic}.

\vspace{0.5 em}
\noindent \textbf{The trilinear and quartic bounds}.
A careful check of the localized symbol bounds in Lemma \ref{lem:Q} reveals that the potentially restrictive fractional powers are $2-\frac{5\alpha}{2}$, $1-\alpha$, $3-3\alpha$, and $2\alpha-1$. Over the interval $\frac{2}{3} \le \alpha \le \frac{10}{13}$, all these exponents remain strictly positive
\begin{equation}
2-\frac{5\alpha}{2}\geq\frac1{13},\quad
1-\alpha\geq\frac3{13},\quad
3-3\alpha\geq\frac9{13},\quad
2\alpha-1\geq\frac13.\nonumber
\end{equation}
Because these convergence factors remain strictly positive, the absolute sums over intermediate frequency scales converge just as they did in Lemma \ref{lem:Q}. Consequently, the estimate \eqref{eq:Q-top} holds verbatim.

For the quartic term ($f = \Bop(\theta,\theta)$), substituting the low-pass bound \eqref{eq:upper-B-lowpass} into \eqref{eq:Q-top} implies
\[ \|\Lambda^{-\frac{\alpha}{2}}\Qop(u_B;\theta,\theta)\|_{L^{p}} \lesssim A_\theta^2 \sum_{M \ge 1} M^{\frac{4-7\alpha}{2}} \left( C_\rho M^\rho A_\theta \|\theta\|_{L^{p}} \right) = C_\rho A_\theta^3 \|\theta\|_{L^{p}} \sum_{M \ge 1} M^{\frac{4-7\alpha}{2}+\rho}. \]
If we choose $\rho <\frac{1}{3}$ (and our choice $\rho <\frac{1}{39}$ satisfies this), the exponent $\frac{4-7\alpha}{2}+\rho<0$. The geometric series converges, proving \eqref{eq:upper-quartic}.

\vspace{0.5em}
\noindent \textbf{The resonant forcing \eqref{eq:upper-resonance}}.
Because the exponent is strictly negative, the geometric sum over high-output blocks converges without the need to borrow any positive fractional regularity from $\theta$. Therefore, the baseline $L^p \times L^\infty$ bound from Lemma \ref{lem:res} applies exactly:
\[ \|\Lambda^{-\frac{\alpha}{2}}\Cres(\theta,\theta)\|_{L^{p}} \lesssim \|\theta\|_{L^{\infty}} \|\theta\|_{L^{p}}. \]
This completes the proof of the lemma.
\end{proof}

\section{$L^m$ estimates for the corrected vorticity}
\label{sec:energy}

We estimate the corrected-vorticity equation by signed-power testing.
The strict temperature interpolation gaps allow all coefficients
to be controlled by the $L^2$ temperature dissipation energy.

\subsection{Coercivity and the fractional chain rule}

For $m\geq2$, we denote
\begin{equation*}
E_m=\norm{H}{L^m}^m,\qquad
\Phi_m=\sgn(H)|H|^{\frac{m}{2}},\qquad
D_m=\norm{\Lambda^{\frac{\alpha}{2}}\Phi_m}{L^2}.
\end{equation*}

\begin{lemma}
\label{lem:signed}
For real-valued $H$, it holds
\begin{align}
\langle\Lambda^\alpha H,|H|^{m-2}H\rangle\triangleq\int_{\mathbb{R}^2}\Lambda^\alpha H |H|^{m-2}H\,dx
&\geq \frac{4(m-1)}{m^2}D_m^2,\label{eq:coercive}\\
\norm{H}{L^{r_m}}^{\frac{m}{2}}&\lesssim D_m,\qquad
r_m=\frac{2m}{2-\alpha}.\label{eq:embedding}
\end{align}
If $m\leq q\leq r_m$ and
$\frac{1}{\ell}=\frac{1}{2}+\frac{\frac{m}{2}-1}{q}$, then
\begin{equation}\label{eq:chain}
\norm{\Lambda^{\frac{\alpha}{2}}(|H|^{m-2}H)}{L^{\ell}}
\lesssim D_m\norm{H}{L^q}^{\frac{m}{2}-1}.
\end{equation}
\end{lemma}

\begin{proof}
Set $\phi(z)=\sgn(z)|z|^{\frac{m}{2}}$ and $\psi(z)=|z|^{m-2}z$.
For $a\geq b$, Cauchy--Schwarz gives
\begin{align*}
|\phi(a)-\phi(b)|^2
&=\frac{m^2}{4}\left|\int_b^a|z|^{\frac{m}{2}-1}\,dz\right|^2\\
&\leq\frac{m^2}{4}(a-b)\int_b^a|z|^{m-2}\,dz\\
&=\frac{m^2}{4(m-1)}(a-b)[\psi(a)-\psi(b)].
\end{align*}
Interchanging $a$ and $b$ gives the same inequality for all real
$a,b$, including opposite signs. The symmetric kernel formula implies
\begin{align*}
\langle\Lambda^\alpha H,\psi(H)\rangle
&=\frac{c_\alpha}{2}\iint_{\RR^2\times\RR^2}
\frac{[H(x)-H(y)][\psi(H(x))-\psi(H(y))]}
{|x-y|^{2+\alpha}}\,dx\,dy\\
&\geq\frac{4(m-1)}{m^2}\frac{c_\alpha}{2}
\iint\frac{|\phi(H(x))-\phi(H(y))|^2}
{|x-y|^{2+\alpha}}\,dx\,dy\\
&=\frac{4(m-1)}{m^2}
\|\Lambda^{\frac{\alpha}{2}}\Phi_m\|_2^2,
\end{align*}
which yields \eqref{eq:coercive}. The two-dimensional Sobolev inequality gives
\[
\|H\|_{L^{r_m}}^{\frac{m}{2}}
=\|\Phi_m\|_{\frac{4}{2-\alpha}}
\lesssim\|\Lambda^{\frac{\alpha}{2}}\Phi_m\|_{L^2}=D_m,
\]
which proves \eqref{eq:embedding}.
For the chain rule, we set
\[
F(z)=\sgn(z)|z|^{2-\frac{2}{m}},\qquad
F(\Phi_m)=\psi(H),\qquad
|F'(\Phi_m)|=(2-\frac{2}{m})|H|^{\frac{m}{2}-1}.
\]
Let $r=\frac{q}{\frac{m}{2}-1}$. The exponents satisfy
\[
\frac1\ell=\frac12+\frac1r,\qquad
r>1,\qquad1<\ell<2,\qquad0<\frac{\alpha}{2}<1.
\]
The fractional chain rule in \cite[Proposition 1.1]{Fujiwara}
and \cite{CW} therefore yields
\begin{align*}
\|\Lambda^{\frac{\alpha}{2}}(|H|^{m-2}H)\|_{L^{\ell}}
&=\|\Lambda^{\frac{\alpha}{2}}F(\Phi_m)\|_{L^{\ell}}\\
&\lesssim\|F'(\Phi_m)\|_{L^{r}}\|\Lambda^{\frac{\alpha}{2}}\Phi_m\|_{L^{2}}\\
&\lesssim\|H\|_{L^{q}}^{\frac{m}{2}-1}D_m,
\end{align*}
which is \eqref{eq:chain}.
Smooth approximation of $F$ at zero justifies the application. Hence, we complete the proof of Lemma \ref{lem:signed}.
\end{proof}

\subsection{The corrected-vorticity energy inequality}

\begin{proposition}[Corrected-vorticity estimate]\label{prop:energy}
Assume
\begin{equation}
\frac8{13}<\alpha<\frac23,\qquad \frac{4-4\alpha-3\alpha^{2}}{\alpha}<\beta<\alpha,\nonumber
\end{equation}
and suppose $m>2$ can be chosen with
\begin{align}\label{eq:bF}
\frac{1}{m}>\max\left\{\frac{3(2+\alpha)(2-3\alpha)-\alpha\beta}{4\beta},\ \frac{(4+2\alpha)[(\frac{2}{\alpha}-\frac{3}{2})\beta-\alpha]-\alpha\beta}{4\beta}
\right\}\triangleq A(\alpha,\beta).
\end{align}
Then every smooth solution on $[0,T]$ satisfies
\begin{equation}\label{eq:energy-conclusion}
\sup_{0\leq t\leq T}\norm{H(t)}{L^m}^m
+\int_0^T D_m(t)^2\,dt
\leq C(T,u_{0},\theta_{0}).
\end{equation}
\end{proposition}

\begin{proof}
We multiply \eqref{eq:H-equation} by $|H|^{m-2}H$ and integrate over space. Because the velocity field $u$ is divergence-free, the transport term vanishes $\int (u \cdot \nabla H) |H|^{m-2}H \, dx = 0$.
Applying the signed coercivity bound from Lemma~\ref{lem:signed} to the dissipative term $\Lambda^\alpha H$, we obtain the basic energy inequality
\begin{equation}\label{eq:test-energy}
\frac{1}{m} \frac{dE_m}{dt}
+\frac{4(m-1)}{m^2}D_m^2
\leq \sum_F|\langle F,|H|^{m-2}H\rangle|,
\end{equation}
where the sum is over the seven forcing terms $F$ on the right-hand side of \eqref{eq:H-equation}.

For any forcing term $F$, we apply H$\rm\ddot{o}$lder's inequality and the fractional chain rule (Lemma~\ref{lem:signed}). If we choose parameters $q$ and $\ell$ such that $\frac{1}{\ell}=\frac{1}{2}+\frac{\frac{m}{2}-1}{q}$, and let $\ell'$ be the H$\rm\ddot{o}$lder conjugate of $\ell$, then we have
\begin{align}
|\langle F,|H|^{m-2}H\rangle| &= \left| \langle \Lambda^{-\frac{\alpha}{2}}F, \Lambda^{\frac{\alpha}{2}}(|H|^{m-2}H)\rangle \right| \nonumber \\
&\leq \norm{\Lambda^{-\frac{\alpha}{2}}F}{L^{\ell'}} \norm{\Lambda^{\frac{\alpha}{2}}(|H|^{m-2}H)}{L^{\ell}} \nonumber \\
&\lesssim \norm{\Lambda^{-\frac{\alpha}{2}}F}{L^{\ell'}} D_m \norm{H}{L^{q}}^{\frac{m}{2}-1}. \label{eq:transfer}
\end{align}
We now bound each of the forcing terms using this transfer estimate.

\vspace{0.5 em}
\noindent \textbf{The mixed commutator $F = [ R_{\alpha} , u_H \cdot \nabla]\theta$} \\
We aim to place the temperature in $L^{b_H}$. We let
$$\ell' =\frac{mb_H}{b_H+m-2}>2.$$
By the conjugate relation, this forces $\frac{1}{q_H}=\frac{1}{\ell'}-\frac{1}{b_H}$, which yields
$$q_H = \frac{mb_H}{b_H-2}.$$
By interpolation between $L^m$ and $L^{\frac{2m}{2-\alpha}}$, we have
\[ \norm{H}{L^{q_H}} \lesssim \norm{H}{L^{m}}^{1-\lambda_H} \norm{H}{L^{\frac{2m}{2-\alpha}}}^{\lambda_H} \lesssim (E_m^{\frac{1}{m}})^{1-\lambda_H} (D_m^{\frac{2}{m}})^{\lambda_H}, \]
where the interpolation parameter is uniquely determined as $\lambda_H=\frac{4}{\alpha b_H}$.
For the $F$ part, Lemma \ref{lem:CH} gives
$$\|\Lambda^{-\frac{\alpha}{2}}F\|_{L^{\ell'}} \lesssim \|H\|_{L^{q_H}} \|\theta\|_{\Bes{b_H}{2}{\frac{2-3\alpha}{2}}}.$$
Lemma \ref{lem:theta-GN} allows us to derive
$$\|\theta\|_{\Bes{b_H}{2}{\frac{2-3\alpha}{2}}} \lesssim 1 + Z^{\frac{2}{b_H}},$$
where we have used the restriction
\begin{align}\label{asdfvb11}
\frac{2-3\alpha}{2}<\frac{\beta}{b_H}\Leftrightarrow b_H<\frac{2\beta}{2-3\alpha}.
\end{align}
Combining these, we get
\begin{align}
|\langle F, |H|^{m-2}H \rangle| &\lesssim \|H\|_{L^{q_H}} \|\theta\|_{\Bes{b_H}{2}{\frac{2-3\alpha}{2}}}D_m \norm{H}{L^{q_H}}^{\frac{m}{2}-1}\nonumber\\&\lesssim (1 + Z^{\frac{2}{b_H}})D_m \norm{H}{L^{q_H}}^{\frac{m}{2}}
\nonumber\\&\lesssim (1 + Z^{\frac{2}{b_H}})D_m \left((E_m^{\frac{1}{m}})^{1-\lambda_H} (D_m^{\frac{2}{m}})^{\lambda_H}\right)^{\frac{m}{2}}
\nonumber\\&\lesssim (1 + Z^{\frac{2}{b_H}}) E_m^{\frac{1-\lambda_H}{2}} D_m^{1+\lambda_H}\nonumber\\&\leq
\varepsilon D_m^2 + C_\varepsilon (1+Z^{\frac{2}{b_H}})^{\frac{2}{1-\lambda_H}} E_m\nonumber\\&\leq
\varepsilon D_m^2 + C_\varepsilon (1+Z^{2}) E_m,\nonumber
\end{align}
where we have used
$$\frac{2}{b_H}\frac{2}{1-\lambda_H} = \frac{4}{b_H-\frac{4}{\alpha}}\leq2\Leftrightarrow b_H\geq 2+\frac{4}{\alpha}.$$
This restriction along with \eqref{asdfvb11} yields
\begin{align}
2+\frac{4}{\alpha}\leq b_H<\frac{2\beta}{2-3\alpha},\nonumber
\end{align}
which requires
$$\beta>\frac{4-4\alpha-3\alpha^{2}}{\alpha}.$$

\vspace{0.5 em}
\noindent \textbf{The cubic terms and resonance $F_\theta$}. For simplicity, we denote
$$F_\theta\triangleq\Qop(\Top\theta;\theta,\theta) + [ R_{\alpha} ,u_B\cdot\nabla]\theta + \Cres(\theta,\theta).$$
We place the temperature in $L^{b_F}$ and set $\ell' = b_F$.
The conjugate relation forces $q_F = \frac{(m-2)b_F}{b_F-2}$. Interpolating $L^{q_F}$ between $L^m$ and $L^{\frac{2m}{2-\alpha}}$ gives
$$\lambda_F = \frac{4(m-b_F)}{\alpha b_F(m-2)}.$$
The following restriction ensures $\lambda_F\in [0,1]$,
\begin{align}\label{qwxmj011}
\frac{4m}{\alpha(m-2)+4}\leq b_F\leq m.
\end{align}
The $H$-dependent part of \eqref{eq:transfer} becomes
\[ D_m \norm{H}{L^{q_F}}^{\frac{m}{2}-1} \lesssim D_m \left( E_m^{\frac{1-\lambda_F}{m}} D_m^{\frac{2\lambda_F}{m}} \right)^{\frac{m}{2}-1} = E_m^{\frac{(m-2)(1-\lambda_F)}{2m}} D_m^{1+\frac{(m-2)\lambda_F}{m}}. \]
It follows from \eqref{eq:Q-cubic} that
\begin{align}\label{ADXZ1}
\norm{\Lambda^{-\frac{\alpha}{2}}\Qop(\Top\theta;\theta,\theta)}{L^{b_F}}
&\lesssim \norm{\theta}{L^{\infty}}^2\norm{\theta}{B_{b_F,1}^{\frac{3(2-3\alpha)}{2}}}.
\end{align}
According to \eqref{eq:other-cubic}, one has
\begin{equation}\label{ADXZ2}
\norm{\Lambda^{-\frac{\alpha}{2}}[ R_{\alpha} ,u_B\cdot\nabla]\theta}{L^{b_F}}
\lesssim\norm{\theta}{L^{\infty}}^2\norm{\theta}{B_{b_F,1}^{\frac{3(2-3\alpha)}{2}}}.
\end{equation}
By \eqref{eq:res-Besov}, we get
\begin{equation}\label{ADXZ3}
\norm{\Lambda^{-\frac{\alpha}{2}}\Cres(\theta,\theta)}{L^{b_F}}
\lesssim\norm{\theta}{L^\infty}\norm{\theta}{B_{b_F,1}^{\frac{(4-3\alpha)
\beta-2\alpha^{2}}{2\alpha}}}.
\end{equation}
For simplicity, we denote
$$\Gamma(\alpha,\beta)\triangleq \max\left\{\frac{3(2-3\alpha)}{2},\, \frac{(4-3\alpha)
\beta-2\alpha^{2}}{2\alpha}\right\}>0.$$
Combining \eqref{ADXZ1}, \eqref{ADXZ2} and \eqref{ADXZ3}, and using Lemma \ref{lem:theta-GN}, we obtain
\begin{align}
\|\Lambda^{-\frac{\alpha}{2}}F_\theta\|_{L^{b_F}} &\lesssim \norm{\theta}{B_{b_F,1}^{\frac{3(2-3\alpha)}{2}}}+ \norm{\theta}{B_{b_F,1}^{\frac{(4-3\alpha)
\beta-2\alpha^{2}}{2\alpha}}}
\nonumber\\&\lesssim \norm{\theta}{B_{b_F,1}^{\Gamma(\alpha,\beta)}}
\nonumber\\&\lesssim 1 + Z^{\frac{2}{b_F}},\nonumber
\end{align}
where in the last line we have used
\begin{align}\label{qwxmj012}
b_F<\frac{\beta}{\Gamma(\alpha,\beta)}.
\end{align}
Combining above estimates yields
\begin{align}\label{askhg97}
|\langle F_\theta, |H|^{m-2}H \rangle| &\lesssim (1+Z^{\frac{2}{b_F}}) E_m^{\frac{(m-2)(1-\lambda_F)}{2m}} D_m^{1+\frac{(m-2)\lambda_F}{m}}\nonumber\\
&\leq \varepsilon D_m^2 + C_\varepsilon (1+Z^{\gamma_F}) E_m^{\frac{(m-2)(1-\lambda_F)}{m - (m-2)\lambda_F}}
\nonumber\\
&\leq \varepsilon D_m^2 + C_\varepsilon (1+Z^{\gamma_F})(1+E_m),
\end{align}
where $\gamma_F$ is given by
\[ \gamma_F = \frac{2}{b_F}\, \frac{2}{1 - \frac{(m-2)\lambda_F}{m}}= \frac{4\alpha m}{b_F(\alpha m + 4) - 4m}. \]
To apply Gronwall's inequality, we absolutely must have $\gamma_F\leq2$, precisely,
\begin{align}\label{qwxmj013}
\frac{4\alpha m}{b_F(\alpha m + 4) - 4m}\leq2 \iff   b_F\geq \frac{m(4+2\alpha)}{\alpha m + 4}.
\end{align}
It thus follows from \eqref{askhg97} that
\begin{align*}
|\langle F_\theta, |H|^{m-2}H \rangle|
 \leq \varepsilon D_m^2 + C_\varepsilon (1+Z^{2})(1+E_m).
\end{align*}
Collecting \eqref{qwxmj011}, \eqref{qwxmj012} and \eqref{qwxmj013} implies
$$\max\left\{\frac{4m}{\alpha(m-2)+4},\ \frac{m(4+2\alpha)}{\alpha m + 4}\right\}=\frac{m(4+2\alpha)}{\alpha m + 4}\leq b_F\leq m,\quad b_F<\frac{\beta}{\Gamma(\alpha,\beta)}.$$
Therefore, to ensure the validity of $b_F$, it requires
$$\frac{m(4+2\alpha)}{\alpha m + 4}< \min\left\{m,\ \frac{\beta}{\Gamma(\alpha,\beta)}\right\}.$$
The above inequality is true as long as
\begin{align*}
\frac{1}{m}>\max\left\{\frac{3(2+\alpha)(2-3\alpha)-\alpha\beta}{4\beta},\ \frac{(4+2\alpha)[(\frac{2}{\alpha}-\frac{3}{2})\beta-\alpha]-\alpha\beta}{4\beta}
\right\},
\end{align*}
which is \eqref{eq:bF}.

\vspace{0.5 em}
\noindent \textbf{The normal-form term containing $u_H$}.
Making use of \eqref{eq:Q-mixed} gives
\begin{align} \label{erpgh834}
|\langle \Qop(u_H;\theta,\theta), |H|^{m-2}H \rangle|
&\lesssim
\norm{\Lambda^{-\frac{\alpha}{2}}\Qop(u_H;\theta,\theta)}{L^{m}}
\norm{\Lambda^{\frac{\alpha}{2}}(|H|^{m-2}H)}{L^{\frac{m}{m-1}}}
\nonumber\\ &\lesssim \|\theta\|_{L^{\infty}}^2 \|H\|_{L^{m}} D_m \|H\|_{L^{m}}^{\frac{m}{2}-1} \nonumber\\ &\lesssim D_m E_m^{\frac{1}{2}}
\nonumber\\ &\leq \varepsilon D_m^2 + C_\varepsilon E_m.
\end{align}

\vspace{0.5 em}
\noindent \textbf{Quartic and low-Output forcing}.
The remaining terms depend only on the temperature, not on the velocity. Their transferred $L^m$ norms are bounded by a constant $C_0$ depending only on initial data (via the maximum principle for $\theta$). More precisely, we deduce from \eqref{eq:low-bound} and \eqref{eq:Q-quartic} that
\begin{align*}
&|\langle F_{\mathrm{low}}, |H|^{m-2}H \rangle|+|\langle \Qop(u_B;\theta,\theta), |H|^{m-2}H \rangle|\\
&\lesssim
(\norm{\Lambda^{-\frac{\alpha}{2}}F_{\mathrm{low}}}{L^{m}}
+\norm{\Lambda^{-\frac{\alpha}{2}}\Qop(u_B;\theta,\theta)}{L^{m}})
\norm{\Lambda^{\frac{\alpha}{2}}(|H|^{m-2}H)}{L^{\frac{m}{m-1}}}
\\ &\lesssim (1+\|\theta\|_{L^{\infty}}^3)\norm{\theta}{L^m}  D_m \|H\|_{L^{m}}^{\frac{m}{2}-1} \\ &\lesssim C_0 D_m E_m^{\frac{m-2}{2m}} \\ &\leq \varepsilon D_m^2 + C_\varepsilon E_m^{\frac{m-2}{m}} \\ &\leq \varepsilon D_m^2 + C_\varepsilon (1 + E_m).
\end{align*}

Finally, returning to the energy equality \eqref{eq:test-energy}, we select $\varepsilon$ sufficiently small to arrive at the differential inequality
\begin{equation}
\frac{d}{dt}E_m+c_m D_m^2
\leq C_0(1+Z^2)(1+E_m).\nonumber
\end{equation}
Applying Gronwall's inequality and using \eqref{eq:theta-energy}, it immediately yields \eqref{eq:energy-conclusion}. Therefore, we complete the proof of Proposition \ref{prop:energy}.
\end{proof}

\subsection{Energy closure on the remaining critical interval $\frac{2}{3}\leq\alpha\leq\frac{10}{13}$}
\label{sec:upper-critical-energy}

In the case when $\frac{2}{3}\leq\alpha\leq\frac{10}{13}$, some inequalities in Proposition~\ref{prop:energy}
are replaced by positive-$\rho$ interpolation gaps. The following
choice works on the entire additional interval.

\begin{proposition}[Critical-interval energy estimate]
\label{prop:upper-critical-energy}
Let $\beta=1-\alpha$ with $\frac{2}{3}\leq\alpha\leq\frac{10}{13}$.
We take
\begin{equation}
m=\frac{13}{4},\qquad b_H=9,\qquad b_F=3,\qquad
\rho=\frac1{100}.\nonumber
\end{equation}
Then every smooth solution on $[0,T]$ satisfies \eqref{eq:energy-conclusion} with $m=\frac{13}{4}$, namely,
\begin{equation}\label{esdsdw23}
\sup_{0\leq t\leq T}\norm{H(t)}{L^{\frac{13}{4}}}^{\frac{13}{4}}
+\int_0^T D_{\frac{13}{4}}(t)^2\,dt
\leq C(T,u_{0},\theta_{0}).
\end{equation}
\end{proposition}
\begin{remark}
Under the assumption $\beta=1-\alpha$ with $\frac{2}{3}\leq\alpha\leq\frac{10}{13}$, we can more easily obtain the following estimate (the detailed proof is omitted for brevity)
$$\sup_{0\leq t\leq T}\norm{H(t)}{L^2}^2
+\int_0^T D_2(t)^2\,dt
\leq C(T,u_{0},\theta_{0}).$$
\end{remark}
\begin{proof}[The proof of Proposition \ref{prop:upper-critical-energy}]
Since $\beta\geq\frac{3}{13}$, the interpolation gap is
$\rho<\frac{\beta}{9}$, which along with Lemma~\ref{lem:theta-GN} therefore gives
\[
\|\theta\|_{B_{9,2}^{\rho}}\leq C_0(1+Z^{\frac{2}{9}}),
\qquad
\|\theta\|_{B_{3,1}^{\rho}}\leq C_0(1+Z^{\frac{2}{3}}).
\]
Use the signed-power test \eqref{eq:test-energy} and the transfer
estimate \eqref{eq:transfer}. For the mixed commutator, choose
\[
q_H=\frac{mb_H}{b_H-2}=\frac{117}{28},\qquad
\lambda_H=\frac4{9\alpha},\qquad
\ell'_H=\frac{mb_H}{b_H+m-2}=\frac{117}{41}>2.
\]
Here $\frac{1}{\ell'_H}=\frac{1}{q_H}+\frac{1}{9}$, and
\[
m<q_H<\frac{39}{8}\leq r_m\equiv\frac{2m}{2-\alpha},
\qquad
\|H\|_{q_H}\lesssim E_m^{\frac{1-\lambda_H}{m}}
D_m^{\frac{2\lambda_H}{m}}.
\]
By \eqref{eq:upper-CH}, we have
\begin{align}
|\langle[ R_{\alpha} ,u_H\cdot\nabla]\theta,|H|^{m-2}H\rangle|
&\lesssim(1+Z^{\frac{2}{9}})E_m^{\frac{1-\lambda_H}{2}}
D_m^{1+\lambda_H}\nonumber\\
&\leq\varepsilon D_m^2+
C_\varepsilon(1+Z^{\gamma_H})E_m,
\nonumber\\
&\leq\varepsilon D_m^2+
C_\varepsilon(1+Z^{2})E_m,\nonumber
\end{align}
where
$$\gamma_H=\frac4{9-\frac{4}{\alpha}} \leq\frac43<2.$$
Let $F_\theta$ be the sum of the two cubic temperature terms and
$\Cres(\theta,\theta)$ in \eqref{eq:H-equation}.
By means of Lemma~\ref{lem:upper-critical-multipliers}, we get
\[
\|\Lambda^{-\frac{\alpha}{2}}F_\theta\|_3
\leq C_0\|\theta\|_{B_{3,1}^{\rho}}
\leq C_0(1+Z^{\frac{2}{3}}).
\]
For the chain rule, we take
\[
q_F=\frac{(m-2)b_F}{b_F-2}=\frac{15}{4},\qquad
\lambda_F=\frac4{15\alpha},\qquad \ell'_F=3.
\]
Again $m<q_F<\frac{39}{8}\leq r_m$ and $0<\lambda_F<1$.
The same interpolation and Young's inequality as in
\eqref{askhg97} yield
\begin{align}
|\langle F_\theta,|H|^{m-2}H\rangle|
&\leq\varepsilon D_m^2+
C_\varepsilon(1+Z^{\gamma_F})E_m^{\varrho_F}
\nonumber\\ &\leq
\varepsilon D_m^2 + C_\varepsilon (1+Z^{2})(1+E_m),\nonumber
\end{align}
where
$$\gamma_F=\frac{52\alpha}{39\alpha-4}
 \leq\frac{52}{33}<2,\qquad
\varrho_F=\frac{15\alpha-4}{39\alpha-4}\in(0,1).$$
In particular, the resonance needs no positive temperature
regularity on this interval.

The estimate \eqref{eq:Q-top} with $f=H$ and $\alpha>\frac{4}{7}$
gives \eqref{eq:Q-mixed}. Thus the term containing $u_H$
contributes at most $\varepsilon D_m^2+C_\varepsilon E_m$,
as in \eqref{erpgh834}.
The quartic term is bounded by \eqref{eq:upper-quartic};
the low-output terms by \eqref{eq:low-bound} and the fixed-cutoff
bound for $(1-\chi_h(D))\partial_1\theta$.
Their contributions are at most
$\varepsilon D_m^2+C_\varepsilon(1+E_m)$.
Absorbing the dissipative terms, we are still able to derive
\begin{equation*}
\frac{d}{dt}E_m+c_mD_m^2
\leq C_0(1+Z^2)(1+E_m).
\end{equation*}
Gronwall's inequality proves \eqref{esdsdw23}. We thus complete the proof of Proposition \ref{prop:upper-critical-energy}.
\end{proof}

\subsection{Higher regularity estimate of $H$}
\label{secasdfgv56y}

To achieve this goal, our first objective here is to evaluate the precise, unadulterated scaling of the mixed commutator $[ R_{\alpha} , \phi\cdot \nabla]\theta$ by controlling it directly with the fractional dissipation of the velocity field, $\Lambda^\alpha \phi$, and the $L^\infty$ amplitude of the temperature field.

\begin{lemma}
Assume $\frac{1}{2}< \alpha < 1$ and $1 < r < \infty$. For any divergence-free vector field $\phi$ and scalar field $\theta$, we have
\begin{equation}\label{eq:mixed_direct}
\norm{\Delta_k [ R_{\alpha} ,\phi\cdot\nabla]\theta}{L^{r}} \leq C 2^{(2-2\alpha)k} \norm{\Lambda^{\alpha} \phi}{L^{r}} \norm{\theta}{L^{\infty}}.
\end{equation}
\end{lemma}

\begin{proof}
Using Bony's paraproduct decomposition, we decompose the product $\phi \cdot \nabla \theta$ into three parts
\begin{equation*}
    \phi \cdot \nabla \theta = \sum_{j} S_{j-1}\phi \cdot \nabla \Delta_j \theta + \sum_{j} \Delta_j \phi \cdot \nabla S_{j-1}\theta + \sum_{j} \Delta_j \phi \cdot \nabla \widetilde{\Delta}_j \theta.
\end{equation*}
Applying the Littlewood-Paley projection $\Delta_k$ and utilizing frequency localization, we split the commutator into three terms: $I_1 + I_2 + I_3$. Let us first deal with the term $I_1$, namely,
\begin{equation*}
    I_1 = \Delta_k \sum_{|j-k| \le 4} [ R_{\alpha} , S_{j-1}\phi \cdot \nabla] \Delta_j \theta
\end{equation*}
Using the first-order Taylor expansion, we obtain
\begin{align*}
    \|[ R_{\alpha} , S_{j-1}\phi \cdot \nabla] \Delta_j \theta\|_{L^r} &\lesssim \||x|K_j\|_{L^1} \|\nabla S_{j-1}\phi\|_{L^r} \|\nabla \Delta_j \theta\|_{L^\infty}\\& \lesssim 2^{-\alpha j} \|\nabla S_{j-1}\phi\|_{L^r} \|\nabla \Delta_j \theta\|_{L^\infty}.
\end{align*}
By Bernstein's inequality,
$$\|\nabla \Delta_j \theta\|_{L^\infty} \lesssim 2^j \|\theta\|_{L^\infty}.$$
 Since $\alpha < 1$, we use Bernstein's inequality on $\phi$ to get
\begin{align*}
    \|\nabla S_{j-1}\phi\|_{L^r} &\le \sum_{m \le j-2} 2^m \|\Delta_m \phi\|_{L^r} \\& \approx \sum_{m \le j-2} 2^{(1-\alpha)m} 2^{\alpha m} \|\Delta_m \phi\|_{L^r} \\&\lesssim \sum_{m \le j-2} 2^{(1-\alpha)m} \|\Lambda^\alpha \phi\|_{L^r}\\& \lesssim 2^{(1-\alpha)j} \|\Lambda^\alpha \phi\|_{L^r}.
\end{align*}
Substituting these back and noting $j \sim k$, we get
\begin{equation*}
    \|I_1\|_{L^r}   \lesssim 2^{(2-2\alpha)k} \|\Lambda^\alpha \phi\|_{L^r} \|\theta\|_{L^\infty}.
\end{equation*}
Regarding the term $I_2$, we will not use the commutator structure
\begin{equation*}
    I_2 = \Delta_k \sum_{|j-k| \le 4} \Big(  R_{\alpha} (\Delta_j \phi \cdot \nabla S_{j-1}\theta) - \Delta_j \phi \cdot \nabla  R_{\alpha} (S_{j-1}\theta) \Big).
\end{equation*}
Since $ R_{\alpha} $ behaves as an operator of order $1-\alpha$, we thus have
\begin{equation*}
    \| R_{\alpha} (\Delta_j \phi \cdot \nabla S_{j-1}\theta)\|_{L^r} \lesssim 2^{(1-\alpha)k} \|\Delta_j \phi\|_{L^r} \|\nabla S_{j-1}\theta\|_{L^\infty}.
\end{equation*}
Using $\|\Delta_j \phi\|_{L^r} \approx 2^{-\alpha j} \|\Lambda^\alpha \phi\|_{L^r}$ with $j\geq0$ and $\|\nabla S_{j-1}\theta\|_{L^\infty} \lesssim 2^j \|\theta\|_{L^\infty}$, with $j \sim k$, it yields
$$\| R_{\alpha} (\Delta_j \phi \cdot \nabla S_{j-1}\theta)\|_{L^r} \lesssim 2^{(2-2\alpha)k} \|\Lambda^\alpha \phi\|_{L^r} \|\theta\|_{L^\infty}.$$
By H\"{o}lder's inequality,
$$\|\Delta_j \phi \cdot \nabla  R_{\alpha} (S_{j-1}\theta)\|_{L^r} \le \|\Delta_j \phi\|_{L^r} \|\nabla  R_{\alpha} (S_{j-1}\theta)\|_{L^\infty}.$$
Since $\nabla  R_{\alpha} $ behaves as an operator of order $2-\alpha > 0$, we get
\begin{equation*}
    \|\nabla  R_{\alpha} (S_{j-1}\theta)\|_{L^\infty} \lesssim \sum_{m \le j-2} 2^{(2-\alpha)m} \|\Delta_m \theta\|_{L^\infty} \lesssim 2^{(2-\alpha)j} \|\theta\|_{L^\infty}.
\end{equation*}
As a result, it gives
\begin{equation*}
    \|I_2\|_{L^r}   \lesssim 2^{(2-2\alpha)k} \|\Lambda^\alpha \phi\|_{L^r} \|\theta\|_{L^\infty}.
\end{equation*}
Finally, for the high-high frequency interactions, output frequencies are mainly $\le 2^{j+1}$, thus $j \ge k-2$. Moreover, since $\phi$ is divergence-free, we can pull the divergence operator outside
\begin{equation*}
    I_3 = \Delta_k \sum_{j \ge k-2} \Big(  R_{\alpha} \nabla \cdot (\Delta_j \phi \widetilde{\Delta}_j \theta) - \nabla \cdot (\Delta_j \phi R_{\alpha} (\widetilde{\Delta}_j \theta)) \Big)
\end{equation*}
Keeping in mind that the operator $\Delta_k  R_{\alpha}  \nabla \cdot$ has total order $2-\alpha$, we are able to conclude
\begin{align*}
    \sum_{j \ge k-2}\|\Delta_k  R_{\alpha}  \nabla \cdot (\Delta_j \phi \widetilde{\Delta}_j \theta)\|_{L^r} &\lesssim \sum_{j \ge k-2}2^{(2-\alpha)k} \|\Delta_j \phi\|_{L^r} \|\widetilde{\Delta}_j \theta\|_{L^\infty} \\&\lesssim \sum_{j \ge k-2}2^{(2-\alpha)k} 2^{-\alpha j} \|\Lambda^\alpha \phi\|_{L^r} \|\theta\|_{L^\infty}\\& = 2^{(2-2\alpha)k}\sum_{j \ge k-2} 2^{-\alpha (j-k)} \|\Lambda^\alpha \phi\|_{L^r} \|\theta\|_{L^\infty}
    \\&\lesssim 2^{(2-2\alpha)k} \|\Lambda^\alpha \phi\|_{L^r} \|\theta\|_{L^\infty}.
\end{align*}
Also, one may derive
\begin{align*}
\sum_{j \ge k-2}\|\Delta_k \nabla \cdot (\Delta_j \phi  R_{\alpha}  (\widetilde{\Delta}_j \theta))\|_{L^r}   &\lesssim \sum_{j \ge k-2}2^k \|\Delta_j \phi\|_{L^r} \| R_{\alpha}  (\widetilde{\Delta}_j \theta)\|_{L^\infty}\\&\lesssim \sum_{j \ge k-2}2^k 2^{-\alpha j} 2^{(2-2\alpha)k} \sum_{j \ge k-2} 2^{(1-2\alpha)(j-k)} \|\Lambda^\alpha \phi\|_{L^r} \|\theta\|_{L^\infty}
\\&\lesssim 2^{(2-2\alpha)k} \|\Lambda^\alpha \phi\|_{L^r} \|\theta\|_{L^\infty},
\end{align*}
where in the last line we have used the assumption $\alpha>\frac{1}{2}$.
Combining the bounds for $I_1, I_2$, and $I_3$, we conclude that
\begin{equation*}
    \|\Delta_k [ R_{\alpha} , \phi\cdot\nabla]\theta\|_{L^{r}} \leq C 2^{(2-2\alpha)k}\|\Lambda^{\alpha} \phi\|_{ L^{r}}\|\theta\|_{L^{\infty}}
\end{equation*}
This completes the proof of the lemma.
\end{proof}

\begin{lemma}[Full-frequency trilinear convective commutator estimate]
 	\label{lem:Q-general-v}
 	Let $\frac12<\alpha<1$, $0<\beta<\alpha$, and $1<r<\infty$.
 	For every divergence-free vector field $v$ and every $j\geq-1$,
 	\begin{equation}\label{esdadfdq1}
 		\|\Delta_j\Qop(v;\theta,\theta)\|_{L^r}
 		\lesssim 2^{(3-4\alpha)j}
 		\|\Lambda^\alpha v\|_{L^r}\|\theta\|_{L^\infty}^2.
 	\end{equation}
 	In particular, we have
 	\begin{align}
 		\|\Delta_j\Qop(u_H;\theta,\theta)\|_{L^r}
 		&\lesssim 2^{(3-4\alpha)j}
 		\|\Lambda^\alpha u_H\|_{L^r}\|\theta\|_{L^\infty}^2,
 		\label{eq:delta-bound}\\
 		\|\Delta_j\Qop(\Top\theta;\theta,\theta)\|_{L^r}
 		&\lesssim 2^{(3-4\alpha)j}
 		\|\Lambda^\alpha\Top\theta\|_{L^r}\|\theta\|_{L^\infty}^2.
 		\label{eq:delta-T-bound}
 	\end{align}
 	The implicit constants may depend on $\alpha,\beta,r$ and the
 	fixed cutoffs, but not on $j$ or the input functions.
 \end{lemma}

 \begin{proof}
 	Put $w=\Lambda^\alpha v$, $W=\|w\|_{L^r}$ and
 	$A=\|\theta\|_{L^{\infty}}$. On homogeneous velocity annuli,
 	$\|v_U\|_{L^r}\lesssim U^{-\alpha}W$.
 	The dyadic estimates \eqref{eq:LP-Q-low} and
 	\eqref{eq:LP-Q-high}, proved for $\frac{1}{2}<\alpha<1$, imply
 	for $N=2^j\geq1$ that
 	\begin{align*}
 		\|\Delta_N\Qop(v;\theta,\theta)\|_{L^r}
 		&\lesssim A^2W
 		\left[
 		N^{2-3\alpha}\sum_{U\ll N}U^{1-\alpha}
 		+N^{2-2\alpha}\sum_{U\gtrsim N}U^{1-2\alpha}
 		+N\sum_{U\gtrsim N}U^{2-4\alpha}
 		\right]\\
 		&\lesssim A^2W N^{3-4\alpha}.
 	\end{align*}

 	For $j=-1$, the two shifted bilinear terms vanish by the
 	fixed output cutoff. The unshifted term has internal
 	bilinear output comparable to a velocity annulus
 	$U\geq c_0>0$. Using incompressibility, the integrable kernel
 	of $\Delta_{-1}\operatorname{div}$, and Lemma~\ref{lem:B},
 	\[
 	\|\Delta_{-1}\Qop(v;\theta,\theta)\|_{L^r}
 	\lesssim A^2W\sum_{U\geq c_0}U^{2-4\alpha}
 	\lesssim A^2W.
 	\]
 	Thus \eqref{esdadfdq1} holds for every $j\geq-1$.
 	The two stated applications follow from the divergence-free
 	property of $u_H$ and $\Top\theta$.
 \end{proof}

With the above lemmas in hand, we are in a position to derive a higher regularity estimate of the combined quantity $H$.
\begin{proposition}\label{pr027}
Let $\beta=1-\alpha$ and $\frac{\sqrt{57}-5}{4}<\alpha\leq\frac{10}{13}$. Suppose $H$ admits the following bound
$$\sup_{0\leq t\leq T}\|H(t)\|_{L^{m}}<\infty, \qquad m>\frac{2}{\alpha},$$
where we may assume $m<\frac{2}{1-\alpha}$.
Then it holds
$$\sup_{0\leq t\leq T}\|H(t)\|_{B_{\frac{2m}{2-(1-\alpha)m},\infty}^{3\alpha-2}}<\infty.$$
\end{proposition}
\begin{proof}
First, we separate the lowest output block. The fixed cutoffs give
\[
 \Delta_{-1} R_{\alpha} \theta=0,\qquad
 \Delta_{-1}\Bop(\theta,\theta)=0.
\]
By Bernstein's inequality and the velocity energy bound,
\[
 \|\Delta_{-1}H(t)\|_{L^r}
 =\|\Delta_{-1}\omega(t)\|_{L^r}
 \lesssim\|u(t)\|_2\leq C(T,u_0,\theta_0).
\]
Every use of the annular dissipative lower bound below is
restricted to $k\geq0$.
We denote
\begin{align*}
\Upsilon\triangleq&
[ R_{\alpha} ,u_H\cdot\nabla]\theta+\Cres(\theta,\theta)
+[ R_{\alpha} ,u_B\cdot\nabla]\theta +\Qop(u_H;\theta,\theta)
+\Qop(\Top\theta;\theta,\theta)
\nonumber\\
&+\Qop(u_B;\theta,\theta)+F_{\mathrm{low}},
\end{align*}
then \eqref{eq:H-equation} reduces to
\begin{align}\label{tspeq268}
(\Dt+\Lambda^\alpha)H=\Upsilon.
\end{align}
Applying $\Delta_{k}$ to (\ref{tspeq268}) gives
\begin{eqnarray}\label{Cri005}
\partial_{t}\Delta_{k}H+\Lambda^{\alpha}\Delta_{k}H=
\Delta_{k}\Upsilon-\Delta_{k}(u\cdot\nabla
H).
\end{eqnarray}
For simplicity, we set
$$r=\frac{2m}{2-(1-\alpha)m}.$$
Multiplying (\ref{Cri005}) by
$|\Delta_{k}H|^{r-2}\Delta_{k}H$, integrating the result over space
$\mathbb{R}^{2}$ and using the divergence-free condition, we get
\begin{eqnarray}\label{Cri006}
\frac{1}{r}\frac{d}{dt}\|\Delta_{k}H(t)\|_{L^{r}}^{r}
+\int_{\mathbb{R}^{2}}
(\Lambda^{\alpha}\Delta_{k}H)|\Delta_{k}H|^{r-2}\Delta_{k}H\,dx=I_{1}^{k}
+I_{2}^{k},
\end{eqnarray}
where
$$I_{1}^{k}=-\int_{\mathbb{R}^{2}}
{\Delta_{k}(u\cdot\nabla H)\,\,|\Delta_{k}H|^{r-2}\Delta_{k}H\,dx},$$
$$I_{2}^{k}=\int_{\mathbb{R}^{2}}
{\Delta_{k}\Upsilon\,\,|\Delta_{k}H|^{r-2}\Delta_{k}H\,dx}.$$
We recall the following lower bound
$$\int_{\mathbb{R}^{2}}
(\Lambda^{\alpha}\Delta_{k}H)|\Delta_{k}H|^{r-2}\Delta_{k}H\,dx\geq c2^{\alpha k}\|\Delta_{k}H\|_{L^{r}}^{r},\quad k\geq0 $$
with an absolute constant $c>0$ independent of $k$.
Keeping in mind $r=\frac{2m}{2-(1-\alpha)m}$, we are able to conclude
\begin{align}\label{Cri007}
 \|\Lambda^{\alpha}u\|_{L^{r}} &\leq
\|\Lambda^{\alpha}u_{H}\|_{L^{r}}+\|\Lambda^{\alpha}u_{B}\|_{L^{r}}
+\|\Lambda^{\alpha}\Top\theta\|_{L^{r}}
\nonumber\\&\leq
C+C\|\Lambda^{\alpha-1}H\|_{L^{r}}+C\| \theta\|_{L^{r}}
\nonumber\\&\leq
C+C\|H\|_{L^{m}} +C\|\theta\|_{L^{r}}
\nonumber\\&\leq C(T,\,u_{0},\,\theta_{0}).
\end{align}
where we have used the following fact due to \eqref{eq:B-dyadic}
\begin{align}
\|\Lambda^{\alpha}u_{B}\|_{L^{r}}&\lesssim \|\Lambda^{\alpha-1}\Bop(\theta,\theta)\|_{L^{r}}\nonumber\\
&\lesssim \sum_{j\geq0}\|\Lambda^{\alpha-1}\Delta_j\Bop(\theta,\theta)\|_{L^{r}}
\nonumber\\
&\lesssim \sum_{j\geq0}2^{j(\alpha-1)}\|\Delta_j\Bop(\theta,\theta)\|_{L^{r}}
\nonumber\\
&\lesssim \sum_{j\geq0}2^{j(\alpha-1)}\norm{\theta}{L^{\infty}}\,2^{j(2-3\alpha)}
\sum_{v\geq j-C}2^{-\alpha(v-j)}
\norm{\widetilde\Delta_v\theta}{L^{r}}\nonumber\\
&\lesssim \sum_{j\geq0}2^{j(\alpha-1)}\norm{\theta}{L^{\infty}}\,2^{j(2-3\alpha)}
\norm{\theta}{L^{r}}
\nonumber\\
&\lesssim \sum_{j\geq0} 2^{-j(2\alpha-1)}
\nonumber\\ &\leq C.\nonumber
\end{align}
Noticing \eqref{Cri007}, it follows from the proof of the estimate (7.17) in \cite{JMWZ} that
\begin{align}\label{Cri009}
|I_{1}^{k}| \leq& C\|\Lambda^{\alpha}u\|_{L^{r}}2^{(1-\alpha+\frac{2}{r})k}\|\Delta_{k}
H\|_{L^{r}}^{r-1}\Big(\|\Delta_{k}H\|_{L^{r}}+
\sum_{l\leq k-1}2^{(1+\frac{2}{r})(l-k)}\|\Delta_{l}H\|_{L^{r}}\nonumber\\&
+
\sum_{l\geq k-1}2^{(\alpha-\frac{2}{r})(k-l)}\|\Delta_{l}H\|_{L^{r}}\Big)\nonumber\\
 \leq&
C2^{(1-\alpha+\frac{2}{r})k}\|\Delta_{k}
H\|_{L^{r}}^{r-1}\Big(\|\Delta_{k}H\|_{L^{r}}+
\sum_{l\leq k-1}2^{(1+\frac{2}{r})(l-k)}\|\Delta_{l}H\|_{L^{r}}\nonumber\\&
+
\sum_{l\geq k-1}2^{(\alpha-\frac{2}{r})(k-l)}\|\Delta_{l}H\|_{L^{r}}\Big).
\end{align}
It follows from \eqref{eq:mixed_direct} and \eqref{Cri007} that
\begin{align}\label{sdfgh01}
\norm{\Delta_k [ R_{\alpha} ,u_{H}\cdot\nabla]\theta}{L^{r}}+\norm{\Delta_k [ R_{\alpha} ,u_{B}\cdot\nabla]\theta}{L^{r}}&\lesssim 2^{(2-2\alpha)k} (\norm{\Lambda^{\alpha}u_{H}}{L^{r}} +\norm{\Lambda^{\alpha}u_{B}}{L^{r}} )\norm{\theta}{L^{\infty}}\nonumber\\& \lesssim 2^{(2-2\alpha)k}.
\end{align}
If \eqref{eq:res-condition} holds, then we derive from \eqref{eq:res-bound} that
\begin{equation}\label{sdfgh02}
\norm{\Lambda^{-\frac{\alpha}{2}}\Cres(\theta,\theta)}{L^r}
\lesssim\norm{\theta}{L^{\infty}}\norm{\theta}{L^r}\leq C.
\end{equation}
Thanks to \eqref{eq:delta-bound} and \eqref{Cri007}, we have
\begin{equation}\label{sdfgh03}
\norm{ \Delta_{k}\Qop(u_H;\theta,\theta)}{L^{r}} \lesssim 2^{(3-4\alpha)k} \norm{\Lambda^\alpha u_H}{L^{r}} \norm{\theta}{L^{\infty}}^{2}\lesssim 2^{(3-4\alpha)k}.
\end{equation}
Making use of \eqref{eq:delta-T-bound} and \eqref{Cri007}, one also has
\begin{equation}\label{sdfgh04}
\norm{\Delta_k \Qop(\Top\theta;\theta,\theta)}{L^{r}}
\lesssim 2^{(3-4\alpha)k} \norm{\Lambda^\alpha \Top\theta}{L^{r}}\norm{\theta}{L^{\infty}}^2 \lesssim 2^{(3-4\alpha)k}.
\end{equation}
By \eqref{eq:Q-quartic} and \eqref{eq:upper-quartic}, one gets
\begin{equation}\label{sdfgh05}
\norm{\Lambda^{-\frac{\alpha}{2}}\Qop(u_B;\theta,\theta)}{L^{r}}
\lesssim \norm{\theta}{L^{\infty}}^3\norm{\theta}{L^{r}}\leq C.
\end{equation}
We deduce from \eqref{eq:low-bound} that
\begin{equation}\label{sdfgh06}
\norm{\Lambda^{-\frac{\alpha}{2}}F_{\mathrm{low}}}{L^r}
\lesssim (1+\norm{\theta}{L^\infty})\norm{\theta}{L^r}\leq C,
\end{equation}
Combining \eqref{sdfgh01}, \eqref{sdfgh02}, \eqref{sdfgh03}, \eqref{sdfgh04}, \eqref{sdfgh05}, \eqref{sdfgh06} and using the Bernstein inequality, we derive
\begin{align}
\|\Delta_{k}\Upsilon\|_{L^{r}}  \lesssim 2^{(2-2\alpha)k}+2^{\frac{\alpha k}{2}}
+2^{(3-4\alpha)k}\lesssim 2^{2(1-\alpha)k},\nonumber
\end{align}
Finally, it follows from the Bernstein inequality
\begin{align}\label{Cri010}
|I_{2}^{k}| \lesssim&
 \| \Delta_{k}\Upsilon\|_{L^{r}}
\|\Delta_{k}H\|_{L^{r}}^{r-1}
\nonumber\\  \lesssim&
 2^{2(1-\alpha)k}\|\Delta_{k}H\|_{L^{r}}^{r-1}.
\end{align}
Inserting \eqref{Cri009} and \eqref{Cri010} into (\ref{Cri006}) yields
\begin{align}\label{Cri011}
\frac{d}{dt}\|\Delta_{k}H(t)\|_{L^{r}}
+c2^{\alpha k}\|\Delta_{k}H\|_{L^{r}}\leq C2^{2(1-\alpha)k}+  C2^{(1-\alpha+\frac{2}{r})k}L(t),
\end{align}
where $L(t)$ is given by
\begin{eqnarray}
L(t)\triangleq\|\Delta_{k}H(t)\|_{L^{r}}+
\sum_{l\leq k-1}2^{(1+\frac{2}{r})(l-k)}\|\Delta_{l}H(t)\|_{L^{r}}
+\sum_{l\geq k-1}2^{(\alpha-\frac{2}{r})(k-l)}\|\Delta_{l}H(t)\|_{L^{r}}.\nonumber
\end{eqnarray}
Integrating (\ref{Cri011}) in time, we have
\begin{align}\label{Cri012}
\|\Delta_{k}H(t)\|_{L^{r}} \leq&  e^{-ct2^{\alpha k}}\|\Delta_{k}H_{0}\|_{L^{r}}+ C2^{(2-3\alpha)k}  +C2^{(1-\alpha+\frac{2}{r})k}\int_{0}^{t}{ e^{-c(t-\tau)2^{\alpha k}}L(\tau) \,d\tau}.
\end{align}
Multiplying (\ref{Cri012}) by $2^{sk}$ with $s\triangleq3\alpha-2$ and taking sup with respect to $k$, it is obvious to see that
$$\|H(t)\|_{B_{r,\infty}^{s}}\leq C+M_{1}+M_{2}+M_{3},$$
where
$$M_{1}=C\sup_{k\geq-1}\Big(2^{(1-\alpha+\frac{2}{r})k}2^{sk}\int_{0}^{t}{ e^{-c(t-\tau)2^{\alpha k}}\|\Delta_{k}H(\tau)\|_{L^{r}} \,d\tau}\Big),$$
$$M_{2}=C\sup_{k\geq-1}\Big(2^{(1-\alpha+\frac{2}{r})k}2^{sk}\int_{0}^{t}{ e^{-c(t-\tau)2^{\alpha k}}\sum_{l\leq k-1}2^{(1+\frac{2}{r})(l-k)}\|\Delta_{l}H(\tau)\|_{L^{r}} \,d\tau}\Big),$$
$$M_{3}=C\sup_{k\geq-1}\Big(2^{(1-\alpha+\frac{2}{r})k}2^{sk}\int_{0}^{t}{ e^{-c(t-\tau)2^{\alpha k}}\sum_{l\geq k-1}2^{(\alpha-\frac{2}{r})(k-l)}\|\Delta_{l}H(\tau)\|_{L^{r}} \,d\tau}\Big).$$
Thanks to $m>\frac{2}{\alpha}$, it is easy to see
 $$r\equiv \frac{2m}{2-(1-\alpha)m}>\frac{2}{2\alpha-1}.$$
This allows us to choose $\epsilon$ such that
$$0<\epsilon<2\alpha-1-\frac{2}{r}.$$
According to the Bernstein inequality and the convolution Young inequality, we derive
\begin{align}
M_{1} =&
C\sup_{k\geq -1
}\Big(2^{(1-\alpha+\frac{2}{r}+\epsilon)k}\int_{0}^{t}{
e^{-c(t-\tau)2^{\alpha k}}2^{(s-\epsilon)k}\|\Delta_{k}H(\tau)\|_{L^{r}}
\,d\tau}\Big)
\nonumber\\ \leq&
C\sup_{k\geq -1
}\Big(2^{(1-\alpha+\frac{2}{r}+\epsilon)k}\int_{0}^{t}{
e^{-c(t-\tau)2^{\alpha k}}\|H(\tau)\|_{B_{r,\infty}^{s-\epsilon}}
\,d\tau}\Big)
\nonumber\\ \leq&
C\sup_{k\geq-1
}2^{-(2\alpha-1-\frac{2}{r}-\epsilon)k}\|H\|_{L^{\infty}(0,\, t;B_{r,\infty}^{s-\epsilon})}\nonumber\\ \leq&
C\sup_{k\geq-1
}2^{-(2\alpha-1-\frac{2}{r}-\epsilon)k}\|H\|_{L^{\infty}(0,\, T;B_{r,\infty}^{s-\epsilon})}\nonumber\\ \leq&
C \|H\|_{L^{\infty}(0,\, T;B_{r,\infty}^{s-\epsilon})}.\nonumber
\end{align}
Similarly, it follows that
\begin{align}
M_{2} =& C\sup_{k\geq-1}\Big(2^{(1-\alpha+\frac{2}{r}+\epsilon)k} \int_{0}^{t}{ e^{-c(t-\tau)2^{\alpha k}}\sum_{l\leq k-1}2^{(1+\frac{2}{r}-s+\epsilon)(l-k)}2^{(s-\epsilon)l}\|\Delta_{l}H(\tau)\|_{L^{r}} \,d\tau}\Big) \nonumber\\ \leq&
C\sup_{k\geq-1}\Big(2^{(1-\alpha+\frac{2}{r}+\epsilon)k}\int_{0}^{t}{
e^{-c(t-\tau)2^{\alpha k}}\Big(\sum_{l\leq
k-1}2^{(1+\frac{2}{r}-s+\epsilon)(l-k)}\Big)\|H(\tau)\|_{B_{r,\infty}^{s-\epsilon}}
\,d\tau}\Big)
\nonumber\\ \leq&C\sup_{k\geq-1}\Big(2^{(1-\alpha+\frac{2}{r}+\epsilon)k}\int_{0}^{t}{
e^{-c(t-\tau)2^{\alpha k}} \|H(\tau)\|_{B_{r,\infty}^{s-\epsilon}}
\,d\tau}\Big)\nonumber\\ \leq&
C \|H\|_{L^{\infty}(0,\, T;B_{r,\infty}^{s-\epsilon})}.\nonumber
\end{align}
Regarding the term $M_{3}$, one also has
\begin{align}
M_{3} =&C\sup_{k\geq-1}\Big(2^{(1-\alpha+\frac{2}{r}+\epsilon)k}\int_{0}^{t}{ e^{-c(t-\tau)2^{\alpha k}}\sum_{l\geq k-1}2^{(\alpha-\frac{2}{r}+s-\epsilon)(k-l)}2^{(s-\epsilon)l}\|\Delta_{l}H(\tau)\|_{L^{r}} \,d\tau}\Big) \nonumber\\ \leq&C\sup_{k\geq-1}\Big(2^{(1-\alpha+\frac{2}{r}+\epsilon)k}\int_{0}^{t}{ e^{-c(t-\tau)2^{\alpha k}}\sum_{l\geq k-1}2^{(\alpha-\frac{2}{r}+s-\epsilon)(k-l)}\|H(\tau)\|_{B_{r,\infty}^{s-\epsilon}}\,d\tau}\Big)
\nonumber\\ \leq&C\sup_{k\geq-1}\Big(2^{(1-\alpha+\frac{2}{r}+\epsilon)k}\int_{0}^{t}{
e^{-c(t-\tau)2^{\alpha k}} \|H(\tau)\|_{B_{r,\infty}^{s-\epsilon}}
\,d\tau}\Big)\nonumber\\ \leq&
C \|H\|_{L^{\infty}(0,\, T;B_{r,\infty}^{s-\epsilon})}.\nonumber
\end{align}
Therefore, we obtain
\begin{align}\label{OKE3}
\|H\|_{L^{\infty}(0,\, T;B_{r,\infty}^{s})} \leq C
+
C \|H\|_{L^{\infty}(0,\, T;B_{r,\infty}^{s-\epsilon})}.
\end{align}
Using the Bernstein inequality, we have for $\epsilon>0$ that
\begin{align}\label{OKE4}
C\|H\|_{B_{r,\infty}^{s-\epsilon}}\leq& C\sup_{-1\leq j\leq L}2^{(s-\epsilon)j}\|\Delta_{j}H\|_{L^{r}}
+C\sup_{ j\geq L+1}2^{(s-\epsilon)j}\|\Delta_{j}H\|_{L^{r}}\nonumber\\
\leq& C\sup_{-1\leq j\leq L}2^{(s-\epsilon)j}2^{2j(\frac{1}{2}-\frac{1}{r})}
\|\Delta_{j}H\|_{L^{2}}
+C\sup_{ j\geq L+1}2^{-j\epsilon}2^{js}\|\Delta_{j}H\|_{L^{r}}
\nonumber\\
\leq&  C\sup_{-1\leq j\leq L}2^{j(2+s-\epsilon-\frac{2}{r})}
+C\sup_{ j\geq L+1}2^{-j\epsilon}\|H\|_{B_{r,\infty}^{s}}\nonumber\\
\leq&  C2^{L(2+s-\epsilon-\frac{2}{r})}
+C2^{-L\epsilon} \|H\|_{B_{r,\infty}^{s}},
\end{align}
where we have used the following fact due to \eqref{eq:B-dyadic}
\begin{align}
\|\Delta_{j}H\|_{L^{2}}\leq&  \|\Delta_{j}\omega\|_{L^{2}}+\|\Delta_{j} R_{\alpha} \theta\|_{L^{2}}
+\|\Delta_{j}\Bop(\theta,\theta)\|_{L^{2}}\nonumber\\
\leq& C2^{j}\|u\|_{L^{2}}+ C2^{(1-\alpha)j}\|\theta\|_{L^{2}}
+ C2^{(2-3\alpha)j}\|\theta\|_{L^{2}}\nonumber\\
\leq& C2^{j}.\nonumber
\end{align}
Now fixing $L$ as
$$\frac{1}{100}\leq C2^{-L\epsilon}\leq \frac{1}{2},$$
we deduce from \eqref{OKE4} that
\begin{align}\label{OKEsd671}
C\|H\|_{B_{r,\infty}^{s-\epsilon}}\leq  C
+\frac{1}{2}\|H\|_{B_{r,\infty}^{s}}.
\end{align}
Inserting \eqref{OKEsd671} into \eqref{OKE3} gives
\begin{align}\label{OKE5}
\|H\|_{L^{\infty}(0,\, T;B_{r,\infty}^{s})}\leq& C+\|H_{0}\|_{B_{r,\infty}^{s}}+
\frac{1}{2}\|H\|_{L^{\infty}(0,\, T;B_{r,\infty}^{s})}.
\end{align}
Therefore, it follows from \eqref{OKE5} that
$$\sup_{0\leq t\leq T}\|H(t)\|_{B_{r,\infty}^{s}}\leq C(T,u_{0},\theta_{0}). $$
This completes the proof of Proposition \ref{pr027}.
\end{proof}

\begin{proposition}[The range of $\alpha$]\label{dflkeuyr11}
To ensure the validity of Proposition \ref{prop:energy}, Proposition \ref{prop:upper-critical-energy} and Proposition \ref{pr027}, in the critical case $\alpha+\beta=1$, it requires
\begin{eqnarray}\label{sdcpkh98}
\frac{\sqrt{57}-5}{4}<\alpha\leq \frac{10}{13}.
\end{eqnarray}
\end{proposition}
\begin{proof}
To ensure the validity of both Proposition \ref{prop:energy} and Proposition \ref{pr027}, it requires
\begin{equation}
\frac8{13}<\alpha<\frac23,\qquad \frac{4-4\alpha-3\alpha^{2}}{\alpha}<1-\alpha<\alpha,\nonumber
\end{equation}
and
\begin{align} \label{askop99}
\max\left\{A(\alpha,1-\alpha),\, \frac{1-\alpha}{2}\right\}<\frac{1}{m}<\frac{\alpha}{2}.
\end{align}
where $A(\alpha,1-\alpha)$ is given by
\begin{align*}
\max\left\{\frac{3(2+\alpha)(2-3\alpha)-\alpha(1-\alpha)}{4(1-\alpha)},\ \frac{(4+2\alpha)[(\frac{2}{\alpha}-\frac{3}{2})(1-\alpha)-\alpha]-\alpha(1-\alpha)}
{4
(1-\alpha)}
\right\}.
\end{align*}
It is not hard to check that the $m$ of \eqref{askop99} is valid as long as
$$\frac{\sqrt{57}-5}{4}<\alpha<\frac{2}{3}.$$

When $\alpha\in[\frac{2}{3},\,\frac{10}{13}]$, we have that
$m=\frac{13}{4}>\frac{2}{\alpha}$, which ensures the validity of both Proposition \ref{prop:upper-critical-energy} and Proposition \ref{pr027}. Therefore, the lower bound of $\alpha$ stated in \eqref{sdcpkh98} follows immediately.
\end{proof}

\section{The proof of Theorem \ref{thm:critical}}
\label{secscgd9871}
We are ready to prove Theorem \ref{thm:critical}. For $\alpha>\frac{10}{13}$ the result is contained in \cite{WXXY}; it therefore suffices to treat $\alpha_0<\alpha\le\frac{10}{13}$.
It follows from Proposition \ref{dflkeuyr11} (consider the case $\alpha_0<\alpha\leq\frac{10}{13}$ only) that for some $m>\frac{2}{\alpha}$,
$H$ satisfies
\begin{eqnarray}\label{Yes01}
\sup_{0\leq t\leq T}\|H(t)\|_{B_{\frac{2m}{2-(1-\alpha)m},\,\infty}
^{3\alpha-2}}\leq C(T,\,u_{0},\,\theta_{0})<\infty.
\end{eqnarray}
Set $\sigma\triangleq 2\alpha-\frac2m$. Since $m>\frac2\alpha$ we have $\sigma>\alpha$, and since $m<\frac{2}{1-\alpha}\le\frac{2}{2\alpha-1}$ (resp.\ $m=\frac{13}{4}<\frac{2}{2\alpha-1}$ on $[\frac23,\frac{10}{13}]$), we have $\sigma<1$. Thus, making use of \eqref{Yes01}, we have
\begin{align}
\| u_{H}\|_{C^{\sigma}}&=  \|\nabla^{\perp}
\Delta^{-1}H\|_{C^{\sigma}}\nonumber\\&\approx
\|\nabla^{\perp}\Delta^{-1}H\|_{B_{\infty,\infty}^{\sigma}}
\nonumber\\&\leq C+C\|H\|_{B_{\infty,\infty}^{\sigma-1}}
\nonumber\\& \leq
C +C\|H\|_{B_{\frac{2m}{2-(1-\alpha)m},\,\infty}^{3\alpha-2}} \nonumber\\&\leq
C(T,\,u_{0},\,\theta_{0}).\nonumber
\end{align}
This yields the global $C^{\sigma}$ bound on $ u_{H}$. According to the bound $\|\theta(t)\|_{L^{\infty}}\leq
\|\theta_{0}\|_{L^{\infty}}$, we
 conclude
\begin{eqnarray}\|\Top\theta\|_{C^{\alpha}}
\leq  C +C\|\theta\|_{L^{\infty}}\leq
C(T,\,u_{0},\,\theta_{0}).\nonumber
\end{eqnarray}
Recalling \eqref{eq:uB-holder}, we have
\begin{align}
\norm{u_B}{B_{\infty,\infty}^{3\alpha-1}}\leq
C \norm{\theta}{L^{\infty}}^2\leq
C(T,\,u_{0},\,\theta_{0}).\nonumber
\end{align}
Combining the above three bounds gives
\begin{align}
\|u\|_{C^{\alpha}}&\leq \|u_{H}\|_{C^{\alpha}}+\|\Top\theta\|_{C^{\alpha}}+\|u_{B}\|_{C^{\alpha}}
\nonumber\\&\leq\|u_{H}\|_{C^{\sigma}}+\|\Top\theta\|_{C^{\alpha}}+C\norm{u_B}{B_{\infty,\infty}^{3\alpha-1}} \nonumber\\&\leq C(T,\,u_{0},\,\theta_{0}).\nonumber
\end{align}

Fix $t_0>0$ inside the local existence interval, depending only
on the initial $H^s$ norms. Local well-posedness gives
\[
 \int_0^{\min\{T,t_0\}}\|\nabla\theta(t)\|_{L^{\infty}}\,dt
 \leq C(u_0,\theta_0,t_0).
\]
For $T>t_0$, apply the positive-time H\"older estimate in
\cite[Theorem 1.1]{Silvestre} to
\[
 \partial_t\theta+u\cdot\nabla\theta+\Lambda^\beta\theta=0.
\]
For some $\eta>0$,
\[
 \sup_{\frac{t_0}{2}\leq t\leq T}\|\theta(t)\|_{C^\eta}
 \leq C(T,t_0,u_0,\theta_0).
\]
Choose $0<\eta_0<\min\{\eta,1-\alpha\}$. Since $\Top$ has
order $-\alpha$ and vanishes near zero,
\[
 \sup_{\frac{t_0}{2}\leq t\leq T}
       \|\Top\theta(t)\|_{C^{\alpha+\eta_0}}
 \leq C(T,t_0,u_0,\theta_0).
\]
We take
\[
 1-\beta\equiv\alpha<\varrho<\min\{\sigma,\,3\alpha-1,\,\alpha+\eta_0\}.
\]
The velocity decomposition yields
\[
 \sup_{\frac{t_0}{2}\leq t\leq T}\|u(t)\|_{C^\varrho}
 \leq C(T,t_0,u_0,\theta_0).
\]
Now take $0<\delta<\min\{\varrho-(1-\beta),\beta\}$.
The differentiability estimate (see
\cite{Sil20121,XueYepjm,Schauder}), with diffusion exponent
$\frac{\beta}{2}$ and zero source, applies on time cylinders contained
in $[\frac{t_0}{2},T]$. A covering by cylinders of fixed length depending
on $t_0$, we thus derive
\[
 \sup_{t_0\leq t\leq T}\|\nabla\theta(t)\|_{L^{\infty}}
 \leq C(T,t_0,u_0,\theta_0).
\]
Combining this with the local estimate, and recalling that
$t_0$ is fixed by the initial data, gives
\begin{eqnarray}\label{KEY1}
\int_{0}^{T}{\|\nabla \theta(t)\|_{L^{\infty}}\,dt}\leq C(T,\,u_{0},\,\theta_{0}),\end{eqnarray}
and
\begin{eqnarray}\label{KEY2}\|\omega(t)\|_{L^{\infty}}\leq \|\omega_{0}\|_{L^{\infty}}+\int_{0}^{t}{\|\nabla \theta(\tau)\|_{L^{\infty}}\,d \tau}<\infty.
\end{eqnarray}
The standard energy method allows us to derive
\begin{align}\label{cvbmr82}
&\frac{d}{dt}(\|\Lambda^{s}u(t)\|_{L^{2}}^{2}+
\|\Lambda^{s}\theta(t)\|_{L^{2}}^{2})+\|\Lambda^{s+\frac{\alpha}{2}}
u\|_{L^{2}}^{2}+\|\Lambda^{s+\frac{\beta}{2}}\theta\|_{L^{2}}^{2}\nonumber\\
&\leq C (1+\|\nabla u\|_{L^{\infty}}+\|\nabla \theta\|_{L^{\infty}})(\|\Lambda^{s}u\|_{L^{2}}^{2}+
\|\Lambda^{s}\theta\|_{L^{2}}^{2}).
\end{align}
Recalling the following logarithmic Sobolev embedding inequality
\begin{eqnarray}
\|\nabla u\|_{L^{\infty}}
\leq C\Big(1+\|u\|_{L^{2}}+
\|\omega\|_{L^{\infty}} \ln\big(e+\|
\Lambda^{s}u\|_{L^{2}}\big)\Big),\nonumber
\end{eqnarray}
we deduce from \eqref{cvbmr82} that
\begin{align}\label{cvbmr83}
&\frac{d}{dt}(\|\Lambda^{s}u(t)\|_{L^{2}}^{2}+
\|\Lambda^{s}\theta(t)\|_{L^{2}}^{2})+\|\Lambda^{s+\frac{\alpha}{2}}
u\|_{L^{2}}^{2}+\|\Lambda^{s+\frac{\beta}{2}}\theta\|_{L^{2}}^{2}\nonumber\\
&\leq C (1+\|\omega\|_{L^{\infty}}+\|\nabla \theta\|_{L^{\infty}})\ln(e+\|\Lambda^{s}u\|_{L^{2}}^{2}+
\|\Lambda^{s}\theta\|_{L^{2}}^{2})(\|\Lambda^{s}u\|_{L^{2}}^{2}+
\|\Lambda^{s}\theta\|_{L^{2}}^{2}).
\end{align}
Keeping in mind (\ref{KEY1}) and (\ref{KEY2}), one obtains by applying the Log-Gronwall inequality to \eqref{cvbmr83}
\begin{eqnarray}
\|\Lambda^{s}u(t)\|_{L^{2}}^{2}+
\|\Lambda^{s}\theta(t)\|_{L^{2}}^{2}+\int_{0}^{t}{(\|\Lambda^{s+\frac{\alpha}{2}}
u\|_{L^{2}}^{2}+\|\Lambda^{s+\frac{\beta}{2}}\theta\|_{L^{2}}^{2})(\tau)\,d\tau}\leq C(t,\,u_{0},\,\theta_{0}).\nonumber
\end{eqnarray}
This completes the proof of Theorem \ref{thm:critical}.

\vskip .3in
\appendix

\section{The proof of the localized kernel principle (Lemma \ref{lem:kernel})}
\label{app:Q-proof11}

\begin{proof}
The proof is divided into six systematic steps.

\medskip
\noindent\textbf{Step 1:} \emph{Multilinear spatial convolution representation.} \\
Let $\xi = (\xi_1, \ldots, \xi_r) \in (\mathbb R^2)^r = \mathbb R^d$. By Fourier inversion, the action of the multilinear multiplier operator $T_{\mathfrak{m}}$ on Schwartz functions $f_1, \ldots, f_r \in \mathcal{S}(\mathbb R^2)$ is defined pointwise by
\begin{align}\label{eq:spatial-integral-rep}
T_{\mathfrak{m}}(f_1, \ldots, f_r)(x) &= \frac{1}{(2\pi)^{2r}} \int_{(\mathbb R^2)^r} e^{ix \cdot (\xi_1 + \cdots + \xi_r)} \mathfrak{m}(\xi_1, \ldots, \xi_r) \prod_{\iota=1}^r \widehat{f}_\iota(\xi_\iota) \,d\xi_1 \cdots d\xi_r \nonumber\\
&= \int_{(\mathbb R^2)^r} \Psi
(y_1, \ldots, y_r) \prod_{\iota=1}^r f_\iota(x - y_\iota) \,dy_1 \cdots dy_r,
\end{align}
where the convolution kernel $\Psi \in \mathcal{S}(\mathbb R^d)$ is given by the inverse Fourier transform of $\mathfrak{m}$ on $\mathbb R^d$:
\begin{equation}
\Psi(Y) \triangleq (\mathcal{F}_{\mathbb R^d}^{-1} \mathfrak{m})(Y) = \frac{1}{(2\pi)^d} \int_{\mathbb R^d} e^{i Y \cdot \xi} \mathfrak{m}(\xi) \,d\xi, \qquad Y = (y_1, \ldots, y_r) \in (\mathbb R^2)^r = \mathbb R^d.\nonumber
\end{equation}

\medskip
\noindent\textbf{Step 2:} \emph{Invariance of the $L^1$ norm under linear coordinate changes.} \\
By hypothesis, there exists an invertible linear map $\mathrm{L} \in \mathrm{GL}(d, \mathbb R)$ such that the adapted coordinates $\zeta = \mathrm{L}\xi$ decouple the frequency interactions. We define the relocated symbol on $\mathbb R^d$ by
\[
\widetilde{\mathfrak{m}}(\zeta) \triangleq (\mathfrak{m} \circ L^{-1})(\zeta) = \mathfrak{m}(L^{-1}\zeta), \qquad \zeta \in \mathbb R^d.
\]
Let $\widetilde{\Psi} \triangleq \mathcal{F}_{\mathbb R^d}^{-1} \widetilde{\mathfrak{m}}$ denote the inverse Fourier transform of $\widetilde{\mathfrak{m}}$. Using the change of variables $\zeta = \mathrm{L}\xi$ (where $d\zeta = |\det \mathrm{L}| \,d\xi$):
\begin{align*}
\Psi(Y) &= \frac{1}{(2\pi)^d} \int_{\mathbb R^d} e^{i Y \cdot \xi} \widetilde{\mathfrak{m}}(\mathrm{L}\xi) \,d\xi \\
&= \frac{1}{(2\pi)^d |\det \mathrm{L}|} \int_{\mathbb R^d} e^{i Y \cdot (\mathrm{L}^{-1}\zeta)} \widetilde{\mathfrak{m}}(\zeta) \,d\zeta \\
&= \frac{1}{(2\pi)^d |\det \mathrm{L}|} \int_{\mathbb R^d} e^{i (\mathrm{L}^{-T} Y) \cdot \zeta} \widetilde{\mathfrak{m}}(\zeta) \,d\zeta \\
&= \frac{1}{|\det \mathrm{L}|} \widetilde{\Psi}(\mathrm{L}^{-T} Y),
\end{align*}
where $\mathrm{L}^{-T} \triangleq (\mathrm{L}^{-1})^T$ is the inverse transpose matrix.
We now evaluate the $L^1(\mathbb R^d)$ norm of $\Psi$ by applying the spatial linear substitution $Z = \mathrm{L}^{-T} Y$ (with volume element $dZ = |\det \mathrm{L}^{-T}| \,dY = |\det \mathrm{L}|^{-1} \,dY$, hence $dY = |\det \mathrm{L}| \,dZ$):
\begin{align}\label{eq:L1-invariance-L}
\|\Psi\|_{L^1(\mathbb R^d)} &= \int_{\mathbb R^d} |\Psi(Y)| \,dY = \int_{\mathbb R^d} \frac{1}{|\det \mathrm{L}|} \big|\widetilde{\Psi}(\mathrm{L}^{-T} Y)\big| \,dY \nonumber\\
&= \int_{\mathbb R^d} \frac{1}{|\det \mathrm{L}|} |\widetilde{\Psi}(Z)| |\det \mathrm{L}| \,dZ = \int_{\mathbb R^d} |\widetilde{\Psi}(Z)| \,dZ = \|\widetilde{\Psi}\|_{L^1(\mathbb R^d)}.
\end{align}
Notice that the Jacobian factor $|\det \mathrm{L}|$ cancels out completely. Thus, the $L^1$ norm of the convolution kernel is strictly invariant under any invertible linear map $\mathrm{L}$.

\medskip
\noindent\textbf{Step 3:} \emph{Anisotropic dyadic rescaling to unit frequency.} \\
The relocated symbol $\widetilde{\mathfrak{m}}(\zeta)$ is localized at dyadic scales $M_1, \ldots, M_r \in 2^{\mathbb Z}$. We define the diagonal dilation matrix on $\mathbb R^d = \mathbb R^2 \times \cdots \times \mathbb R^2$ by
\[
\mathsf D \triangleq \operatorname{diag}(M_1 I_2, M_2 I_2, \ldots, M_r I_2) \in \mathrm{GL}(d, \mathbb R),
\]
where $I_2$ is the $2 \times 2$ identity matrix. Its determinant is $\det \mathsf D = \prod_{\iota=1}^r M_\iota^2$.
We define the dimensionless profile $a(\eta)$ on $\mathbb R^d$ by
\[
a(\eta) \triangleq \widetilde{\mathfrak{m}}(\mathsf D \eta) = \mathfrak{m}(L^{-1} \mathsf D \eta), \qquad \eta = (\eta_1, \ldots, \eta_r) \in (\mathbb R^2)^r.
\]
By the localization hypothesis \eqref{eq:mudsp4689}, each variable $\eta_\iota = \frac{\zeta_\iota}{M_\iota}$ is localized on a fixed compact ball:
\[
\operatorname{supp} a \subset \Omega \triangleq \prod_{\iota=1}^r \big\{\eta_\iota \in \mathbb R^2 : |\eta_\iota| \leq C_0\big\} \subset \mathbb R^d.
\]
Applying the chain rule, for any multi-index $\nu = (\nu_1, \ldots, \nu_r) \in \mathbb N_0^d$ with $\nu_\iota \in \mathbb N_0^2$,
\[
\partial_\eta^\nu a(\eta) = \left( \prod_{\iota=1}^r M_\iota^{|\nu_\iota|} \right) (\partial_\zeta^\nu \widetilde{\mathfrak{m}})(\mathsf D \eta).
\]
Consequently, the derivative assumption \eqref{eq:multiplier-hypo-sharp} translates directly into the scale-invariant uniform bound
\begin{equation}\label{eq:a-derivative-uniform}
\sum_{|\nu| \leq 12} \|\partial_\eta^\nu a\|_{L^\infty(\mathbb R^d)} \leq \sum_{|\nu| \leq 12} C_\nu \mathfrak{A} \leq C_d \mathfrak{A},
\end{equation}
where $C_d \triangleq \sum_{|\nu| \leq 12} C_\nu < \infty$.

\medskip
\noindent\textbf{Step 4:} \emph{Integration by parts and uniform spatial kernel decay.} \\
Let $k \triangleq \mathcal{F}_{\mathbb R^d}^{-1} a \in \mathcal{S}(\mathbb R^d)$ denote the inverse Fourier transform of $a$:
\[
k(z) = \frac{1}{(2\pi)^d} \int_{\mathbb R^d} e^{i z \cdot \eta} a(\eta) \,d\eta, \qquad z \in \mathbb R^d.
\]
Observing that $(1 - \Delta_\eta)^N e^{i z \cdot \eta} = (1 + |z|^2)^N e^{i z \cdot \eta}$ for any integer $N \geq 1$, integrating by parts $2N$ times on $\mathbb R^d$ yields
\begin{equation}\label{eq:parts-identity}
(1 + |z|^2)^N k(z) = \frac{1}{(2\pi)^d} \int_{\mathbb R^d} (1 - \Delta_\eta)^N \big(e^{i z \cdot \eta}\big) a(\eta) \,d\eta = \frac{1}{(2\pi)^d} \int_{\mathbb R^d} e^{i z \cdot \eta} (1 - \Delta_\eta)^N a(\eta) \,d\eta.
\end{equation}
The boundary terms at infinity vanish identically because $a \in C_c^\infty(\mathbb R^d)$ has compact support in $\Omega$.
We fix $N = 6$. The differential operator $(1 - \Delta_\eta)^6 = \sum_{j=0}^6 \binom{6}{j} (-\Delta_\eta)^j$ involves derivatives of order at most $2N = 12$. It follows from \eqref{eq:a-derivative-uniform} that
\[
\big\|(1 - \Delta_\eta)^6 a\big\|_{L^\infty(\mathbb R^d)} \leq C \sum_{|\nu| \leq 12} \|\partial_\eta^\nu a\|_{L^\infty(\mathbb R^d)} \leq C \mathfrak{A}.
\]
Taking absolute values in \eqref{eq:parts-identity} with $N=6$, and using the fact that $|\operatorname{supp} a| \leq |\Omega| = (\pi C_0^2)^r < \infty$, we obtain the pointwise algebraic decay
\begin{equation}\label{eq:pointwise-decay-k}
|k(z)| \leq \frac{|\Omega|}{(2\pi)^d} \big\|(1 - \Delta_\eta)^6 a\big\|_{L^\infty} (1 + |z|^2)^{-6} \leq C(d, C_0) \mathfrak{A} (1 + |z|^2)^{-6}, \qquad \forall\, z \in \mathbb R^d.
\end{equation}
Since $d = 2r \leq 6$, the exponent satisfies $2N = 12 > d$. Integrating \eqref{eq:pointwise-decay-k} in polar coordinates on $\mathbb R^d$:
\begin{align}\label{eq:k-L1-finite}
\|k\|_{L^1(\mathbb R^d)} &= \int_{\mathbb R^d} |k(z)| \,dz \leq C(d, C_0) \mathfrak{A} \int_{\mathbb R^d} (1 + |z|^2)^{-6} \,dz \nonumber\\
&= C(d, C_0) \mathfrak{A} \omega_{d-1} \int_0^\infty (1 + \rho^2)^{-6} \rho^{d-1} \,d\rho \nonumber\\
&\leq C(d, C_0) \mathfrak{A},
\end{align}
where $\omega_{d-1}$ is the surface area of the unit sphere in $\mathbb R^d$.

\medskip
\noindent\textbf{Step 5:} \emph{Scale-invariance of the $L^1$ kernel norm.} \\
We now relate the normalized kernel $k$ to the relocated kernel $\widetilde{\Psi}$. By definition, $\widetilde{\mathfrak{m}}(\zeta) = a(\mathsf D^{-1} \zeta)$. Thus, by the scaling property of the inverse Fourier transform:
\[
\widetilde{\Psi}(Z) = \frac{1}{(2\pi)^d} \int_{\mathbb R^d} e^{i Z \cdot \zeta} a(\mathsf D^{-1}\zeta) \,d\zeta = \frac{\det \mathsf D}{(2\pi)^d} \int_{\mathbb R^d} e^{i Z \cdot (\mathsf D \eta)} a(\eta) \,d\eta = (\det \mathsf D) k(\mathsf D Z).
\]
Performing the substitution $W = \mathsf D Z$ (where $dW = (\det \mathsf D) \,dZ$):
\begin{equation}\label{eq:rescaling-L1-invariance}
\|\widetilde{\Psi}\|_{L^1(\mathbb R^d)} = \int_{\mathbb R^d} (\det \mathsf D) |k(\mathsf D Z)| \,dZ = \int_{\mathbb R^d} |k(W)| \,dW = \|k\|_{L^1(\mathbb R^d)}.
\end{equation}
Combining \eqref{eq:L1-invariance-L}, \eqref{eq:k-L1-finite}, and \eqref{eq:rescaling-L1-invariance}, we conclude that
\begin{equation}\label{eq:kernel-final-L1}
\|\Psi\|_{L^1((\mathbb R^2)^r)} = \|\widetilde{\Psi}\|_{L^1(\mathbb R^d)} = \|k\|_{L^1(\mathbb R^d)} \leq C(d, C_0) \mathfrak{A}.
\end{equation}
This establishes \eqref{eq:kernel-L1-bound-sharp}.

\medskip
\noindent\textbf{Step 6:} \emph{Minkowski and generalized H\"older estimates.} \\
With the universal $L^1$ spatial kernel bound \eqref{eq:kernel-final-L1} in hand, we now estimate the operator $T_{\mathfrak{m}}$ in Lebesgue spaces.
From the spatial integral representation \eqref{eq:spatial-integral-rep}, applying Minkowski's integral inequality yields
\begin{align}\label{eq:minkowski-step}
\|T_{\mathfrak{m}}(f_1, \ldots, f_r)\|_{L^p(\mathbb R^2)} &= \left( \int_{\mathbb R^2} \left| \int_{(\mathbb R^2)^r} \Psi(y_1, \ldots, y_r) \prod_{\iota=1}^r f_\iota(x - y_\iota) \,dy_1 \cdots dy_r \right|^p \,dx \right)^{\frac{1}{p}} \nonumber\\
&\leq \int_{(\mathbb R^2)^r} |\Psi(y_1, \ldots, y_r)| \left( \int_{\mathbb R^2} \prod_{\iota=1}^r |f_\iota(x - y_\iota)|^p \,dx \right)^{\frac{1}{p}} \,dy_1 \cdots dy_r.
\end{align}
To bound the inner $L^p$ norm, let $I_{L^{\infty}} \triangleq \{\iota \in \{1,\ldots,r\} : p_\iota = \infty\}$ and $I_{\mathrm{fin}} \triangleq \{\iota \in \{1,\ldots,r\} : p_\iota < \infty\}$.
For any $\iota \in I_{L^{\infty}}$, we pull out the bounded factors: $|f_\iota(x - y_\iota)| \leq \|f_\iota\|_{L^\infty(\mathbb R^2)}$ almost everywhere.
For the remaining indices $\iota \in I_{\mathrm{fin}}$, the reciprocal relation $\frac{1}{p} = \sum_{\iota=1}^r \frac{1}{p_\iota}$ becomes $\sum_{\iota \in I_{\mathrm{fin}}} \frac{p}{p_\iota} = 1$. Applying generalized H\"older's inequality with exponents $q_\iota \triangleq \frac{p_\iota}{p}\in [1, \infty)$ (which satisfy $\sum_{\iota \in I_{\mathrm{fin}}} \frac{1}{q_\iota} = 1$):
\begin{align*}
\left( \int_{\mathbb R^2} \prod_{\iota=1}^r |f_\iota(x - y_\iota)|^p \,dx \right)^{\frac{1}{p}} &\leq \left( \prod_{\iota \in I_{L^{\infty}}} \|f_\iota\|_{L^\infty} \right) \left( \int_{\mathbb R^2} \prod_{\iota \in I_{\mathrm{fin}}} |f_\iota(x - y_\iota)|^p \,dx \right)^{\frac{1}{p}} \\
&\leq \left( \prod_{\iota \in I_{L^{\infty}}} \|f_\iota\|_{L^\infty} \right) \prod_{\iota \in I_{\mathrm{fin}}} \left( \int_{\mathbb R^2} |f_\iota(x - y_\iota)|^{p \cdot \frac{p_\iota}{p}} \,dx \right)^{\frac{1}{p \cdot \frac{p_\iota}{p}}} \\
&= \prod_{\iota=1}^r \|f_\iota(\cdot - y_\iota)\|_{L^{p_\iota}(\mathbb R^2)} = \prod_{\iota=1}^r \|f_\iota\|_{L^{p_\iota}(\mathbb R^2)}.
\end{align*}
Substituting this back into \eqref{eq:minkowski-step} implies
\begin{align*}
\|T_{\mathfrak{m}}(f_1, \ldots, f_r)\|_{L^p(\mathbb R^2)} &\leq \left( \int_{(\mathbb R^2)^r} |\Psi(y_1, \ldots, y_r)| \,dy_1 \cdots dy_r \right) \prod_{\iota=1}^r \|f_\iota\|_{L^{p_\iota}(\mathbb R^2)} \\
&= \|\Psi\|_{L^1((\mathbb R^2)^r)} \prod_{\iota=1}^r \|f_\iota\|_{L^{p_\iota}(\mathbb R^2)} \\
&\leq C(d, C_0) \mathfrak{A} \prod_{\iota=1}^r \|f_\iota\|_{L^{p_\iota}(\mathbb R^2)},
\end{align*}
which gives the desired estimate \eqref{eq:general-multilinear-bound}.
Finally, taking $(p_1, p_2) = (p, \infty)$ for $r=2$ yields \eqref{eq:bilinear-Lp-bound}, and taking $(p_1, p_2, p_3) = (p, \infty, \infty)$ for $r=3$ yields \eqref{eq:trilinear-Lp-bound}. This completes the proof of Lemma~\ref{lem:kernel}.
\end{proof}

\vskip .1in \noindent \textbf{Acknowledgments}\,\,
Wu was partially supported by the National Science Foundation of
the United States under DMS 2104682 and DMS 2309748.
Xu was partially supported by the National Natural Science Foundation of China (No. 12571244;
12171040) and the National Key R$\&$D Program of China (No. 2020YFA0712900).
Ye was supported by the National Natural Science Foundation of China (Grant No. 12671274).

\vskip .1in
\textbf{Statements and Declarations:} The authors declare that they have no conflict of interest.

\vskip .1in
\textbf{Data availability:} Since this work is of abstract theoretical nature, no data sets are generated or analyzed. One
can obtain the relevant materials from the reference list.

\vskip .1in
\textbf{Declaration on the use of AI:}
The human authors conceived and led the research, formulated
the mathematical problem, developed the main results, and independently verified all arguments. AI assisted in refining the notation and improving the language and presentation. All mathematical arguments and proofs in the
final manuscript were checked and written by the human authors.


\vskip .3in
\begin{thebibliography}{00} \frenchspacing


\bibitem{BCD}
H. Bahouri, J.-Y. Chemin, R. Danchin,
\emph{Fourier Analysis and Nonlinear Partial Differential Equations},
Grundlehren der mathematischen Wissenschaften, vol.~343,
Springer, Heidelberg, 2011.

\bibitem{Cannon1980}
J. R. Cannon, E. DiBenedetto,
\emph{The initial value problem for the Boussinesq equations with data in $L^p$},
Lecture Notes in Math., vol.~771, Springer, Berlin, 1980, pp. 129--144.

\bibitem{CW2013}
C. Cao, J. Wu,
\emph{Global regularity for the 2D anisotropic Boussinesq equations with vertical dissipation},
Arch. Ration. Mech. Anal. \textbf{208} (2013), no.~3, 985--1004.

\bibitem{Chae}
D. Chae, \emph{Global regularity for the 2D Boussinesq equations with partial viscosity terms},
Adv. Math. \textbf{203} (2006), no.~2, 497--513.

\bibitem{ChaeNam}
D. Chae, H. Nam,
\emph{Local existence and blow-up criterion for the Boussinesq equations},
Proc. Roy. Soc. Edinburgh Sect. A \textbf{127} (1997), no.~5, 935--946.


\bibitem{Schauder}
P. Chaudru de Raynal, S. Menozzi, E. Priola, \emph{Schauder estimates for drifted fractional operators in the supercritical case},
J. Funct. Anal. \textbf{278} (2020), no.~8, 108425.


\bibitem{ChenHou2024}
J. Chen, T. Y. Hou,
\emph{Stable nearly self-similar blowup of the 2D Boussinesq and 3D Euler equations with smooth data},
Ann. of Math. (2) \textbf{199} (2024), no.~1, 213--419.

\bibitem{CW}
F. Christ,  M. Weinstein,
\emph{Dispersion of small amplitude solutions of the generalized Korteweg-de Vries equation},
J. Funct. Anal. \textbf{100} (1991), no.~1, 87--109.

\bibitem{CM}
R. Coifman, Y. Meyer,
\emph{Au-del\`a des op\'erateurs pseudo-diff\'erentiels},
Ast\'erisque, no.~57, Soci\'et\'e Math\'ematique de France, Paris, 1978.

\bibitem{CMT1994}
P. Constantin, A. Majda, E. Tabak, \emph{Formation of strong fronts in the 2-D quasigeostrophic thermal active scalar}, Nonlinearity \emph{7} (1994) 1495--1533.

\bibitem{CV}
P. Constantin, V. Vicol,
\emph{Nonlinear maximum principles for dissipative linear nonlocal operators and applications},
Geom. Funct. Anal. \textbf{22} (2012), no.~5, 1289--1321.

\bibitem{CC}
A. C\'ordoba, D. C\'ordoba,
\emph{A maximum principle applied to quasi-geostrophic equations},
Comm. Math. Phys. \textbf{249} (2004), no.~3, 511--528.

\bibitem{DP}
R. Danchin, M. Paicu,
\emph{Global well-posedness issues for the inviscid Boussinesq system with Yudovich's type data},
Comm. Math. Phys. \textbf{290} (2009), no.~1, 1--14.

\bibitem{Elgindi2021}
T. Elgindi,
\emph{Finite-time singularity formation for $C^{1,\alpha}$ solutions to the incompressible Euler equations on $\mathbb{R}^3$},
Ann. of Math. (2) \textbf{194} (2021), no.~3, 647--727.

\bibitem{El1jo20}
T. Elgindi, I. Jeong, \emph{Finite-time singularity formation for strong solutions to the Boussinesq system}, Ann. PDE \textbf{6} (2020), no. 1, Paper No. 5, 50 pp.

\bibitem{EShu1994}
W. E, C. Shu,
\emph{Small-scale structures in Boussinesq convection},
Phys. Fluids \textbf{6} (1994), no.~1, 49--58.


\bibitem{Fujiwara}
K. Fujiwara,
\emph{Remark on the chain rule of fractional derivative in the Sobolev framework},
Math. Inequal. Appl. \textbf{24} (2021), 1113--1124.


\bibitem{HK}
T. Hmidi, S. Keraani,
\emph{On the global well-posedness of the Boussinesq system with zero viscosity},
Indiana Univ. Math. J. \textbf{58} (2009), no.~4, 1591--1618.

\bibitem{HKRNS}
T. Hmidi, S. Keraani, and F. Rousset,
\emph{Global well-posedness for a Boussinesq--Navier--Stokes system with critical dissipation},
J. Differential Equations \textbf{249} (2010), no.~9, 2147--2174.

\bibitem{HKREuler}
T. Hmidi, S. Keraani, F. Rousset,
\emph{Global well-posedness for Euler--Boussinesq system with critical dissipation},
Comm. Partial Differential Equations \textbf{36} (2011), no.~3, 420--445.

\bibitem{HZ}
T. Hmidi, M. Zerguine,
\emph{On the global well-posedness of the Euler--Boussinesq system with fractional dissipation},
Phys. D \textbf{239} (2010), 1387--1401.

\bibitem{HL}
T. Hou, C. Li,
\emph{Global well-posedness of the viscous Boussinesq equations},
Discrete Contin. Dyn. Syst. \textbf{12} (2005), no.~1, 1--12.

\bibitem{JMWZ}
Q. Jiu, C. Miao, J. Wu, Z. Zhang,
\emph{The 2D incompressible Boussinesq equations with general critical dissipation},
SIAM J. Math. Anal. \textbf{46} (2014), no.~5, 3426--3454.

\bibitem{LLT2013}
A. Larios, E. Lunasin, E. S. Titi,
\emph{Global well-posedness for the 2D Boussinesq system with anisotropic viscosity and without heat diffusion},
J. Differential Equations \textbf{255} (2013), no.~9, 2636--2654.


\bibitem{LuoHou2014}
G. Luo, T. Hou,
\emph{Toward the finite-time blowup of the 3D axisymmetric Euler equations: a numerical investigation},
Multiscale Model. Simul. \textbf{12} (2014), no.~4, 1722--1776.


\bibitem{Majda2003}
A. J. Majda,
\emph{Introduction to PDEs and Waves for the Atmosphere and Ocean},
Courant Lecture Notes in Mathematics, vol.~9, AMS, Providence, RI, 2003.

\bibitem{MB}
A. J. Majda, A. L. Bertozzi,
\emph{Vorticity and Incompressible Flow},
Cambridge Texts in Applied Mathematics, vol.~27,
Cambridge University Press, Cambridge, 2002.

\bibitem{MX}
C. Miao, L. Xue,
\emph{On the global well-posedness of a class of Boussinesq--Navier--Stokes systems},
NoDEA Nonlinear Differential Equations Appl. \textbf{18} (2011), 707--735.

\bibitem{Pedlosky}
J. Pedlosky,
\emph{Geophysical Fluid Dynamics}, 2nd ed.,
Springer, New York, 1987.

\bibitem{SZK1993}
M. Shlesinger, G. Zaslavsky, J. Klafter,  \emph{Strange kinetics}, Nature, \textbf{363}(6424), (1993), 31--37.



\bibitem{Silvestre}
L. Silvestre,
\emph{H\"older estimates for advection fractional-diffusion equations},
Ann. Sc. Norm. Super. Pisa Cl. Sci. (5) \textbf{11} (2012), no.~4, 843--855.


\bibitem{Sil20121}
L. Silvestre, \emph{On the differentiability of the solution to an equation with drift and fractional diffusion}, Indiana Univ. Math. J. \textbf{61} (2012), 557--584.


\bibitem{SW}
A. Stefanov, J. Wu,
\emph{A global regularity result for the 2D Boussinesq equations with critical dissipation},
J. Anal. Math. \textbf{137} (2019), 269--290.

\bibitem{SWXY}
A. Stefanov, J. Wu, X. Xu, Z. Ye,
\emph{Global regularity results of the 2D fractional Boussinesq equations},
Math. Ann. \textbf{391} (2025), 5965--6012.

\bibitem{SWXY2026}
A. Stefanov, J. Wu, X. Xu, Z. Ye,
\emph{Global regularity of the 2D fractional Boussinesq equations with subcritical dissipation},
preprint, arXiv:2606.03680, 2026.

\bibitem{WXXY}
J. Wu, X. Xu, L. Xue, Z. Ye,
\emph{Regularity results for the 2D Boussinesq equations with critical or supercritical dissipation},
Commun. Math. Sci. \textbf{14} (2016), no.~7, 1963--1997.

\bibitem{XueYepjm}
L. Xue, Z. Ye, \emph{On the differentiability issue of the drift-diffusion equation with nonlocal L$\rm\acute{e}$vy-type diffusion}, Pacific J. Math. \textbf{293} (2018), 471--510.

\bibitem{Ye2017}
Z. Ye,
\emph{Global smooth solution to the 2D Boussinesq equations with fractional dissipation},
Math. Methods Appl. Sci. \textbf{40} (2017), no.~12, 4595--4612.

\bibitem{Ye2020}
Z. Ye,
\emph{An alternative approach to global regularity for the 2D Euler--Boussinesq equations with critical dissipation},
Nonlinear Anal. \textbf{190} (2020), 111591.

\bibitem{Yejde24}
Z. Ye, \emph{Global well-posedness for the 2D Euler-Boussinesq-B\'enard equations with critical dissipation}, J. Differential Equations \textbf{392} (2024), 209--254.

\bibitem{YX}
Z. Ye, X. Xu,
 \emph{Global well-posedness of the 2D Boussinesq equations with fractional Laplacian dissipation},
J. Differential Equations \textbf{260} (2016), no.~8, 6716--6744.

\bibitem{ZLSWYZ}
D. Zhou, Z. Li, H. Shang, J. Wu, B. Yuan, J. Zhao,
\emph{Global well-posedness for the 2D fractional Boussinesq equations in the subcritical case},
Pacific J. Math. \textbf{298} (2019), no.~1, 233--255.


\end{thebibliography}
\end{document}